\documentclass[11pt]{article}
\usepackage{ae,lmodern} 
\usepackage[latin1]{inputenc}
\usepackage[T1]{fontenc}
\usepackage{amsmath}
\usepackage{amsthm}
\usepackage{mathrsfs}
\usepackage{float}
\usepackage{capt-of}
\usepackage{booktabs}
\usepackage{varwidth}

\usepackage[font=small]{caption}
\usepackage{mathtools}
\usepackage{booktabs}
\usepackage{amscd}
\usepackage{amssymb}
\usepackage{amsfonts}
\usepackage[dvipsnames]{xcolor}
\usepackage{mathtools}
\usepackage{cancel}
\usepackage{stmaryrd}
\SetSymbolFont{stmry}{bold}{U}{stmry}{m}{n}
\usepackage{geometry}
\usepackage{cite}
\usepackage{enumitem}
\usepackage[textwidth=2cm, textsize=scriptsize,color=green!40]{todonotes}
\usepackage[percent]{overpic}
\usepackage{grffile}%
\usepackage{tikz,pgfplots}
\usepackage[export]{adjustbox} 
\usepackage{algorithm}
\usepackage{algorithmicx}
\usepackage{algpseudocode}
\usepackage{printlen} 
\usepackage{upgreek}
\usepackage{subcaption}
\usepackage[percent]{overpic}
\usepackage{pict2e}
\usepackage{hyperref}
\usepackage{textcomp}
\usepackage[OT1]{fontenc}
\usepackage[norefs, nocites]{refcheck}
\newcommand{\silentref}[1]{\sbox0{\ref{#1}}}
\hypersetup{
	colorlinks   = true, 
	urlcolor     = blue, 
	linkcolor    = blue, 
	citecolor   = red 
}

\pgfplotsset{compat=1.16}

\setenumerate[2]{label=(\alph*)}
\setenumerate[1]{label=\textup{(\roman*)}}
\setenumerate[3]{label=(\arabic*)}

\newcommand{\normW}[1]{{\vert\kern-0.25ex\vert\kern-0.25ex\vert #1 
		\vert\kern-0.25ex\vert\kern-0.25ex\vert}}
\newcommand{\normWbig}[1]{{\Big\vert\kern-0.25ex\Big\vert\kern-0.25ex\Big\vert #1 
		\Big\vert\kern-0.25ex\Big\vert\kern-0.25ex\Big\vert}}

\makeatletter
\newcommand\notsotiny{\@setfontsize\notsotiny{6.8}{7.1828}}
\makeatother

\newcommand{\C}{\mathbb{C}}
\newcommand{\N}{\mathbb{N}}

\newcommand{\R}{\mathbb{R}}
\newcommand{\K}{\mathbb{K}}

\newcommand{\T}{\mathbb{T}}

\newcommand{\Z}{\mathbb{Z}}
\newcommand{\E}{\mathbb{E}}

\newcommand{\boC}{\mathcal{C}}
\newcommand{\boD}{\mathcal{D}}
\newcommand{\boE}{\mathcal{E}}

\newcommand{\boG}{\mathcal{G}}
\newcommand{\boH}{\mathcal{H}}
\newcommand{\boI}{\mathcal{I}}
\newcommand{\boJ}{\mathcal{J}}

\newcommand{\boO}{\mathcal{O}}

\newcommand{\boQ}{\mathcal{Q}}
\newcommand{\boR}{\mathcal{R}}
\newcommand{\boS}{\mathcal{S}}
\newcommand{\boT}{\mathcal{T}}

\newcommand{\gs}{\mathfrak{s}}

\newcommand{\ka}{\upkappa}

\newcommand{\Nt}{N}
\newcommand{\dt}{\uptau}
\newcommand{\dx}{h}
\newcommand{\J}{J}
\newcommand{\tnmed}{t_{n+\frac{1}{2}}}

\newcommand{\Dxx}{D_{xx}^\dx}
\newcommand{\Dx}{D_x^{\dx}}
\newcommand{\Dxp}{D_x^{\dx^+}}
\newcommand{\dbra}[2]{\llbracket #1, #2 \rrbracket}

\newcommand{\ptl}{{\partial}}

\newcommand{\Dns}{\mathcal{Q}}
\newcommand{\dns}{q}

\newcommand{\Rh}{\boG_\dx}
\newcommand{\Fh}{\boH_\dx}

\providecommand{\abs}[1]{|#1 |}
\providecommand{\abssqr}[1]{|#1 |^2}
\providecommand{\norm}[1]{\lVert#1 \rVert}
\providecommand{\infnorm}[1]{\lVert#1 \rVert_\infty}

\providecommand{\CHnorm}[3]{\lVert#1 \rVert_{C\IfValueT{#3}{^{#3}}(H^#2)}}
\providecommand{\CLnorm}[2]{\lVert#1 \rVert_{C\IfValueT{#2}{^{#2}}(L^\infty)}}

\providecommand{\ipCp}[2]{\ipp{#1, #2}_{\dx^+}}
\providecommand{\ipC}[2]{\ipp{#1, #2}_{\dx}}
\providecommand{\ipCm}[2]{\ipp{#1, #2}_{\dx^-}}

\providecommand{\ipRp}[2]{\langle #1, #2 \rangle_{\dx^+}}
\providecommand{\ipR}[2]{\langle #1, #2\rangle_{\dx}}
\providecommand{\ipRm}[2]{\langle #1, #2 \rangle_{\dx^-}}

\providecommand{\normpm}[2]{\norm{#1}_{#2, \dx^-}}
\providecommand{\normp}[2]{\norm{#1}_{#2, \dx}}
\providecommand{\normpp}[2]{\norm{#1}_{#2, \dx^+}}

\providecommand{\ipp}[1]{\langle\negthickspace\langle #1 \rangle\negthickspace\rangle}       
\providecommand{\ip}[2]{\langle #1, #2 \rangle}       

\renewcommand{\Re}{\operatorname{Re}}
\renewcommand{\Im}{\operatorname{Im}}

\newcommand{\atanh}{\operatorname{atanh}}
\newcommand{\atan}{\operatorname{atan}}

\renewcommand{\emph}[1]{\textit{#1}}

\newcommand{\conj}[1]{\overline{#1}}

\theoremstyle{plain}
\newtheorem{theorem}{Theorem}[section]
\newtheorem{proposition}[theorem]{Proposition}
\newtheorem{lemma}[theorem]{Lemma}
\newtheorem{corollary}[theorem]{Corollary}

\theoremstyle{definition}

\newtheorem{remark}[theorem]{Remark}

\theoremstyle{remark}

\newtheorem*{merci}{Acknowledgments}

\numberwithin{equation}{section}

\everymath=\expandafter{\the\everymath\displaystyle}

\begin{document}

\title{Numerical analysis of a conservative Crank--Nicolson scheme
for quasilinear Gross--Pitaevskii equations}    
	\author{
		\renewcommand{\thefootnote}{\arabic{footnote}}
		Sebasti\'an Tapia--Mandiola\footnotemark[1]}
	\footnotetext[1]{
		Univ.\ Lille, CNRS, Inria, UMR 8524 - Laboratoire Paul Painlev\'e, Inria, F-59000 Lille, France.\\
		E-mail: {\tt sebastian.tapia@univ-lille.fr}}
        \date{}
	\maketitle

\begin{abstract}
We propose and analyze a fully discrete Crank--Nicolson finite-difference scheme for one-dimensional quasilinear Hamiltonian Gross--Pitaevskii equations with homogeneous Neumann boundary conditions. We establish the existence and uniqueness of numerical solutions and an $H^1$-stability estimate.

For sufficiently regular exact solutions, we obtain second-order convergence in both time and space in the discrete $H^1$-norm. 
The analysis covers the focusing and defocusing cubic nonlinearities and both signs of the quasilinear coefficient $\ka$, under a uniform discrete ellipticity condition when $\ka<0$.

The main difficulty in the error analysis comes from the derivatives introduced by the quasilinear interaction.
To overcome this problem, we introduce a discrete action combining the discrete energy and mass,
and establish a coercivity inequality in the discrete $H^1$-norm that allows us to control the numerical error. Together with consistency estimates and a discrete Gronwall inequality, this yields second-order convergence. Numerical experiments with bright solitons corroborate the predicted convergence rates.

\end{abstract}


\section{Introduction}
For $T>0, L>0$ and $\mathcal{I} = [0, T], \mathcal{J} = [-L, L]$, we consider the one-dimensional initial-boundary value problem for the quasilinear Gross--Pitaevskii equation (QGP) given by
\begin{equation}\label{QGP}\tag{QGP}
\begin{aligned}
i\partial_t \psi = \partial_{xx} \psi + \mathfrak{s}\psi|\psi|^2 +  \upkappa \psi\partial_{xx}|\psi|^2, & \quad (x, t) \in \mathcal{J} \times \mathcal{I}\\
\ptl_x\psi(t, -L) = \ptl_x\psi(t, L) = 0, & \quad t \in \mathcal{I} \\
\psi(0, x) = \psi_0(x), & \quad x \in \mathcal{J},
\end{aligned}
\end{equation}
where $ \psi: [-L, L]\times [0, T] \to \C$ 
denotes the complex-valued unknown wave function, and $\psi_0$ the initial condition. The parameter $ \mathfrak{s} \in \mathbb{R}\setminus\{0\} $ distinguishes between the focusing ($ \mathfrak{s} > 0 $) and defocusing ($ \mathfrak{s} < 0 $) regimes of the cubic semilinear response, while $ \upkappa \in \mathbb{R} $ measures the strength of the quasilinear interaction. By time and space rescaling, we can assume that $\gs$ equals $1$ or $-1$.

In what follows, we introduce the notation $q=|\psi|^2$ for the density of the wave function. Equation \eqref{QGP} has a Hamiltonian structure, and it preserves in time the energy and the mass given by
\begin{equation}
\label{eq:energie}
   E(t) := \frac{1}{2}\int_{\mathcal{J}} |\partial_x \psi(x,t)|^2 dx  -\frac{\gs }{4}    \int_{\mathcal{J}}q(x, t)^2 dx+ \frac{\upkappa}{4}\int_{\mathcal{J}} \left( \partial_x q(x, t)\right)^2 dx = E(0),
\end{equation}
and 
\begin{equation}\label{eq:mass}
    N(t) := \int_{\mathcal{J}} q(x,t) dx = N(0),
\end{equation}
respectively.

This equation has been introduced as a more general model of the cubic nonlinear Schr\"odinger equation (NLS) ($\ka=0$) in several physical contexts when complex effects are taken into account. For instance, we can mention \cite{rutledge} for superfluidity in thin helium films, \cite{litvak} and \cite{shukla} for the dynamics of Langmuir oscillations in isotropic plasmas, and \cite{de_laire_exotic_2026,krolikowski_solitons_2000,horikisWnonloc} for the evolution of optical beams in weakly nonlocal Kerr-like media. Lastly, it arises in modeling Bose--Einstein condensates with finite-range corrections as explained in \cite{FuWangGao2003,CollinMassignanPethick2007, SgarlataMazzarellaSalasnich2015, LiDuanWangYang2026}. See Subsection~1.1 in the companion paper \cite{de_laire_simul_2026} for further physical explanations.

We identify various regimes based on the signs of the coefficients $\gs$ and $\upkappa$. The coefficient $\gs$ determines the character of the local cubic interaction: \(\mathfrak{s}>0\) corresponds to an attractive case, whereas \(\mathfrak{s}<0\) corresponds to a repulsive case. For example, in Bose--Einstein condensates, this term models attractive interatomic interactions in the focusing case and repulsive interactions in the defocusing case. For the quasilinear term, the coefficient $\ka$ plays a similar role, as observed from the energy \eqref{eq:energie}. Indeed, from the energy, we see that when $\ka>0$, the quasilinear term penalizes spatial variations of the gradient of $q$ and is therefore stabilizing, while $\ka<0$ favors concentrations. 

Equation \eqref{QGP} is locally well-posed in $H^p(\mathcal J)$, with 
 $p>13$, as established in Corollary~\ref{lucario}. 
Precisely, defining the loss of ellipticity function:
	\begin{equation*}
		\boD_L(\psi)
:=\inf_{x\in[-L,L]}
\bigl(1+2{\upkappa}|\psi(x)|^2\bigr),
\quad \text{for }  \psi\in C([-L,L];\C),
	\end{equation*}
we can show that for any $\psi_0\in H^p(\mathcal J)$, 
satisfying the compatibility conditions in 
\eqref{comp_cond_init} and $ \boD_{L}(\psi_0)>0,$
 \eqref{QGP} admits a unique maximal solution $\psi$
on $[0,T^*)$, with $T^*\in(0,\infty]$,
such that, for every $0<T<T^*$,
$
 \psi\in C([0,T];H^p(\mathcal J))
 \cap C^1([0,T];H^{p-2}(\mathcal J)),
$
and 
\begin{equation}
\label{the:condition}
\boD_{L}(\psi(t))>0, \quad\text{  for all }t\in[0,T^*).
\end{equation}
 
 Remark that the function  $\boD_{L}$ is strictly positive for $\ka\geq 0$, 
 so that, in this case,  \eqref{the:condition} is  trivially satisfied.
 Hence, the loss of ellipticity can only occur when  $\ka<0$, in which case $\boD_{L}(\psi_0)>0$ and \eqref{the:condition} are meaningful.
 
 The proof of this local well-posedness result is a consequence of the theorems
 by Feola and Iandoli \cite{feola2022local}, 
 for quasilinear Hamiltonian Schr\"odinger equations on the torus.


\subsection{Main result}

Let $T>0$ and $\Nt\in \N\backslash \{0\}$ be such that $\Nt+1$ corresponds to the total number of time steps, and  $\dt = T / \Nt$ as the size of the time step. The discrete time instants are given by $t_n = n \dt$ for $n \in \dbra{0}{\Nt}$. We also use the notation $\tnmed$ for the intermediate state $t_n+\dt/2$. In space, we divide the interval $\mathcal{J}$ into $\J +1$ grid points with $\J\in \N\backslash \{0\}$. The spatial stepsize is $\dx = 2L / \J$, and the grid nodes correspond to $x_j = j\dx - L$ for $j \in \dbra{0}{\J}$.

Let $\psi$ be a solution of \eqref{QGP} in  $C^4(\boI; H^5(\boJ))$, with initial condition 
$\psi_0$. We denote by $\psi_j^n = \psi(x_j, t_n)$ and $\psi^n = (\psi_j^n)_{j\in \dbra{0}{\J}}$ in its vector representation; and the discrete representation of the density by $q^n = (q_j^n)_{j\in \dbra{0}{\J}}$ with $q_j^n = |\psi_j^n|^2$.

We propose the following fully-discrete Crank--Nicolson (CN) scheme to approximate  $\psi^n$ by a discrete complex vector $\Psi^n : = (\Psi_j^n)_{j\in \dbra{0}{\J}}\in \C^{\J+1}$, and the density $q^n$ by $\boQ^n = |\psi^n|^2$:
\begin{equation}
\label{schme}
   \frac{i\left(\Psi^{n+1} -\Psi^{n}\right)}{\dt} = \Dxx\Psi^{n+\frac{1}{2}} + \gs \Psi^{n+\frac{1}{2}} \odot \boQ^{n+\frac{1}{2}} + \ka \Psi^{n+\frac{1}{2}} \odot \Dxx \boQ^{n+\frac{1}{2}},
\end{equation}
where $\odot$ denotes the pointwise product of two vectors, $\Psi^0 = (\psi_0(x_j))_{j\in \dbra{0}{\J}}$, $\Psi^{n+\frac{1}{2}}=(\Psi^{n+1}+\Psi^n)/2$, $\boQ^{n+\frac{1}{2}}= (\boQ^n + \boQ^{n+1})/2$, and $\Dxx$ is the centered second-order finite-difference operator, defined in \eqref{eq:spc_disc_ops}, with homogeneous Neumann boundary conditions implemented as $(\Dxx U)_0 = (2U_1 - 2U_0)/\dx^2$ and $(\Dxx U)_\J = (2U_{\J-1} - 2U_\J)/\dx^2$. This method is implicit in time, and it requires solving a nonlinear algebraic system in each time step.

The main contribution of this work is the following theorem on the estimate of the error $e^n=\Psi^n - \psi^n$ in the discrete $H^1$-norm (defined in \eqref{eq:norms_H1}), which establishes a second-order convergence in both space and time for the CN scheme:
\begin{theorem}[Error convergence]\label{thm:error_convergence}
    Let $\gs \in \{-1, 1\}$ and $\ka\in \R$. 
Assume that there is a solution   $\psi\in C^4(\boI; H^5(\boJ))$  of \eqref{QGP},
with vector representation $\psi^n = (\psi_j^n)_{j\in \dbra{0}{\J}}$.
Consider  $(\Psi^n)_{n\in \dbra{0}{\Nt}}$ a solution of the scheme \eqref{schme}, 
with $\Psi^0 = \psi^0$, satisfying 
    the discrete ellipticity condition 
\begin{equation}\label{eq:disc_ellipticity}
    \boD_n:=1+2\ka\max_{0\leq j\leq \J} \abs{\Psi_j^n}^2\geq \delta>0, \qquad \text{ for all }\quad 0\leq n\leq N,
\end{equation}
for some $\delta >0$.
There exist 
$C>0$,  such that for every $\dx, \dt>0$ sufficiently small,
we have the error estimate
\begin{equation}\label{eq:global-stability}
    \max_{0\leq n\leq N}\norm{\Psi^n - \psi^n}_{H^1_\dx} \leq C(\dt^2+\dx^2).
\end{equation}
Notice that condition \eqref{eq:disc_ellipticity} is trivially satisfied if $\ka\geq 0$.
\end{theorem}
Let us make two remarks concerning the existence of the $\psi$ and $\Psi^n$ in this theorem.
As explained above, for any $\psi_0 \in H^p(\boJ)$, $p>13$, satisfying 
 $ \boD_{L}(\psi_0)>0$, Corollary~\ref{lucario} provides the existence of a unique maximal 
solution  $\psi$ on $[0,T^*)$, satisfying  \eqref{the:condition}.
Moreover, a bootstrap argument yields $\psi \in C^4([0, T]; H^5(\boJ))$ for $T<T^*$.

On the other hand, the existence of the discrete solution $\Psi^n$ is established in Subsection~\ref{subec:discrete_existence} in 
Proposition~\ref{prop:discrete-existence}, 
in which we rely on a fixed-point argument as in \cite{akrivis_1993}, while the uniqueness of the scheme is proved in Proposition~\ref{prop:cn-uniqueness}, 
under a condition between the steps $\dt$ and $\dx$. Notice that \eqref{eq:disc_ellipticity} is the 
 discrete version of the ellipticity the condition
\eqref{the:condition}. In constrat to the continuous case, condition \eqref{eq:disc_ellipticity} is not required for the solution to exist in the discrete setting, but, in the meaningful case $\ka<0$, it will be crucial to guarantee the stability of the scheme proved in Proposition~\ref{prop:stability}.

The key elements to prove our theorem are the second-order consistency of the scheme, established in Proposition~\eqref{prop:second_order} in Section~\ref{sec:consistency}, and the discrete $H^1$-stability estimate, established in Proposition~\eqref{prop:stability} in Subsection~\ref{sec:stability}. However, since the quasilinear term induces extra derivatives in the analysis, the convergence is not directly obtainable from both consistency and stability, nor by using energy estimates as in NLS equations; see, e.g.\ \cite{akrivis_1993, sanz_serna_1984}. To overcome this difficulty, we introduce the \textit{discrete action} operator defined by
\begin{equation}\label{eq:action_op:intro}
    \boS_{\dx,\mu}(\Psi) = E_\dx(\Psi) + \frac{\mu}{2} N_\dx(\Psi),
\end{equation}
where $E_\dx$ and $N_\dx$ are the discrete energy and mass, defined in \eqref{eq:full_energy} and \eqref{eq:full_mass}, respectively; and $\mu>0$ is a constant. This type of operator has been widely used in the continuous setting over $\R$ to analyze soliton stability. Indeed, for a soliton $u$ and $\omega>0$, a Taylor expansion of the continuous action of the soliton $\boS(u) := E(u) + \omega N(u)$ around $u$ with perturbation $v$ gives
\begin{equation*}
    \boS(u+v) = \boS(u) + D\boS(u)[v] + D^2\boS(u)[v, v] + \mathcal{O}(\norm{v}_{H^1}^3), 
\end{equation*}
where $D^2\boS(u)[v, v]$ is the Hessian of the action $\boS$ at $\phi$ in the direction $v$. It is easy to check that 
$D\boS(u)[v]=0$ for all perturbation $v$ since it corresponds to the Euler--Lagrange equation for the soliton $u$ given by:
\begin{equation}\label{eq:soliton_eq}
    u'' + u^3 - \omega u + \ka u (u^2)'' = 0.
\end{equation}
The nontrivial positive solution to equation \eqref{eq:soliton_eq} is explicitly given in \eqref{bright} and \eqref{eq:implibright}, where $u = u_{\omega, \ka}$. By continuity, the Hessian satisfies $|D^2\boS(u)[v, v]| \leq C\norm{v}_{H^1}^2$, and the main challenge is to impose conditions on the soliton to have the coercivity of $D^2\boS$; see e.g., \cite{quiniou_stability_2025,grillakisshatah} and the references therein for more details.

Following this idea, we can write the discrete solution $\Psi^n$ as a perturbation of the exact solution $\psi^n$ by the error $e^n$, i.e., $\Psi^n = \psi^n + e^n$. Then, a Taylor expansion of the discrete action $\boS_\dx(\Psi^n)$ around $\psi^n$ gives
\begin{equation}\label{eq:taylor_action}
    \boS_{\dx, \mu}(\Psi^n) = \boS_{\dx, \mu}(\psi^n) + D\boS_{\dx, \mu}(\psi^n)[e^n] + \boH_{\dx, \mu}(\psi^n, e^n).
\end{equation}
Here, the term $\boH_{\dx, \mu}$ corresponds to the remainder defined as 
\begin{equation}\label{eq:F_remainder}
\boH_{\dx, \mu}(U, V) = \int_0^1 (1 - \theta) D^2\boS_{\dx, \mu}(U+\theta V)[V, V]d\theta,    
\end{equation}
and $D\boS_{\dx, \mu}(U)[V]$ and $D^2\boS_{\dx, \mu}(U)[V, V]$ are the 
derivatives of $\boS_{\dx, \mu}$ with respect to $U$ in the direction $V$. In the particular case of $U = \psi^n$ and $V = e^n$, we write $\boH_{\dx, \mu}^n = \boH_{\dx, \mu}(\psi^n, e^n)$. Equivalently, we can write $\boH_{\dx, \mu}^n$ as
\begin{equation*}
    \boH_{\dx, \mu}^n = \boS_{\dx, \mu}(\Psi^n) - \boS_{\dx, \mu}(\psi^n) - D\boS_{\dx, \mu}(\psi^n)[e^n].
\end{equation*}
The essential point of $\boH_{\dx, \mu}^n$ and why it is important here is its equivalence with the discrete $H^1$-norm of the error $e^n$ established in the following proposition:
\begin{proposition}\label{prop:equivalence}
    Let $(\Psi^n)_{n\in \dbra{0}{\Nt}}$ a solution of \eqref{scheme}. If $\ka<0$ assume that $\boD_n\geq \delta>0$ for all $n\in \dbra{0}{\Nt}$. Then, 
    there exists $\mu^*>0$ such that for all $\mu \geq \mu^*$, there exists constants $c, C>0$ such that for all $\dx, \dt>0$ and $n\in \dbra{0}{\Nt}$,
        \begin{equation}\label{ineq:coercivity_H}
            c\norm{e^n}_{H_h^1}^2 \leq \boH_{\dx, \mu}^n \leq C\norm{e^n}_{H_h^1}^2.
        \end{equation}
\end{proposition}
The estimate in \eqref{ineq:coercivity_H} allows us to control the numerical error and conclude the convergence of the scheme.
\begin{remark}
The proof of Theorem~\ref{thm:error_convergence} using our action strategy is valid for any $\ka\in \R$. In particular, for $\ka=0$, it gives a novel proof of convergence for the CN scheme applied to the cubic NLS equation.
\end{remark}

\paragraph{Organization of the paper.}
The paper organizes as follows: In Section~\ref{juanito}, we introduce a series of auxiliary identities and inequalities. In Section~\ref{pablito}, we formulate the Crank--Nicolson scheme, and derive the results on the stability, existence and uniqueness of the approximation. Section~\ref{sec:consistency} contains the estimations of the consistency error. In Section~\ref{sec:convergence}, we give the proof of the convergence result in Theorem~\ref{thm:error_convergence}, and finally in Section~\ref{sec:numerics}, we perform numerical experiments to validate the convergence result by comparing theoretical bright solitons against our scheme. More experiments can be found in our companion paper \cite{de_laire_simul_2026}.

\subsection{State of the art}
In the literature, CN schemes have been widely used for NLS equations due to their second-order accuracy in time and mass- and energy-conservation properties. For instance, we can mention the works on nonlinear Schr\"odinger-type plasma models in \cite{bruneau_numerical_1999}, and on the stability of solitary waves for a quasilinear Schr\"odinger equation in \cite{Bahhi2026}.

For the cubic NLS ($\ka=0$), in \cite{sanz_serna_1984} and in \cite{akrivis_1993} the authors established a rigorous convergence analysis for the CN scheme for NLS equations. Other conservative finite-difference schemes were also studied in \cite{delfour_fortin_payre_1981} and in \cite{sanz_serna_verwer_1986}, while in \cite{akrivis_dougalis_karakashian_1991} they analyzed Galerkin approximations of second-order temporal accuracy, then extended in \cite{henning_peterseim_2017} to general nonlinearities, and in \cite{bao_cai_2013, HenningMalqvist2017} to the Gross--Pitaevskii equation with angular-momentum rotation. In the Galerkin setting, 
$L^2$-error estimates of optimal order in space were derived in \cite{KarakashianAkrivisDougalis1993} for implicit Runge--Kutta time discretizations, and in \cite{KarakashianMakridakis1998DG, KarakashianMakridakis1999CG} for space-time finite element methods. We refer the reader to \cite{antoine_computational_2013} for a review of numerical methods for NLS equations.

The convergence of energy-preserving schemes for Schr\"odinger equations usually relies on energy estimates as in \cite{akrivis_1993, henning_warnegard_2021, Besse2004, BesseDescombesDujardinLacroixViolet2021}. In the case of the relaxation method introduced in \cite{Besse1998, Besse2004}, the estimates are performed on a relative energy preserved by the scheme. In this line, in \cite{Zouraris2023} they derived a optimal-order error estimates in a discrete $L^\infty(H^1)$-norm for the relaxation scheme in NLS, and in \cite{LiuYiYin2026} the authors proved optimal error estimates for the same method applied to the Schr\"odinger--Poisson equation. Similarly, in \cite{PoulainSchratz2022} they established the convergence of the scalar auxiliary variable (SAV) method for the NLS equation, which was also used in \cite{AntoineShenTang2021} for the dynamics of Gross--Pitaevskii equations. 
In the case of splitting methods, one can apply Duhamel expansions and commutator estimates as in \cite{lubich_splitting_2008}. More recently, for the derivative NLS, in \cite{Vougiantzis_2026} they obtained a second-order convergence in the discrete $L^2$ norm by analyzing an auxiliary modified CN scheme.

For the quasilinear setting $\ka\neq0$, in \cite{lu_strang_2015} a second-order convergence in time is proved of a Strang splitting method for sufficiently small and regular initial data, and in \cite{krolikowski_solitons_2000} a split-step Fourier method is used to investigate bright and dark solitons. The numerical analysis of our CN scheme for \eqref{QGP} has not yet been addressed.

Our discrete action strategy is also called \textit{relative energy method} in the context of fluid dynamics; see e.g.\ \cite{Brenier2000, Puel2002Incompressible, Puel2002Magnetic, JungelWang2003, AlazardCarles2009, KimMoon2022, KimMoon2023, KimMoon2024, GiesselmannLattanzioTzavaras2017} where this method has been used in the continous level. In fact, the remainder $\boH_{\dx, \mu}^n$ in \eqref{eq:taylor_action} is precisely the discrete relative energy of the discrete action $\boS_{\dx, \mu}$ at $\Psi^n$ with respect to $\psi^n$, and the coercivity estimate \eqref{ineq:coercivity_H} plays the role of the convexity of the energy in the relative energy method. In the discrete level, it has been used to prove the convergence of numerical schemes to, for instance, compressible Navier--Stokes equations in \cite{GallouetHerbinMalteseNovotny2016, FeireislHosekMalteseNovotny2017, MizerovaShe2020, GallouetMerigotNatale2022}, for Cahn--Hilliard-type systems in \cite{BrunkEggerHabrichLukacova2023CH, BrunkEggerHabrichLukacova2023CHNS, BrunkEggerHabrich2024}, and for the analysis of gas transport in pipe networks in \cite{EggerGiesselmann2023, EggerGiesselmannKunkelPhilippi2023}.

However, to the best of our knowledge, our discrete action strategy has not previously been used to prove convergence of numerical schemes for Schr\"odinger equations. The closest approach is that in \cite{bambusi_faou_grebert_2013}, who used a discrete variational coercivity to establish orbital stability near NLS ground states.


\subsection{Notation} Throughout the paper, $c$ and $C$ denote generic positive constants. We denote by $H^p(\R)$ the Sobolev space of complex-valued functions on $\R$ with weak derivatives up to order $p\in [1, \infty)$ in $L^2(\R)$. For $m\geq 0,\, p\in [1, \infty)$ and $I\subset \R$ an interval, we denote $C^m(I, H^p(\R))$ the space of $H^p(\R)$-valued functions on $I$ of class $C^m$. For notational simplicity, this space will usually be written $C^m(H^p)$. We equip this space with the norm
\begin{equation*}
    \norm{\psi}_{C^m(H^p)} = \sum_{j=0}^{m}\sup_{t\in I}\norm{\ptl_t^{j}\psi(t)}_{H^p}.
\end{equation*}

For $z_1, z_2 \in \C$, we define the inner product in $\C$ by $\ip{z_1}{z_2}=\Re(z_1)\Re(z_2)+\Im(z_1)\Im(z_2) = \Re(z_1\conj{z_2}) = \Re(\conj{z_1}z_2)$, where $\conj{z_1}$ denotes the complex conjugate of $z_1$.

For $N_1, N_2 \in \N$, $N_1<N_2$, we also introduce the notation $\dbra{N_1}{N_2}$, which represents the set of integers $\{N_1, N_1+1, ..., N_2\}$..

For given $\J\in \N\backslash\{0\}$ and $u = (u_0, u_1, ..., u_{\J-1})^T, v=(v_0, v_1, ..., v_{\J-1})^T\in \C^\J$, we define the pointwise product $\odot$ between $u$ and $v$ as $u\odot v = (u_0 v_0, u_1 v_1, ..., u_{\J-1} v_{\J-1})^T$. Here, $(\cdot)^T$ denotes the transpose operator.

For any function $F: \K_1 \to \K_2$, for $\K_1, \K_2 \in \{\R, \C\}$, we extend its definition to a function $\mathcal{F}: \K_1^\J\to \K_2^\J$, which applies $F$ pointwise, i.e.\ for $u\in \K_1^\J$, $(\mathcal{F}(u))_j = F(u_j), j\in \dbra{0}{\J-1}$. For the sake of simplicity, we keep their original notation when there is no ambiguity. For instance, the real part function $\Re(\cdot): \C \to \R$ and the absolute value $|\cdot|: \C \to \R$ are extended to vectors as $\Re(u) = (\Re(u_1), \Re(u_2), ..., \Re(u_\J))^T \in \R^\J$ and $|u| = (|u_1|, |u_2|, ..., |u_\J|)^T\in \R^\J$, respectively.

We also extend functions of the form $G: \K_1 \times \K_2 \to \K_3$, with $\K_3=\R, \C$, to $\mathcal{G}: \big(\K_1 \times \K_2\big)^\J \to \K_3^\J$ applying $G$ pointwise, i.e.\ $(G(u, v))_j = \mathcal{G}(u_j, v_j)$. As before, we keep the original notation when there is no ambiguity. For example, the inner product $\ip{\cdot}{\cdot}$ defined above is extended to $\ip{u}{v}_\odot = (\ip{u_1}{v_1}, \ip{u_2}{v_2}, ..., \ip{u_\J}{v_\J}) = \Re(u \odot \conj{v}) \in \R^\J$. Therefore, the square modulus $|\cdot|^2$ of a vector is well defined in the sense that $(|u|^2)_j = |u_j|^2$, and it satisfies $|u|^2 = \ip{u}{u}_\odot = u\odot \conj{u}$. Notice that the pointwise product $\odot$ can be seen as the extension to $\C^\J$ of the usual complex product $z_1\cdot z_2$.
\paragraph{Discrete operators.}

For $W, Z \in \C^{\J}$ and $U, V \in \C^{\J+1}$, we define the following sesquilinear forms:
\begin{equation}\label{eq:C_inner}
    \ipCm{W}{Z} = \dx \sum_{j=0}^{\J-1} W_j\conj{Z}_j, \quad \ipC{U}{V} = \dx\sum_{j=0}^{\J} U_j\conj{V}_j, \quad \ipCp{U}{V} = \dx \Bigg(\frac{U_0\conj{V}_0}{2} + \sum_{j=1}^{\J-1} U_j\conj{V}_j +\frac{U_J\conj{V}_J}{2} \Bigg),
\end{equation}
and the associated inner products by taking the real part, i.e.
\begin{equation}\label{eq:inner}
    \ipRm{W}{Z} = \Re \ipCm{W}{Z}, \quad \ipR{U}{V} = \Re \ipC{U}{V} \quad \text{ and } \quad  \ipRp{U}{V} = \Re \ipCp{U}{V}.
\end{equation}
We emphasize that the subscripts $\dx^-$, $\dx^+$ and $\dx$ are not merely notational, since they correspond to slightly different discrete products. 
Notice that $\ipR{U}{iV} = -\ipR{iU}{V} = \Im \ipC{U}{V}$, and the same holds for $\ipRp{\cdot}{\cdot}$ and $\ipRm{\cdot}{\cdot}$.

For $p\in \N\setminus \{0\}$, we define the following discrete norms:
\begin{equation}\label{eq:norms}
    \norm{W}_{p, \dx^-}^p = \dx \sum_{j=0}^{\J-1} |W_j|^p, \quad \norm{U}_{p, \dx}^p = \dx \sum_{j=0}^{\J} |U_j|^p, \quad \norm{U}_{p, \dx^+}^p = \dx \left( \frac{1}{2} |U_0|^p + \sum_{j=1}^{\J-1} |U_j|^p + \frac{1}{2} |U_J|^p\right),
\end{equation}
and the discrete infinity-norm and $H^1$-norm:
\begin{equation}\label{eq:norms_H1}
    \norm{U}_\infty = \max_{0\leq j\leq \J}|U_j|, \quad \norm{U}_{H^1_\dx}^2 = \norm{U}_{2, \dx^+}^2 + \norm{\Dxp U}_{2, \dx^-}^2.
\end{equation}

For $p=2$, the norms defined in \eqref{eq:norms} are induced from the inner products in \eqref{eq:inner}. Therefore, they satisfy all properties shared by inner products and norms, such as the Cauchy--Schwarz inequality. For instance,
\begin{equation}\label{Cauchy--Schwarz}
\abs{\ipR{U}{V}} \leq \normp{U}{2}\normp{V}{2},
\quad 
\abs{\ipRm{W}{Z}} \leq \normpm{W}{2}\normpm{Z}{2},
 \quad 
  \abs{\ipRp{U}{V}} \leq \normpp{U}{2}\normpp{V}{2}.
\end{equation}

To approximate the spatial derivatives and formulate the discrete problem, we define for $U\in \C^{\J+1}$ the central finite-difference operators $\Dx$ and $\Dxx$ from $\C^{\J+1}$ to itself by
\begin{equation}\label{eq:spc_disc_ops}
    (\Dx U)_j = \frac{U_{j+1} - U_{j-1}}{2\dx}, \quad \text{ and } \quad (\Dxx U)_j = \dfrac{U_{j-1} - 2U_{j} + U_{j+1}}{\dx^2},\quad j\in \dbra{1}{\J-1},
\end{equation}
in the interior nodes, and the extension to the boundary nodes $j=0$ and $j=\J$ by imposing the ghost points $U_{-1} = U_1$ and $U_{\J+1} = U_{\J-1}$ to implement the Neumann boundary conditions, which yields
\begin{equation*}
    (\Dx U)_0= (\Dx U)_\J = 0, \quad (\Dxx U)_0 = \dfrac{2U_{1} - 2U_{0}}{\dx^2}, \quad \text{ and }\quad (\Dxx U)_\J = \dfrac{2U_{J-1} - 2U_{J}}{\dx^2}.
\end{equation*}
We also define the forward finite-difference operator $\Dxp$ defined from $\C^{\J+1}$ to $\C^\J$ by:
\begin{equation}\label{eq:forward}
    (\Dxp U)_j = \frac{U_{j+1} - U_{j}}{\dx}, \quad j\in \dbra{0}{\J-1}.
\end{equation}

\section{Preliminaries}\label{juanito}
In this section, we provide several useful identities and inequalities for the discrete norms and inner products defined above that we will often use in the numerical analysis of our scheme.
\begin{lemma}\label{lemma:discrete_product}
    Let $U, V, W, Z\in \C^\J$. The following property holds:
    \begin{equation}\label{eq:discrete_product}
        \ip{U}{V}_\odot - \ip{W}{Z}_\odot = \ip{\frac{U+W}{2}}{V-Z}_\odot + \ip{\frac{V+Z}{2}}{U-W}_\odot.
    \end{equation}
\end{lemma}
\begin{proof}
    Expand the right-hand side of \eqref{eq:discrete_product} and use the symmetry $\ip{U}{V}_\odot = \ip{V}{U}_\odot$.
\end{proof}
\begin{corollary}[Product rule]\label{corollary:product_rule} Let $U, V\in \C^{\J+1}$ and define $U_{j+\frac{1}{2}} =(U_{j+1} + U_j)/2$ for $j\in \dbra{0}{\J-1}$. Then, the discrete operator $\Dxp$ satisfies
    \begin{equation}\label{eq:product_rule}
        (\Dxp \ip{U}{V}_\odot)_j = \ip{U_{j+\frac{1}{2}}}{(\Dxp V)_j} + \ip{V_{j+\frac{1}{2}}}{(\Dxp U)_j}, \quad  j\in \dbra{0}{\J-1},
    \end{equation}
    In particular, 
    \begin{equation}\label{eq:product_U^2}
        (\Dxp |U|^2)_j =2\ip{U_{j+\frac{1}{2}}}{(\Dxp U)_j}, \quad j\in \dbra{0}{\J-1}.
    \end{equation}
\end{corollary}
\begin{proof}
    To prove \eqref{eq:product_rule}, use the definition of $\Dxp$ in \eqref{eq:forward} and apply Lemma~\ref{lemma:discrete_product} with $\J=1$ and $U = U_{j+1}$, $V=V_{j+1}$, $W=U_j$ and $Z=V_j$.

    In order to obtain \eqref{eq:product_U^2} simply set $V=U$ in \eqref{eq:product_rule}.
\end{proof}
\begin{lemma}[Summation by parts]\label{lemma:sum_by_parts}
        For $U, V \in \C^{\J+1}$, the product $\ipRp{\cdot}{\cdot}$ satisfies
        \begin{equation*}
            \ipRp{\Dxx U}{V} = -\ipRm{\Dxp U}{\Dxp V}.
        \end{equation*}
    \end{lemma}
\begin{proof}
    Let $U = (U_0, U_1, ..., U_{\J})^T$ and $V = (V_0, V_1, ..., V_{\J})^T$. Using the definition of $\Dxx$ and $\ipCp{\cdot}{\cdot}$, we have
    \begin{align*}
       \ipCp{\Dxx U}{V} &= \dx\left[\frac{U_1 - U_0}{\dx^2}\conj{V}_0+ \sum_{j=1}^{\J-1} \frac{U_{j-1} - 2U_j + U_{j+1}}{\dx^2} \conj{V}_j + \frac{U_{\J-1} - U_{\J}}{\dx^2}\conj{V}_\J \right] \\
        &=  \dx\left[\frac{U_1 - U_0}{\dx^2}\conj{V}_0+ \sum_{j=1}^{\J-1} \frac{(U_{j+1} - U_{j})}{\dx^2} \conj{V}_j - \sum_{j=1}^{\J-1}\frac{(U_j - U_{j-1})}{\dx^2} \conj{V}_j - \frac{U_{\J} - U_{\J-1}}{\dx^2}\conj{V}_\J\right] \\
        &=\dx \left[\sum_{j=0}^{\J-1} \frac{(U_{j+1} - U_{j})}{\dx^2} \conj{V}_j - \sum_{j=1}^{\J}\frac{(U_{j} - U_{j-1})}{\dx^2} \conj{V}_{j}\right] \\
        &=\dx \left[\sum_{j=0}^{\J-1} \frac{(U_{j+1} - U_{j})}{\dx^2} \conj{V}_j - \sum_{j=0}^{\J-1}\frac{(U_{j+1} - U_{j})}{\dx^2} \conj{V}_{j+1}\right] \\
        &= -\dx \sum_{j=0}^{\J-1}  \frac{(U_{j+1} - U_{j})}{\dx}\frac{(\conj{V}_{j+1} - \conj{V}_{j})}{\dx},
    \end{align*}
    and the lemma follows by taking the real part.
\end{proof}

\begin{lemma}[Distributive property]\label{lemma:distributive}
    Let $U\in \R^{\J+1}$ and  $V, W\in \C^{\J+1}$. Then,
    \begin{equation*}
        \ipR{U\odot V}{W} = \ipR{U}{\ip{V}{W}_\odot}.
    \end{equation*}
    In particular,
    \begin{equation*}
        \ipR{U\odot V}{V} = \ipR{U}{|V|^2}.
    \end{equation*}
    The same holds for $\ipRm{\cdot}{\cdot}$ and $\ipRp{\cdot}{\cdot}$
\end{lemma}
\begin{proof}
    Use the definition of $\ip{\cdot}{\cdot}_\odot$ and each inner product in \eqref{eq:inner}.
\end{proof}

\begin{lemma}\label{lemma:Dx_Dxx_bounds}
    Let $U\in \C^{\J+1}$ be such that $\normpp{U}{2}\neq0$. Then, the discrete operators $\Dxp$ and $\Dxx$ satisfy the following upper bounds:
    \begin{equation}\label{ineq:Dx_Dxx_bounds}
        \normpm{\Dxp U}{2} \leq \frac{2}{\dx} \normpp{U}{2} \quad \text{ and }\quad \normpp{\Dxx U}{2} \leq \frac{4}{\dx^2} \normpp{U}{2}.
    \end{equation}
    
\end{lemma}
\begin{proof}
For the first bound in \eqref{ineq:Dx_Dxx_bounds}, the property $|a - b|^2 \leq 2|a|^2 + 2|b|^2$ with the definition of $\Dxp$ in \eqref{eq:forward} and $\normpm{\cdot}{2}$ in \eqref{eq:norms} give
    \begin{equation*}
        \normpm{\Dxp U}{2}^2 = \frac{1}{\dx} \sum_{j=0}^{\J-1} |U_{j+1} - U_j|^2 \leq \frac{2}{\dx} \sum_{j=0}^{\J-1} |U_{j+1}|^2 + |U_j|^2 = \frac{4}{\dx^2} \normpp{U}{2}^2.
    \end{equation*}
For the second, we use Lemma~\ref{lemma:sum_by_parts} and Cauchy--Schwarz \eqref{Cauchy--Schwarz} to obtain
\begin{align*}
    \normpp{\Dxx U}{2}^2 &= -\ipRm{\Dxp U}{\Dxp(\Dxx U)} \leq \normpm{\Dxp U}{2} \normpm{\Dxp(\Dxx U)}{2} \\& \leq \frac{4}{\dx^2} \normpp{U}{2} \normpp{\Dxx U}{2}, 
\end{align*}
which concludes the proof.
\end{proof}

For the subsequent properties, it is useful to introduce the linear interpolation function $u:\mathcal{J}\to \C$ of $U\in \C^{\J+1}$, defined by
    \begin{equation}\label{def:linear_intp}
        u(x) = \sum_{j=0}^{\J} U_j \phi_j(x),~\text{ where }\phi_j(x)= \max\left(0, 1 - \left|\frac{x-x_j}{\dx}\right|\right).
    \end{equation}
\begin{lemma}[Equivalence of norms]\label{lemma:equivalence}
    Let $u$ be the linear interpolation of $U\in\C^{\J+1}$ defined in \eqref{def:linear_intp}. Then
    \begin{equation}\label{eq:equivalence}
        \frac{1}{2^{p+2}} \norm{U}_{p, \dx^+}^p \leq \norm{u}_{L^p(\mathcal{J})}^p \leq \norm{U}_{p, \dx^+}^p,
    \end{equation}
    for all $p\in\N\backslash \{0\}$.
\end{lemma}
\begin{proof}
    We first prove the upper bound in \eqref{eq:equivalence}. Let us fix $x\in[x_j, x_{j+1}]$, so that $u(x) = U_j\phi_j(x) + U_{j+1}\phi_{j+1}(x)$. Noticing that $\phi_j(x), \phi_{j+1}(x)\geq 0$ and $\phi_j(x) + \phi_{j+1}(x) = 1$, the convexity of the function $z\to |z|^p$ yields
    \begin{equation*}
        |u(x)|^p \leq |U_j|^p \phi_j(x) + |U_{j+1}|^p \phi_{j+1}(x).
    \end{equation*}
    Hence, using the definition of the $L^p$-norm, we have
    \begin{equation*}
        \norm{u}_{L^p(\mathcal{J})}^p = \sum_{j=0}^{\J-1} \int_{x_j}^{x_{j+1}} |u(x)|^pdx \leq \sum_{j=0}^{\J-1}  |U_j|^p \int_{x_j}^{x_{j+1}} \phi_j(x)dx + \sum_{j=0}^{\J-1} |U_{j+1}|^p \int_{x_j}^{x_{j+1}}\phi_{j+1}(x)dx.
    \end{equation*}
    Since $\int_{x_j}^{x_{j+1}} \phi_j(x)dx = \int_{x_j}^{x_{j+1}}\phi_{j+1}(x)dx = \frac{\dx}{2}$, it follows that
    \begin{equation*}
        \norm{u}_{L^p(\mathcal{J})}^p \leq \frac{\dx}{2}\sum_{j=0}^{\J-1}  |U_j|^p +  |U_{j+1}|^p = \dx \left( \frac{1}{2} |U_0|^p + \sum_{j=1}^{\J-1} |U_j|^p + \frac{1}{2} |U_\J|^p\right) = \norm{U}_{p, \dx^+}^p.
    \end{equation*}
    To establish the lower bound in \eqref{eq:equivalence}, we fix $j$ and work on $[x_j, x_{j+1}]$. First, we assume that $|U_{j+1}|\leq |U_{j}|$, so for all $x\in[x_j, x_{j+1}]$,
    \begin{equation*}
        (\phi_j(x) -\phi_{j+1}(x))|U_j| \leq |U_j|\phi_j(x) - |U_{j+1}|\phi_{j+1}(x) \leq |u(x)|.
    \end{equation*}
    In the subinterval $\left[x_j, x_j+\dx/4\right]$, we note that $\phi_j(x) -\phi_{j+1}(x) \geq 1/2$. Hence,
    \begin{equation*}
        \frac{1}{2^{p+1}} (|U_j|^p + |U_{j+1}|^p)\leq \frac{1}{2^p}|U_j|^p \leq |u(x)|^p,\quad \text{for all } x\in \left[x_j, x_j+\dx/4\right].
    \end{equation*}
    Integrating from $x_j$ to $x_j+\dx/4$, we conclude that
    \begin{equation}\label{ineq:lower_equivalence_norms}
        \frac{\dx}{2^{p+3}} (|U_j|^p + |U_{j+1}|^p) \leq \int_{x_j}^{x_{j}+\dx/4} |u(x)|^p dx \leq \int_{x_j}^{x_{j+1}} |u(x)|^p dx.
    \end{equation}
    The same inequality holds if $|U_{j}|\leq |U{j+1}|$, in the subinterval $\left[x_{j+1}-\dx/4, x_{j+1}\right]$. Thus, by summing \eqref{ineq:lower_equivalence_norms} from $j=0$ to $j=\J-1$, we obtain the lower inequality in \eqref{eq:equivalence} as before, concluding the proof.
\end{proof}
\begin{lemma}[Discrete Gagliardo--Nirenberg's inequality]\label{lemma:discrete_GN}
    For $U\in \C^{\J+1}$, there exists a constant $C_L>0$ depending only on $L$, such that 
        \begin{equation}\label{eq:GN_disc}
        \norm{U}_{4,\dx,+}^4 \leq C_L\left(\norm{U}_{2,\dx,+}^3\norm{D_x^{h,+}U}_{2,\dx^-}+\norm{U}_{2,\dx,+}^4\right).
        \end{equation}
    \end{lemma}
\begin{proof}
   Let $u$ be the linear interpolation of $U$ in \eqref{def:linear_intp}. Since $u$ is continuous and piecewise affine on $\mathcal{J}$, so that $u\in H^1(\mathcal{J})$. Therefore, the one-dimensional Gagliardo--Nirenberg inequality, established in \cite{Burenkov}, applied with the $L^4, L^2$ and $H^1$-norms yields
    \begin{equation}\label{eq:GN_inter}
        \norm{u}_{L^4(\mathcal{J})}^4 \leq  C_L \left(\norm{u}^3_{L^2(\mathcal{J})} \norm{u'}_{L^2(\mathcal{J})} + \norm{u}_{L^2(\mathcal{J})}^4\right),
    \end{equation}
    where $C_L>0$ depends only on $L$. Given that 
    \begin{equation*}
        \norm{u'}^2_{L^2(\mathcal{J})} = \sum_{j=0}^{\J-1} \int_{x_j}^{x_{j+1}} \left|\frac{u_{j+1} - u_j}{\dx}\right|^2 = \dx \sum_{j=0}^{\J-1} \left|\frac{u_{j+1} - u_j}{\dx}\right|^2 = \norm{D_x^{h,+}U}_{2,\dx^-}^2,
    \end{equation*}
    the proof of \eqref{eq:GN_disc} is a consequence of Lemma~\ref{lemma:equivalence}. Indeed, applying \eqref{eq:equivalence} with $p=2, 4$ in \eqref{eq:GN_inter}, we conclude that
    \begin{align*}
        \norm{U}_{4,\dx,+}^4 \leq \norm{u}_{L^4(\mathcal{J})}^4 \leq&  C_L \left(\norm{u}^3_{L^2(\mathcal{J})} \norm{u'}_{L^2(\mathcal{J})} + \norm{u}_{L^2(\mathcal{J})}^4\right)
        \\\leq& C_L\left(\norm{U}_{2,\dx,+}^3\norm{D_x^{h,+}U}_{2,\dx^-}+\norm{U}_{2,\dx,+}^4\right),
    \end{align*}
    which completes the proof.
\end{proof}

\begin{lemma}[Discrete Sobolev embedding]\label{lemma:discrete_sobolev}
    The following holds for $U\in \C^{\J+1}$:
    \begin{equation}\label{ineq:discrete_sobolev}
        \norm{U}_\infty \leq C_L\norm{U}_{H^1_\dx},
    \end{equation}
    where $C_L>0$ depends only on $L$.
\end{lemma}
\begin{proof}
    Let $u$ be the linear interpolation of $U$. As in Lemma~\ref{lemma:discrete_GN}, $u\in H^1(\boJ)$. Using the embedding $H^1(\boJ)\hookrightarrow L^\infty(\boJ)$ and the equivalence of norms in Lemma~\ref{lemma:equivalence}, we have
    \begin{equation*}
        \norm{U}_{L^\infty} = \norm{u}_{L^\infty (\boJ)}\leq c_L\norm{u}_{H^1(\boJ)} \leq C_L \norm{U}_{H^1_\dx},
    \end{equation*}
    where $c_L$ and $C_L$ depend only on $L$, giving inequality \eqref{ineq:discrete_sobolev} and concluding the proof.
\end{proof}
\begin{lemma}[Discrete $H^1$-algebra property]\label{lemma:discrete_H1_algebra}
    The following holds for $U, V\in \C^{\J+1}$:
    \begin{equation*}
        \norm{UV}_{H^1_\dx} \leq C\norm{U}_{H^1_\dx}\norm{V}_{H^1_\dx},
    \end{equation*}
    where $C>0$ is independent of $\dx$.
\end{lemma}
\begin{proof}
    Let $u, v$ be the linear interpolation of $U$ and $V$, respectively. As in Lemma~\ref{lemma:discrete_GN}, $u, v\in H^1(\boJ)$. Using that $H^1(\boJ)$ is an algebra in dimension 1, and the equivalence of norms in Lemma~\ref{lemma:equivalence}, we have
    \begin{equation*}
        \norm{UV}_{H^1_\dx} \leq C \norm{uv}_{H^1(\boJ)} \leq C\norm{u}_{H^1(\boJ)}\norm{v}_{H^1(\boJ)}\leq C\norm{U}_{H^1_\dx}\norm{V}_{H^1_\dx}
    \end{equation*}
    where $C$ does not depend on $\dx$, concluding the proof.
\end{proof}
\begin{lemma}[Discrete H\"older]\label{lemma:discrete_holder}
    The following holds for $U, V\in \C^{\J+1}$:
    \begin{equation*}
        \normpp{UV}{2} \leq \norm{U}_\infty \normpp{V}{2}.
    \end{equation*}
\end{lemma}
\begin{proof}
    Again, the proof follows directly by using the linear interpolation of $U$ and $V$, the continuous H\"older inequality, and the equivalence of norms in Lemma~\ref{lemma:equivalence}.
\end{proof}

\section{Formulation of the Crank--Nicolson scheme}\label{pablito}
We recall our Crank--Nicolson scheme for \eqref{QGP}:
\begin{equation}\label{scheme}
   \frac{i\left(\Psi^{n+1} -\Psi^{n}\right)}{\dt} = \Dxx\Psi^{n+\frac{1}{2}} + \gs \Psi^{n+\frac{1}{2}} \odot \boQ^{n+\frac{1}{2}} + \ka \Psi^{n+\frac{1}{2}} \odot \Dxx \boQ^{n+\frac{1}{2}}.
\end{equation}
We define the discrete version of the energy in \eqref{eq:energie} by
\begin{equation}\label{eq:full_energy}
    E_\dx(U) = \frac{1}{2}\normpm{\Dxp U}{2}^2 - \frac{\gs}{4}\normpp{U}{4}^4 + \frac{\ka}{4}\normpm{\Dxp \boQ}{2}^2,
\end{equation}
and the discrete version of the mass in \eqref{eq:mass} by 
\begin{equation}\label{eq:full_mass}
    N_\dx(U) = \normpp{U}{2}^2,
\end{equation}
using the norms and inner products defined in \eqref{eq:inner} and \eqref{eq:norms}.
As addressed in our companion paper \cite[Proposition~2.1]{de_laire_simul_2026}, our CN scheme \eqref{scheme} is derived from energy and mass preservation at the continuous level in space and is shown to preserve these quantities in the fully discrete setting. This property is established in the following proposition:
\begin{proposition}\label{prop:full_pres}
    If the sequence $(\Psi^n)_{n=0}^\Nt\subset\C^{J+1}$ satisfies the scheme \eqref{scheme}, then the quantities $E_\dx$ and $N_\dx$ satisfy
    \begin{equation*}
        E_\dx(\Psi^n) = E_\dx(\Psi^0)\quad \text{ and } \quad N_\dx(\Psi^n) = N_\dx(\Psi^0), \quad \text{ for all } n\in \dbra{0}{\Nt}.
    \end{equation*}
\end{proposition}
The proof can be found in our companion paper \cite[Proposition~2.3]{de_laire_simul_2026}. In what follows, we will denote them as $E_\dx$ and $N_\dx$ for simplicity.

\subsection{Stability}\label{sec:stability}

The following proposition establishes the $H^1$-stability of the scheme \eqref{scheme} as a consequence of the discrete preservation of mass and energy. However, when $\ka<0$, the quasilinear contribution to the energy is no longer positive, and energy conservation alone does not provide the estimate. To recover the stability, we need to assume that $\boD_n\geq \delta$ for some $\delta>0$. With this, we can control the discrete $H^1$-norm of the solution.

\begin{proposition}[$H^1$-stability]\label{prop:stability}
    Let $(\Psi^n)_{n=0}^\Nt$ be a sequence. In the case $\ka<0$, assume that the sequence satisfies $\boD_n\geq \delta$ for some $\delta >0$. If $(\Psi^n)_{n=0}^\Nt$ satisfies the scheme \eqref{scheme}, then
    \begin{equation*}
        \norm{\Psi^n}^2_{H^1_\dx}  < C,
    \end{equation*}
    for all $n\in \dbra{0}{\Nt}$, where $C$ is independent of $\dx$ and $\dt$ and it also depends on $\delta$ if $\ka<0$.
\end{proposition}
\begin{proof}
Let $n\in \dbra{0}{\Nt}$. We first prove the case $\ka>0$. We recall that $E_\dx$ denotes the energy of the scheme and $N_\dx$ the mass, both defined in \eqref{eq:full_energy} and \eqref{eq:full_mass}, respectively.
Since $\ka > 0$, we have
    \begin{equation}\label{eq:energy_bound}
        \frac{1}{2}\norm{\Dxp\Psi^n}_{2, \dx^-}^2 - \frac{\gs}{4}\norm{\Psi^n}_{4, \dx^+}^4 \leq E_\dx
    \end{equation}
    Notice that, in the case $\gs=-1$, the $H^1$-stability follows directly. In the case $\gs=1$, applying Lemma~\ref{lemma:discrete_GN} on $\Psi^n$ yields
    \begin{align*}
        \norm{\Psi^n}_{4, \dx^+}^4 &\leq C\left( \norm{\Psi^n}_{2, \dx^+}^3  \norm{\Dxp\Psi^n}_{2, \dx^-} +  \norm{\Psi^n}_{2, \dx^+}^4\right)\\
        &= C(N_\dx^{3/2}\norm{\Dxp\Psi^n}_{2, \dx^-} + N_\dx^2)
    \end{align*}
    where $C$ is independent of $\dx$ and $\dt$. It follows that
    \begin{equation*}
        \frac{1}{2}\norm{\Dxp\Psi^n}_{2, \dx^-}^2 - C N_\dx^{3/2}\norm{\Dxp\Psi^n}_{2, \dx^-} - C N_\dx^2 \leq E_\dx.
    \end{equation*}
    For all $\varepsilon>0$, by the Young inequality $ab \leq \varepsilon a^2 + b^2/(4\varepsilon)$, we can upper bound the second term on the left-hand side as
    \begin{equation}\label{young}
        CN_\dx^{3/2}\norm{\Dxp\Psi^n}_{2, \dx^-} \leq \frac{C^2}{4\varepsilon}\left(N_\dx\right)^3 + \varepsilon \norm{\Dxp\Psi^n}^2_{2, \dx^-}, 
    \end{equation}
    and therefore,
    \begin{equation*}
        \left(\frac{1}{2} - \varepsilon\right)\norm{\Dxp\Psi^n}_{2, \dx^-}^2 - C \left(N_\dx\right)^2  - C_\varepsilon \left(N_\dx\right)^3 \leq E_\dx,
    \end{equation*}
    where $C_\varepsilon = C^2/4\varepsilon$. Taking $\varepsilon=1/4$, we conclude that
    \begin{equation*}
        \norm{\Psi^n}^2_{H^1_\dx} = \norm{\Psi^n}_{2, \dx^+}^2 + \norm{\Dxp\Psi^n}^2_{2, \dx^-} \leq C\left(E_\dx + N_\dx + N_\dx^2 + N_\dx^3\right).
    \end{equation*}
Now, for the case $\ka<0$, we estimate the quasilinear term in the energy to obtain a lower bound as in \eqref{eq:energy_bound}. For all $j\in \dbra{0}{J-1}$, we have
\begin{equation*}
\left|\left(D_x^{\dx,+}\boQ^n\right)_j\right| = \frac{\left|\boQ_{j+1}^n-\boQ_j^n\right|}{\dx} \leq \left(|\Psi_{j+1}^n|+|\Psi_j^n|\right) \frac{|\Psi_{j+1}^n-\Psi_j^n|}{\dx},
\end{equation*}
regarding that $\boQ_j^n = |\Psi_j^n|^2$. Hence,
\begin{equation*}
\left|\left(D_x^{\dx,+}\boQ^n\right)_j\right|\leq2\norm{\Psi^n}_{\infty}\left| \left(D_x^{\dx +}\Psi^n\right)_j\right|.
\end{equation*}
Multiplying by the weights in the discrete norm \eqref{eq:norms} and summing over $j$, we obtain
\begin{equation}\label{ineq:density_gradient_bound}
\norm{D_x^{\dx,+}\boQ^n}_{2,\dx^-}^2 \leq 4\norm{\Psi^n}_{\infty}^2 \norm{D_x^{\dx,+}\Psi^n}_{2,\dx^-}^2.
\end{equation}
Since $\ka<0$, multiplying \eqref{ineq:density_gradient_bound} by $\ka/4$ reverses the inequality. Therefore,
\begin{equation*}
\ka\norm{\Psi^n}_{\infty}^2\norm{D_x^{\dx,+}\Psi^n}_{2,\dx^-}^2 \leq \frac{\ka}{4}\norm{D_x^{\dx,+}\boQ^n}_{2, \dx^-}^2.
\end{equation*}
It follows that
\begin{equation*}
\left( \frac{1}{2} + \ka\norm{\Psi^n}_{\infty}^2 \right) \norm{D_x^{\dx,+}\Psi^n}_{2,\dx^-}^2 - \frac{\gs}{4}\norm{\Psi^n}_{4,\dx,+}^4 \leq E_\dx.
\end{equation*}
By the discrete ellipticity condition \eqref{eq:disc_ellipticity}, we obtain
\begin{equation*}
\frac{\delta}{2} \norm{D_x^{\dx,+}\Psi^n}_{2,\dx^-}^2 - \frac{\gs}{4}\norm{\Psi^n}_{4,\dx,+}^4 \leq E_\dx.
\end{equation*}
For $\gs=-1$, the $H^1$-stability follows directly. For $\gs=1$, we use Lemma~\ref{lemma:discrete_GN} and inequality \eqref{young} with $\varepsilon=\delta/4$ to conclude the proof.

\end{proof}

\subsection{Existence and uniqueness of the discrete problem}\label{subec:discrete_existence}
The existence of a solution of the CN scheme \eqref{scheme} will be based on the following Brouwer-type zero lemma, which is a consequence of Brouwer's fixed-point theorem. We refer to \cite{akrivis_dougalis_karakashian_1991} for the proof, and to \cite{akrivis_1993, Vougiantzis_2026} for some applications in the context of nonlinear Schr\"odinger equations.
\begin{lemma}[Brouwer-type zero lemma]\label{lem:brouwer-zero}
Let $H$ be a finite-dimensional complex inner-product space, with norm
$\norm{\cdot}_H$, and let $\Lambda:H\to H$ be continuous.  If there
exists $\beta>0$ such that
\begin{equation}\label{eq:brouwer-outward}
 \Re(\Lambda(z),z)_H>0 \qquad\text{whenever }\norm{z}_H=\beta.
\end{equation}
Then there exists $\widehat z\in H$, with $\norm{\widehat z}_H\le\beta$, such that $\Lambda(\widehat z)=0$.
\end{lemma}
\begin{proposition}[Existence]
\label{prop:discrete-existence}
Let $\Psi^0\in\C^{J+1}$ be given.  Then, there exists at least one sequence $(\Psi^n)_{n=0}^\Nt\subset\C^{J+1}$ satisfying \eqref{scheme}. 
\end{proposition}

\begin{proof}
For $U\in \C^{\J+1}$ fixed, we define the continuous function $\Lambda_U:\C^{J+1}\to\C^{J+1}$ given by
\begin{equation*}
 \Lambda_U(Z) =Z-U+\frac{i\tau}{2}\left[ \Dxx Z+\gs \mathcal{M}_U(Z) Z+\ka\Big(\Dxx \mathcal{M}_U(Z)\Big)Z \right],
\end{equation*}
for $Z\in \C^{\J+1}$, where $\mathcal{M}_U(Z) = (|U|^2+|2Z-U|^2)/2$. The idea is to apply Lemma~\ref{lem:brouwer-zero} to the function $\Lambda_U$ with the complex operator $\ipCp{\cdot}{\cdot}$ defined in \eqref{eq:C_inner}. To do so, it remains to verify the existence of a radius $\beta>0$ such that \eqref{eq:brouwer-outward} holds.

We claim that
\begin{equation}\label{eq:AU-real}
 \ipCp{\Dxx Z+\gs \mathcal{M}_U(Z) Z+\ka\Big(\Dxx \mathcal{M}_U(Z)\Big)Z}{Z}\in\R \quad\text{for every }Z\in\C^{J+1}.
\end{equation}
Indeed, by Lemma~\ref{lemma:sum_by_parts}, Lemma~\ref{lemma:distributive} and using the fact that $\mathcal{M}_U(Z)$ and $\Dxx \mathcal{M}_U(Z)$ are real-valued, we have
\begin{equation}\label{eq:dxx_cub_real}
 \ipCp{\Dxx Z}{Z} =-\normpm{\Dxp Z}{2}^2\in\R, \quad \ipCp{\mathcal{M}_U(Z)\odot Z}{Z} =\ipCp{\mathcal{M}_U(Z)}{|Z|^2}\in\R,
\end{equation}
and
\begin{equation}\label{eq:quasi_real}
    \ipCp{(\Dxx \mathcal{M}_U(Z))\odot Z}{Z} =\ipCp{\Dxx \mathcal{M}_U(Z)}{|Z|^2}\in\R.
\end{equation}
Hence, \eqref{eq:AU-real} holds. By the definition of the inner product $\ipRp{\cdot}{\cdot}$ in \eqref{eq:inner} and the discrete Cauchy--Schwarz inequality in \eqref{Cauchy--Schwarz}, we obtain 
\begin{equation}\label{ineq:lambda_U_Z}
     \ipRp{\Lambda_U(Z)}{Z} = \Re \ipCp{\Lambda_U(Z)}{Z} = \normpp{Z}{2}^2-\ipRp{U}{Z} \geq\norm{Z}_{2,h,+} \big(\normpp{Z}{2}-\normpp{U}{2}\big).
\end{equation}
Choosing $\beta=1+\normpp{U}{2}$, if $Z$ satisfies $\normpp{Z}{2}=\beta$, then \eqref{ineq:lambda_U_Z} implies $\ipRp{\Lambda_U(Z)}{Z} \geq \beta > 0$. By Lemma~\ref{lem:brouwer-zero}, there exists $\widehat Z\in\C^{J+1}$ such that $\normpp{\widehat Z}{2}\leq \beta$ and $\Lambda_U(\widehat Z)=0$.

Now, let us fix $n\in \dbra{0}{\Nt}$ and suppose that $\Psi^n \in \C^{\J+1}$ is known. For $U = \Psi^n$, the calculations above give us the existence of $\widehat{Z}\in \C^{\J+1}$ such that $\Lambda_U(\widehat{Z}) = 0$. Defining $\Psi^{n+1} = 2\widehat{Z} - \Psi^n$, it is easy to check that $\Psi^{n+1}$ satisfies the scheme \eqref{scheme} from the fact that $\Lambda_U(\widehat{Z}) = 0$. The proof follows by induction on $n$, starting from $n=0$.
\end{proof}

The uniqueness of the discrete solution is established in the following proposition, where it is necessary to assume a smallness condition between the steps $\dt$ and $\dx$ of order $\dt = \boO(\dx^2)$.

\begin{proposition}[Uniqueness]
\label{prop:cn-uniqueness}
Let $\Psi^0\in\C^{J+1}$ be given. If $\ka<0$, assume the discrete ellipticity condition $\boD_n\geq \delta >0$ for some $\delta$. Then, there exists a constant $M>0$ independent of $\dx$ and $\dt$  such that: if $\dt, \dx >0$ satisfy
\begin{equation}\label{eq:cn-uniqueness-condition}
 M\Big(\dt+\frac{4\dt|\ka|}{\dx^2}\Big)<1,
\end{equation}
then, there exists a unique sequence $(\Psi^n)_{n=0}^N\subset\C^{J+1}$ satisfying \eqref{scheme}. 
\end{proposition}

\begin{proof}
Existence follows from Proposition~\ref{prop:discrete-existence}. Let us fix $n\in \dbra{0}{\Nt}$ and suppose that we know $U:= \Psi^n$, the solution at time step $n$. Let $X, Y\in \C^{\J+1}$ be two candidates for the next time step, both satisfying \eqref{scheme}. From Proposition~\ref{prop:full_pres}, we have
\begin{equation*}
    E_\dx(U) = E_\dx(X) = E_\dx(Y) = E_\dx \quad \text{ and } \quad N_\dx(U) = N_\dx(X) = N_\dx(Y) = N_\dx.
\end{equation*}
Let $W=X-Y$. Substracting the respective equations \eqref{scheme} for $X$ and for $Y$ gives
\begin{equation}\label{eq:cn-difference}
 \frac{2i}{\dt} W =\Dxx W+N_X\odot W +(N_X-N_Y)\odot(U+Y),
\end{equation}
where $N_A:=\gs Q_A+\ka \Dxx Q_A$ and $Q_A := \big(|U|^2 + |A^2|\big)/2$ for $A\in \{X, Y\}$. Taking $\ipCp{\cdot}{W}$ and the imaginary part of equation \eqref{eq:cn-difference}, the first and second terms on the right-hand side of \eqref{eq:cn-difference} vanish by the same arguments as in \eqref{eq:dxx_cub_real} and \eqref{eq:quasi_real}. Therefore,
\begin{equation}\label{ineq:W_bound}
\begin{aligned}
    \frac{2}{\dt} \normpp{W}{2}^2 &
    \leq \abs{\ipCp{(N_X - N_Y)\odot (U+Y)}{W}} \leq \normpp{(N_X - N_Y)\odot (U+Y)}{2}\normpp{W}{2}
    \\ & \leq \norm{U+Y}_\infty \normpp{N_X - N_Y}{2} \normpp{W}{2},
\end{aligned}
\end{equation}
where we used 
the Cauchy--Schwarz inequality in \eqref{Cauchy--Schwarz} and Lemma~\ref{lemma:discrete_holder}. Now, we bound each term on the right-hand side of \eqref{ineq:W_bound}.
From Lemma~\ref{lemma:discrete_sobolev} and the $H^1$-stability in Proposition~\eqref{prop:stability}, we have
\begin{equation}\label{eq:bound_U_Y}
    \norm{U+Y}_\infty \leq C_L \norm{U+Y}_{H_\dx^1} \leq C_L\Big(\norm{U}_{H_\dx^1} + \norm{Y}_{H_\dx^1}\Big) \leq C_1.
\end{equation}
On the other hand, from the definition of $N_X$ and $N_Y$ and Lemma~\ref{lemma:Dx_Dxx_bounds}, we obtain
\begin{equation*}
    \normpp{N_X - N_Y}{2} \leq \Bigg( 1 + \frac{4|\ka|}{\dx^2} \Bigg) \normpp{Q_X - Q_Y}{2}.
\end{equation*}
Notice that $Q_X - Q_Y = \big(|X|^2 - |Y|^2\big)/2 = \Re \Big( \conj{(X+Y)}\odot W \Big)/2$. Hence,
\begin{equation}\label{ineq:bound_NX_NY}
    \normpp{N_X - N_Y}{2} \leq \frac{1}{2}\Bigg( 1 + \frac{4|\ka|}{\dx^2} \Bigg) \norm{X + Y}_\infty \normpp{W}{2} \leq \frac{C_2}{2}\Bigg( 1 + \frac{4|\ka|}{\dx^2} \Bigg) \normpp{W}{2}, 
\end{equation}
where we repeated the same process as in \eqref{eq:bound_U_Y} but for $X$ instead of $U$ in the last inequality.
Putting togheter \eqref{eq:bound_U_Y} and \eqref{ineq:bound_NX_NY} into \eqref{ineq:W_bound} give
\begin{equation}\label{ineq:bound_W_2}
    \normpp{W}{2}^2 \leq M \Bigg( \dt + \frac{4\dt|\ka|}{\dx^2}\Bigg) \normpp{W}{2}^2,
\end{equation}
where $M$ contains the constants $C_1$ and $C_2$ and does not depend on $U, X$ nor $Y$. Therefore, choosing this $M$ in \eqref{eq:cn-uniqueness-condition}, the inequality \eqref{ineq:bound_W_2} is only possible if $\normpp{W}{2} = 0$, which gives $X=Y$, concluding the proof.
\end{proof}

\section{Consistency}\label{sec:consistency}
From the local well-posedness result in Section~\ref{appx:well-posedness}, we can take $\psi_0\in H^{p}$, with $p=14$, satisfying $\boD_L(\psi_0)>0$ and the compatibility conditions \eqref{comp_cond_init} to obtain a solution $\psi\in C^4([0,T]; H^5(\boJ))$ of \eqref{QGP}, which satisfies $\boD_L(\psi( t))>0$ and the compatibility conditions \eqref{comp_cond_init} at the boundary, for all $t\in [0, T)$. This choice of $\psi$ gives us sufficient regularity to prove our consistency result.

Before going to the proof, we introduce the following two Taylor estimations well-suited for our purposes:
\begin{lemma}[Time average]\label{lemma:time_av}
    Let $\dt>0$ and $t\in \R$. For $f: [t, t+\dt]\to\C$ of class $C^2$, we have
    \begin{equation}\label{eq:avg}
        \frac{f(t) + f(t+\dt)}{2} = f(t+\dt/2) + \boT^{\,\dt}_f(t), \quad \boT^{\,\dt}_f(t)= \frac{1}{2}\int_{t}^{t+\dt}\left(\dt/2 - \abs{s - (t+\dt/2)}\right) f''(s) ds.
    \end{equation}
Moreover, the remainder $\boT^{\,\dt}_f(t)$ satisfies 
\begin{equation}\label{eq:av_bound}
    |\boT^{\,\dt}_f(t)| \leq \norm{f''}_\infty \cdot \frac{\dt^2}{8}.
\end{equation}
\end{lemma}
\begin{proof}
Performing the following Taylor expansions:
    \begin{align*}
    &f(t) = f(t+\dt/2-\dt/2) = f(t+\dt/2) - \frac{\dt}{2}f'(t+\dt/2) + \int_{t}^{t+\dt/2}\left(s-t\right)f''(s)ds,\\
    &f(t+\dt) = f(t+\dt/2+\dt/2) = f(t+\dt/2) + \frac{\dt}{2}f'(t+\dt/2) + \int_{t+\dt/2}^{t+\dt}\left(t+\dt - s\right)f''(s)ds,
\end{align*}
we obtain equation \eqref{eq:avg} by averaging both terms. The proof of \eqref{eq:av_bound} is straightforward.
\end{proof}
\begin{lemma}[Second-order finite difference approximation]\label{lemma:second_order}
    Let $\dx>0$ and $x\in \R$. For $g:[x-h, x+h]\to \C$ of class $C^4$, we have
    \begin{equation}\label{eq:sec_order_app}
        \Dxx{g} = \frac{g(x-\dx) - 2g(x) + g(x+\dx)}{\dx^2} = g''(x) + \boS^\dx_g(x), \quad \boS^\dx_g(x) = \int_{x-\dx}^{x+\dx}\frac{g^{(4)}(s)}{6\dx^2} (\dx - |s-x|)^3ds.
    \end{equation}
Moreover, the remainder $\boS^\dx_g(x)$ satisfies 
\begin{equation}\label{eq:sec_order_bound}
    \abs{\boS^\dx_g(x)} \leq \norm{g^{(4)}}_\infty \cdot \frac{\dx^2}{12}.
\end{equation}
\end{lemma}
\begin{proof}
    Performing the following Taylor expansions:
    \begin{align*}
    &g(x-\dx) = g(x) - \dx g'(x) + \frac{\dx^2}{2}g''(x) - \frac{h^3}{6}g'''(x) + \frac{1}{6}\int_{x-\dx}^{x}(\dx+s-x)^3g^{(4)}(s),\\
    &g(x+\dx) = g(x) + \dx g'(x) + \frac{\dx^2}{2}g''(x) + \frac{h^3}{6}g'''(x) + \frac{1}{6}\int_{x}^{x+\dx}(\dx+x-s)^3g^{(4)}(s),
\end{align*}
we obtain equation \eqref{eq:sec_order_app} by summing both terms, subtracting $2g(x)$ and dividing by $\dx^2$. The proof of \eqref{eq:sec_order_bound} is straightforward.
\end{proof}

We are now in position to prove the second-order consistency of the scheme:
\begin{proposition}[Second-order consistency]\label{prop:second_order}
Let $\psi\in C^4([0,T]; H^5(\boJ))$.
For $j\in \dbra{0}{\J}$, $n\in \dbra{0}{N-1}$ and $\dx, \dt <1$, define 
the local truncation error
\begin{equation}\label{eq:main_remainder}
    \mathcal{R}^{n+\frac{1}{2}}_j=\frac{i\left(\psi_j^{n} -\psi_j^{n+1}\right)}{\dt} + \big(\Dxx
    \psi^{n+\frac{1}{2}}\big)_j + \gs \psi_j^{n+\frac{1}{2}} q^{n+\frac{1}{2}}_j + \ka \psi_j^{n+\frac{1}{2}} \big(\Dxx q^{n+\frac{1}{2}}\big)_j.
\end{equation}
where $\psi_j^{n+\frac{1}{2}} = \big(\psi_j^{n} + \psi_j^{n+1}\big)/2$, $q_j^{n+\frac{1}{2}} = \big(q_j^n + q_j^{n+1}\big)/2$. Then,
\begin{equation}\label{eq:residual-L2-bounds}
    \max_{0\leq n\leq N-1}\norm{\mathcal{R}^{n+\frac{1}{2}}}_{2,\dx^+}
    \leq C(\dt^2+\dx^2),
\end{equation}
and
\begin{equation}\label{eq:residual-time-increment}
    \max_{1\leq n\leq N-1}
    \norm{\mathcal{R}^{n+\frac{1}{2}}-\mathcal{R}^{n-\frac{1}{2}}}_{2,\dx^+}
    \leq C\dt(\dt^2+\dx^2),
\end{equation}
where $C$ is a constant independent of $\dx$ and $\dt$, and $\mathcal{R}^{n-\frac{1}{2}} =
\mathcal{R}^{{n-1}+\frac{1}{2}}$ corresponds to the local truncation error at the previous time step.
\end{proposition}
\begin{proof}
We start by proving \eqref{eq:residual-L2-bounds}. The proof of \eqref{eq:residual-time-increment} is carried out after that.

\paragraph{Proof of \eqref{eq:residual-L2-bounds}.}
We first estimate $|\mathcal{R}^{n+\frac{1}{2}}_j|$ in the interior points $j\in\dbra{1}{\J-1}$, by using the well-known Taylor expansions in time and space, and then, we extend the estimates to the boundary points $j=0$ and $j=\J$.

Let us fix $n\in\dbra{0}{\Nt-1}$ and $j\in\dbra{1}{\J-1}$. We refer to the first two terms on the right-hand side of \eqref{eq:main_remainder} as linear terms, and the rest as nonlinear terms. 
\paragraph{Linear terms.} For fixed $x\in \R$, noting that $\psi(x, t_n) = \psi(x, \tnmed-\dt/2)$ and $\psi(x, t_{n+1}) = \psi(x, \tnmed+\dt/2)$, we perform the following Taylor expansions:
$$
    \psi(x, t_n) = \psi(x, \tnmed) - \frac{\dt}{2}\ptl_{t}\psi(x, \tnmed) + \frac{\dt^2}{8}\ptl_{tt}\psi(x, \tnmed) + \int_{t_n}^{\tnmed} (\tnmed - \frac{\dt}{2} - s)^2\frac{\ptl_{ttt}\psi(x_j, s)}{2}ds,
$$
$$
    \psi(x, t_{n+1}) = \psi(x, \tnmed) + \frac{\dt}{2}\ptl_{t}\psi(x, \tnmed) +\frac{\dt^2}{8}\ptl_{tt}\psi(x, \tnmed) + \int_{\tnmed}^{t_{n+1}} (\tnmed + \frac{\dt}{2} - s)^2\frac{\ptl_{ttt}\psi(x_j, s)}{2}ds.
$$
Taking $x=x_j$ and subtracting these two expressions, we obtain
\begin{equation*}
    \frac{\psi_j^{n} -\psi_j^{n+1}}{\dt} = -\psi_t(x_j, \tnmed) - R_1(x_j, t_n),
\end{equation*}
where the remainder term $R_1$ is defined by
\begin{equation*}
    R_1(x_j, t_n) = \frac{1}{2\dt}\int_{t_n}^{t_{n+1}} \left(\frac{\dt}{2} - |s - \tnmed|\right)^2\ptl_{ttt}\psi(x_j, s)ds.
\end{equation*}
For simplicity of notation, we ignore the dependence on $x_j$ and $t_n$ in $R_1$ and the same convention will be used in the estimates that follow. It is easy to check that
\begin{equation}\label{eq:R_time}
    |R_1| \leq C \infnorm{\ptl_{ttt}\psi} \dt^2.
\end{equation}

We now estimate the second linear term in \eqref{eq:main_remainder}. Applying Lemma~\ref{lemma:second_order} to the functions $\psi( t_n)$ and $\psi( t_{n+1})$, we get
\begin{align*}
    \Big(\Dxx&\psi^{n+\frac{1}{2}}\Big)_j = \frac{(\Dxx\psi^n)_j + (\Dxx\psi^{n+1})_j}{2} =\frac{\ptl_{xx}\psi(x_j, t_n) + \ptl_{xx}\psi(x_j, t_{n+1})}{2}+\frac{\boS_{\psi( t_n)}^{\dx}+\boS_{\psi( t_{n+1})}^\dx}{2} \\& =\frac{\ptl_{xx}\psi(x_j, t_n) + \ptl_{xx}\psi(x_j, t_{n+1})}{2}+\int_{x_{j-1}}^{x_{j+1}}\frac{(h-\abs{s-x})^3}{6\dx^2} \frac{\ptl_{xxxx}\psi(s, t_n) + \ptl_{xxxx}\psi(s, t_{n+1})}{2} ds.
\end{align*}
Then, applying Lemma~\ref{lemma:time_av} to the functions $\ptl_{xx}\psi(x_j, \cdot)$ and $\ptl_{xxxx}\psi(s, \cdot)$, we have
\begin{equation}\label{eq:second_order_ap}
    \left(\Dxx \psi^{n+\frac{1}{2}}\right)_j =\ptl_{xx}\psi(x_j, \tnmed) + 
     R_2,
\end{equation}
where the remainder term $R_2$ is given by
\begin{equation*}
\begin{aligned}
    R_2=&\boT_{\ptl_{xx}\psi(x_j, \cdot)}^\dx + \boS^\dx_{\psi( \tnmed)}\\ &+\frac{1}{6\dx^2}\int_{x_{j-1}}^{x_{j+1}}\left(h-|s-x_j|\right)^3\left[\int_{t_n}^{t_{n+1}} \left(\abs{\tilde{s}-t_{n+1}} + \dt/2\right)\ptl_{xxxxtt}\psi(s, \tilde{s})d\tilde{s}\right]ds.
\end{aligned}
\end{equation*}
Using the estimates in \eqref{eq:av_bound} and \eqref{eq:sec_order_bound}, we conclude that
\begin{equation}\label{eq:rem_psi_ineq}
\begin{aligned}
    |R_2| &\leq C\left[\norm{\ptl_{xxtt}\psi}_\infty\dt^2  + \norm{\ptl_{xxxx}\psi}_\infty\dx^2 +\norm{\ptl_{xxxxtt}\psi}_\infty\dt^2\dx^2\right] \\ &\leq C \left[\left( \norm{\ptl_{xxtt}\psi}_\infty + \norm{\ptl_{xxxxtt}\psi}_\infty\right)\dt^2  +\norm{\ptl_{xxxx}\psi}_\infty\dx^2\right],
\end{aligned}
\end{equation}
where we used that $\dx<1$. This concludes the estimates for the linear terms. Now, we turn to the nonlinear terms.
\paragraph{Nonlinear terms.} Applying Lemma \ref{lemma:time_av} to $\psi$ and $q$, we obtain the estimates
\begin{equation*}
    \psi_j^{n+\frac{1}{2}} = \frac{\psi^n + \psi^{n+1}}{2} = \psi(x, \tnmed) + \boT_{\psi(x, \cdot)}^{\,\dt},\quad q_j^{n+\frac{1}{2}} = q(x, \tnmed) +\boT_{q(x, \cdot)}^{\,\dt},
\end{equation*}
where the remainders are defined as in \eqref{eq:avg}. It follows that
\begin{equation*}
    \psi_j^{n+\frac{1}{2}} q_j^{n+\frac{1}{2}}= \psi(x_j, \tnmed) q(x_j, \tnmed) + R_3,
\end{equation*}
where
\begin{equation*}
    R_3 = \psi(x_j, \tnmed) \boT_{q(x_j, \cdot)}^{\,\dt}+ q(x_j, \tnmed) \boT_{\psi(x_j, \cdot)}^{\,\dt} + \boT_{q(x_j, \cdot)}^{\,\dt} \boT_{\psi(x_j, \cdot)}^{\,\dt}.
\end{equation*}
Using the estimate in \eqref{eq:av_bound}, we conclude that
\begin{equation}\label{eq:boQ}
    \abs{R_3} \leq C \left(\norm{\psi}_\infty \norm{\ptl_{tt}q}_\infty +  \norm{\ptl_{tt}\psi}_\infty \norm{q}_\infty + \norm{\ptl_{tt}\psi}_\infty \norm{\ptl_{tt}q}_\infty \right)\dt^2.
\end{equation}
Now, repeating the argument used to derive equation \eqref{eq:second_order_ap} but for $q$ instead of $\psi$, we have
\begin{equation*}
    \left(\Dxx q^{n+\frac{1}{2}}\right)_j =\ptl_{xx}q(x_j, \tnmed) + \tilde{R},
\end{equation*}
where the remainder term $\tilde{R}$ is given by
\begin{align*}
    \tilde{R}=&\boT_{\ptl_{xx}q(x_j, \cdot)}^\dx + \boS^\dx_{q(\cdot, \tnmed)}\\ &+\frac{1}{6\dx^2}\int_{x_{j-1}}^{x_{j+1}}\left(h-|s-x_j|\right)^3\left[\int_{t_n}^{t_{n+1}} \left(\abs{\tilde{s}-t_{n+1}} + \dt/2\right)\ptl_{xxxxtt}q(s, \tilde{s})d\tilde{s}\right]ds.
\end{align*}
Multiplying by $\psi_j^{n+\frac{1}{2}}$ 
, we get
\begin{align*}
    \psi_j^{n+\frac{1}{2}} \Big(\Dxx q^{n+\frac{1}{2}}\Big)_j= \psi(x_j, \tnmed) \ptl_{xx}q(x_j, \tnmed) + R_4,
\end{align*}
where
\begin{equation*}
    R_4 = \psi(x_j, \tnmed) \tilde{R}+ \ptl_{xx}q(x_j, \tnmed) \boT^{\,\dt}_{\psi(x_j, \cdot)} + \tilde{R} \boT^{\,\dt}_{\psi(x_j, \cdot)}.
\end{equation*}

Using the estimate in \eqref{eq:av_bound} and applying the same argument as in \eqref{eq:rem_psi_ineq} but for $q$, we conclude that 
\begin{equation}\label{eq:Q_tilde} 
\begin{aligned}
    \abs{R_4} &
    \begin{aligned}[t]
    \leq  C\Big( &\left[\infnorm{\psi}\infnorm{\ptl_{xxtt}q} + \infnorm{\psi}\infnorm{\ptl_{xxxxtt}q}+\infnorm{\ptl_{xx}\abssqr{\psi}}\infnorm{\ptl_{tt}\psi}\right]\dt^2 \\&+ \infnorm{\psi}\infnorm{\ptl_{xxxx}q}\dx^2 + \infnorm{\ptl_{tt}\psi}\infnorm{\ptl_{xxxx}q}\dx^2\dt^2 \\&+ \left[\infnorm{\ptl_{tt}\psi}\infnorm{\ptl_{xxtt}q} + \infnorm{\ptl_{tt}\psi}\infnorm{\ptl_{xxxxtt}q}\right]\dt^4 \Big)
    \end{aligned}
    \\&
    \begin{aligned}[t]
    \leq  C\Big( &\Big[\infnorm{\psi}\infnorm{\ptl_{xxtt}q} + \infnorm{\psi}\infnorm{\ptl_{xxxxtt}q}+\infnorm{\ptl_{xx}\abssqr{\psi}}\infnorm{\psi_{tt}}\\&+\infnorm{\ptl_{tt}\psi}\infnorm{\ptl_{xxtt}q} + \infnorm{\ptl_{tt}\psi}\infnorm{\ptl_{xxxxtt}q}\Big]\dt^2 \\&+ \left[\infnorm{\psi}\infnorm{\ptl_{xxxx}q} + \infnorm{\ptl_{tt}\psi}\infnorm{\ptl_{xxxx}q}\right]\dx^2 \Big),
    \end{aligned}
\end{aligned}
\end{equation}
where we used that $\dt<1$. Putting together: \eqref{eq:R_time}, \eqref{eq:rem_psi_ineq}, \eqref{eq:boQ}, \eqref{eq:Q_tilde} and using that $\psi$ is a solution of \eqref{QGP}, we obtain
\begin{align*}
    \max_{1\leq j\leq J-1}|\mathcal{R}^{n+\frac{1}{2}}_j| & \leq |R_1| + |R_2| + |R_3| + |R_4|
    \\&\begin{aligned}[t]
    \leq &C\Big( \Big[\infnorm{\ptl_{ttt}\psi} + \infnorm{\ptl_{xxtt}\psi} + \infnorm{\ptl_{xxxxtt}\psi} + \infnorm{\psi} \infnorm{\ptl_{tt}q} +  \infnorm{\ptl_{tt}\psi} \infnorm{q} \\&+ \infnorm{\ptl_{tt}\psi} \infnorm{\ptl_{tt}q} + \infnorm{\psi}\infnorm{\ptl_{xxtt}q} + \infnorm{\psi}\infnorm{\ptl_{xxxxtt}q}\\&+\infnorm{\ptl_{xx}q}\infnorm{\ptl_{tt}\psi}+\infnorm{\ptl_{tt}\psi}\infnorm{\ptl_{xxtt}q} + \infnorm{\ptl_{tt}\psi}\infnorm{\ptl_{xxxxtt}q}\Big] \dt^2\\&+ \left[ \infnorm{\ptl_{xxxx}\psi} + \infnorm{\psi}\infnorm{\ptl_{xxxx}q} + \infnorm{\ptl_{tt}\psi}\infnorm{\ptl_{xxxx}q}\right]\dx^2 \Big).
    \end{aligned}
    \end{align*}
Using the embedding $H^1(\mathcal{J})\hookrightarrow L^\infty(\mathcal{J})$ and the fact that $q = \abs{\psi}^2$, we can bound the right-hand side of the previous inequality in terms of Sobolev norms of $\psi$ only. This yields 
\begin{equation}\label{eq:final_R}
\begin{aligned}
    &\begin{aligned}[t]
    \max_{1\leq j\leq J-1}|\mathcal{R}^{n+\frac{1}{2}}_j| \leq &C\Big( \Big[\CHnorm{\psi}{1}{3} + \CHnorm{\psi}{3}{2} + \CHnorm{\psi}{5}{2} + \CHnorm{\psi}{1}{} \CHnorm{\psi}{1}{2}^2 \\&+  \CHnorm{\psi}{1}{2} \CHnorm{\psi}{1}{}^2 + \CHnorm{\psi}{1}{2} \CHnorm{\psi}{1}{2}^2 + \CHnorm{\psi}{1}{}\CHnorm{\psi}{3}{2}^2 \\&+ \CHnorm{\psi}{1}{}\CHnorm{\psi}{5}{2}^2+\CHnorm{\psi}{3}{}^2\CHnorm{\psi}{1}{2}+\CHnorm{\psi}{1}{2}\CHnorm{\psi}{2}{2}^2 \\&+ \CHnorm{\psi}{1}{2}\CHnorm{\psi}{5}{2}^2\Big] \dt^2+ \Big[ \CHnorm{\psi}{5}{} + \CHnorm{\psi}{1}{}\CHnorm{\psi}{5}{}^2\\&+ \CHnorm{\psi}{1}{2}\CHnorm{\psi}{5}{}^2\Big]\dx^2 \Big)
    \\& \leq C \CHnorm{\psi}{5}{3}\left(1 + \CHnorm{\psi}{5}{3}^2\right) (\dt^2+\dx^2).
    \end{aligned}
\end{aligned}
\end{equation}
We now estimate the values at the boundaries. For $j\in\{0, \J\}$, we perform similar arguments as in \eqref{eq:rem_psi_ineq} and \eqref{eq:Q_tilde}. For simplicity, we consider the estimations only for $j=0$, since the case $j=\J$ is symmetric.
To bound $\Dxx(\psi)_0$, we start doing a Taylor expansion on $\psi(x_1, t)$ to get
\begin{equation}\label{eq:taylor_boundary}
\begin{aligned}
    \psi(x_1, t) = &\psi(x_0, t) + \dx \ptl_{x}\psi(x_0, t)  +  \frac{\dx^2}{2}\ptl_{xx}\psi(x_0, t) + \frac{\dx^3}{6}\ptl_{xxx}\psi(x_0, t) \\ &+ \frac{1}{6}\int_{x_0}^{x_1} (h + x - s)^3\ptl_{xxxx}\psi(s, t) ds.
\end{aligned}
\end{equation}
Since the solution $\psi$ satisfies the homogeneous Neumann boundary condition in \eqref{QGP}, the second term in the right-hand side of \eqref{eq:taylor_boundary} vanishes. Moreover, from the compatibility condition \eqref{comp_cond_init} discussed in the introduction, and since $\psi(t, \cdot)\in H^5(\mathcal{J})$, we also have $\psi_{xxx}(x_0, t) = 0$. Therefore, by the definition of $\Dxx(\cdot)$ at the boundary in \eqref{eq:spc_disc_ops}, we get
\begin{equation*}
    (\Dxx\psi)_0 = \ptl_{xx}\psi(x_0, t) + \frac{1}{3\dx^2}\int_{x_0}^{x_1} (h + x - s)^3\ptl_{xxxx}\psi(s, t) ds.
\end{equation*}
Proceeding as in \eqref{eq:final_R}, 
and recalling that $\dt<1$, we finally infer that
\begin{equation}\label{eq:bound_R0}
    |\boR_0^{n+\frac{1}{2}}|
    \leq C \CHnorm{\psi}{5}{3}\left(1 + \CHnorm{\psi}{5}{3}^2\right) (\dt^2+\dx^2).
\end{equation}
From the definition of the norm $\normpp{\cdot}{2}$ in \eqref{eq:norms} and using the estimates \eqref{eq:final_R} and \eqref{eq:bound_R0}, we conclude that
\begin{equation*}
    \norm{\mathcal{R}^{n+\frac{1}{2}}}_{2,\dx^+}^2 \leq 2L  \max_{0\leq j\leq \J}|\mathcal{R}^{n+\frac{1}{2}}_j|^2 \leq C \CHnorm{\psi}{5}{3}\left(1 + \CHnorm{\psi}{5}{3}^2\right) (\dt^2+\dx^2),
\end{equation*}
where we used that $L>0$ is fixed. This concludes the proof of the inequality \eqref{eq:residual-L2-bounds}.

\paragraph{Proof of \eqref{eq:residual-time-increment}.}
For \eqref{eq:residual-time-increment}, it is convenient to define the residual as a continuous function in time.
We first verify that the boundary argument used to prove \eqref{eq:residual-L2-bounds} remains valid for the estimates of the time derivatives appearing below. Let us fix $j\in\{0,\J\}$. Since $\psi$ satisfies \eqref{QGP}, we have $\ptl_x\psi(x_j,t)=0$ for every $t\in\mathcal{I}$. Moreover, from the compatibility condition \eqref{comp_cond_init}
, we also have $\ptl_{xxx}\psi(x_j,t)=0$ for every $t\in\mathcal{I}$. Differentiating these identities with respect to time, we obtain
\begin{equation}\label{eq:boundary_time_psi}
    \ptl_{xt}\psi(x_j,t)=0\quad \text{ and }\quad \ptl_{xxxt}\psi(x_j,t)=0,
\end{equation}
and the same conclusion holds for $q=|\psi|^2$ using that
\begin{equation}\label{eq:boundary_density_derivatives}
    \ptl_x q=2\Re\big(\conj{\psi}\ptl_x\psi\big) \quad \text{ and } \quad \ptl_{xxx}q=2\Re\big(\conj{\psi}\ptl_{xxx}\psi+3\conj{\ptl_x\psi}\ptl_{xx}\psi\big),
\end{equation}
and $H^5(\mathcal{J})$ is an algebra, which gives the time regularity required to differentiate these identities. Hence, using the boundary identities for $\psi$ in \eqref{eq:boundary_density_derivatives}, we obtain
\begin{equation}\label{eq:boundary_time_q}
    \ptl_{xt}q(x_j,t)=0\quad \text{ and }\quad \ptl_{xxxt}q(x_j,t)=0.
\end{equation}
By \eqref{eq:boundary_time_psi} and \eqref{eq:boundary_time_q}, when applying the Taylor expansion \eqref{eq:taylor_boundary} to $\ptl_t\psi$ and $\ptl_tq$, the first- and third-order spatial terms vanish, exactly as for $\psi$ and $q$ in the derivation of \eqref{eq:bound_R0}. Thus, the following spatial estimates also hold at the boundary nodes by the same argument

For $t\in I_\dt := [\dt/2, T-\dt/2]$, we define $\psi_\pm(t) = \psi(\cdot_j, t\pm \dt/2)$,  $q_\pm(t) = |\psi_\pm(t)|^2$, where $\cdot_j$ denotes the evaluation of the space grid;  and the time discrete averages 
$$\psi_\dt = \frac{\psi_+(t) + \psi_-(t)}{2} \quad \text{ and } \quad q_\dt = \frac{q_+(t) + q_-(t)}{2}.$$
Then, we define the residual as a function of time by
\begin{equation*}
    \mathfrak{R}_{\dx, \dt}(t) = i\frac{\psi_-(t) - \psi_+(t)}{\dt} + \Dxx \psi_\dt(t) + \gs q_\dt(t)\odot \psi_\dt(t) + \ka \big(\Dxx q_\dt(t)\big) \odot \psi_\dt(t),
\end{equation*}
which takes values in $\C^{\J+1}$ and is defined in $I_\dt$.
Notice that, by construction, $\mathfrak{R}_{\dx, \dt}\big(\tnmed\big) = \mathcal{R}^{n+\frac{1}{2}}$ and $\mathfrak{R}_{\dx, \dt}\big(t_{n-\frac{1}{2}}\big) = \mathcal{R}^{n-\frac{1}{2}}$. Moreover, since $\psi\in \mathcal{C}^4\big(\mathcal{I}, H^5(\mathcal{J})\big)$, the function $\mathfrak{R}_{\dx, \dt}$ is continuously differentiable in time. Therefore, by the fundamental theorem of calculus, we have
\begin{equation}\label{ineq:residual-time-increment}
    \normpp{\mathcal{R}^{n+\frac{1}{2}} - \mathcal{R}^{n-\frac{1}{2}}}{2} \leq \normpp{\int_{t_{n-\frac{1}{2}}}^{t_{n+\frac{1}{2}}} \ptl_t \mathfrak{R}_{\dx, \dt}(t) dt}{2} \leq \dt \sup_{t\in I_\dt} \normpp{\ptl_t \mathfrak{R}_{\dx, \dt}(t)}{2}.
\end{equation}
It remains to prove that
$$\sup_{t\in I_\dt} \normpp{\ptl_t \mathfrak{R}_{\dx, \dt}(t)}{2} \leq C(\dx^2 + \dt^2).$$
Indeed, taking the time derivative of $\mathfrak{R}_{\dx, \dt}$ yields
\begin{equation}\label{eq:D-time-derivative-increment}
\begin{aligned}
\ptl_t\mathfrak{R}_{\dx,\dt}(t)={}&i\frac{\ptl_t\psi_-(t)-\ptl_t\psi_+(t)}{\dt}+\Dxx\ptl_t\psi_\dt(t) +\gs\big(\ptl_tq_\dt(t)\odot\psi_\dt(t)+q_\dt(t)\odot\ptl_t\psi_\dt(t)\big)\\
&+\ka\big((\Dxx\ptl_tq_\dt(t))\odot\psi_\dt(t)+(\Dxx q_\dt(t))\odot\ptl_t\psi_\dt(t)\big).
\end{aligned}
\end{equation}
 We refer to the first two terms on the right-hand side of \eqref{eq:D-time-derivative-increment} as linear terms and the remaining ones as nonlinear terms. As in the proof of \eqref{eq:residual-L2-bounds}, we estimate the linear and nonlinear terms in \eqref{eq:D-time-derivative-increment} separately.
\paragraph{Linear terms.} The centered Taylor expansion used in \eqref{eq:R_time}, applied to $\ptl_t\psi$ on $[t-\dt/2,t+\dt/2]$, gives
\begin{equation}\label{ineq:time_app_D}
\Big \Vert i\frac{\ptl_t\psi_- - \ptl_t\psi_+}{\dt}+i\ptl_{tt}\psi \Big \Vert_\infty\leq C\infnorm{\ptl_{tttt}\psi}\dt^2.
\end{equation}
For the second term, we repeat the argument leading to \eqref{eq:rem_psi_ineq}, with $\ptl_t\psi$ in place of $\psi$. Using the temporal estimate \eqref{eq:av_bound} and the spatial estimate \eqref{eq:sec_order_bound}, we obtain
\begin{equation}\label{ineq:D-linear-space}
\begin{aligned}
\norm{\Dxx\ptl_t\psi_\dt-\ptl_{xxt}\psi}_\infty&\leq C\big(\infnorm{\ptl_{xxttt}\psi}\dt^2+\infnorm{\ptl_{xxxxt}\psi}\dx^2+\infnorm{\ptl_{xxxxttt}\psi}\dt^2\dx^2\big)\\
&\leq C\big[\big(\infnorm{\ptl_{xxttt}\psi}+\infnorm{\ptl_{xxxxttt}\psi}\big)\dt^2+\infnorm{\ptl_{xxxxt}\psi}\dx^2\big],
\end{aligned}
\end{equation}
where we used that $\dx<1$. The spatial estimates used to derive \eqref{eq:rem_psi_ineq} and \eqref{eq:bound_R0} apply here to \(\partial_t\psi\) in place of \(\psi\), at the interior and boundary nodes, respectively.

\paragraph{Nonlinear terms.} We now estimate separately the remaining four nonlinear terms in \eqref{eq:D-time-derivative-increment}. First, applying the Lemma~\ref{lemma:time_av}, in particular, the inequality \eqref{eq:av_bound} to $\psi$, $\ptl_t\psi$, $q$ and $\ptl_tq$, respectively, gives
\begin{equation}\label{ineq:D-average-psi-t}
\norm{\psi_\dt-\psi}_\infty \leq C\infnorm{\ptl_{tt}\psi}\dt^2,\qquad 
\norm{\ptl_t\psi_\dt-\ptl_t\psi}_\infty\leq C\infnorm{\ptl_{ttt}\psi}\dt^2,
\end{equation}
and
\begin{equation}\label{ineq:D-average-q-t}
\norm{q_\dt-q}_\infty \leq C\infnorm{\ptl_{tt}\abssqr{\psi}}\dt^2, \qquad \norm{\ptl_tq_\dt-\ptl_tq}_\infty \leq C\infnorm{\ptl_{ttt}\abssqr{\psi}}\dt^2.
\end{equation}
Moreover, by the definitions of the time averages, we have
\begin{equation}\label{ineq:D-average-bounds}
\norm{\psi_\dt}_\infty\leq\infnorm{\psi}\qquad \text{ and }\qquad \norm{\ptl_t\psi_\dt}_\infty\leq\infnorm{\ptl_t\psi}.
\end{equation}
We estimate the first cubic product $\ptl_tq_\dt\odot\psi_\dt-\ptl_tq\odot\psi$ by decomposing it as
\begin{equation*}
\ptl_tq_\dt\odot\psi_\dt-\ptl_tq\odot\psi=(\ptl_tq_\dt-\ptl_tq)\odot\psi_\dt+\ptl_tq\odot(\psi_\dt-\psi).
\end{equation*}
Using \eqref{ineq:D-average-psi-t}, \eqref{ineq:D-average-q-t} and \eqref{ineq:D-average-bounds}, we get
\begin{equation}\label{ineq:D-cubic-first}
\norm{\ptl_tq_\dt\odot\psi_\dt-\ptl_tq\odot\psi}_\infty\leq C\big(\infnorm{\psi}\infnorm{\ptl_{ttt}q}+\infnorm{\ptl_tq}\infnorm{\ptl_{tt}\psi}\big)\dt^2.
\end{equation}
Similarly, decomposing the second cubic product $q_\dt\odot\ptl_t\psi_\dt-q\odot\ptl_t\psi$ as
$$q_\dt\odot\ptl_t\psi_\dt-q\odot\ptl_t\psi = (q_\dt-q)\odot\ptl_t\psi_\dt+q\odot(\ptl_t\psi_\dt-\ptl_t\psi),$$
and using \eqref{ineq:D-average-psi-t}, \eqref{ineq:D-average-q-t} and \eqref{ineq:D-average-bounds}, we obtain
\begin{equation}\label{ineq:D-cubic-second}
\norm{q_\dt\odot\ptl_t\psi_\dt-q\odot\ptl_t\psi}_\infty\leq C\big(\infnorm{\ptl_t\psi}\infnorm{\ptl_{tt}q}+\infnorm{q}\infnorm{\ptl_{ttt}\psi}\big)\dt^2.
\end{equation}

For the quasilinear products, we follow \eqref{eq:Q_tilde} and first apply the argument of \eqref{eq:rem_psi_ineq} to $q$ and $\ptl_tq$, as in \eqref{ineq:D-linear-space}. The estimates \eqref{eq:av_bound} and \eqref{eq:sec_order_bound}, together with $\dx<1$, give separately
\begin{equation}\label{ineq:D-space-q}
\norm{\Dxx q_\dt-\ptl_{xx}q}_\infty\leq C\big[\big(\infnorm{\ptl_{xxtt}q}+\infnorm{\ptl_{xxxxtt}q}\big)\dt^2+\infnorm{\ptl_{xxxx}q}\dx^2\big]
\end{equation}
and
\begin{equation}\label{ineq:D-space-q-t}
\norm{\Dxx\ptl_tq_\dt-\ptl_{xxt}q}_\infty\leq C\big[\big(\infnorm{\ptl_{xxttt}q}+\infnorm{\ptl_{xxxxttt}q}\big)\dt^2+\infnorm{\ptl_{xxxxt}q}\dx^2\big].
\end{equation}
We then decompose the first quasilinear product $(\Dxx\ptl_tq_\dt)\odot\psi_\dt-\ptl_{xxt}q\odot\psi$ as
\begin{equation*}
(\Dxx\ptl_tq_\dt)\odot\psi_\dt-\ptl_{xxt}q\odot\psi=(\Dxx\ptl_tq_\dt-\ptl_{xxt}q)\odot\psi_\dt+\ptl_{xxt}q\odot(\psi_\dt-\psi).
\end{equation*}
Using again \eqref{ineq:D-space-q-t}, \eqref{ineq:D-average-psi-t} and \eqref{ineq:D-average-bounds}, we conclude that
\begin{equation}\label{ineq:D-quasilinear-first}
\begin{aligned}
\norm{(\Dxx\ptl_tq_\dt)\odot\psi_\dt-\ptl_{xxt}q\odot\psi}_\infty\leq C\Big(&\big[\infnorm{\psi}\big(\infnorm{\ptl_{xxttt}q}+\infnorm{\ptl_{xxxxttt}q}\big)\\
&+\infnorm{\ptl_{xxt}q}\infnorm{\ptl_{tt}\psi}\big]\dt^2+\infnorm{\psi}\infnorm{\ptl_{xxxxt}q}\dx^2\Big).
\end{aligned}
\end{equation}
For the second quasilinear product $(\Dxx q_\dt)\odot\ptl_t\psi_\dt-\ptl_{xx}q\odot\ptl_t\psi$, we use the decomposition
\begin{equation*}
(\Dxx q_\dt)\odot\ptl_t\psi_\dt-\ptl_{xx}q\odot\ptl_t\psi=(\Dxx q_\dt-\ptl_{xx}q)\odot\ptl_t\psi_\dt+\ptl_{xx}q\odot(\ptl_t\psi_\dt-\ptl_t\psi).
\end{equation*}
Combining \eqref{ineq:D-space-q}, \eqref{ineq:D-average-psi-t} and \eqref{ineq:D-average-bounds}, we get
\begin{equation}\label{ineq:D-quasilinear-second}
\begin{aligned}
\norm{(\Dxx q_\dt)\odot\ptl_t\psi_\dt-\ptl_{xx}q\odot\ptl_t\psi}_\infty\leq C\Big(&\big[\infnorm{\ptl_t\psi}\big(\infnorm{\ptl_{xxtt}q}+\infnorm{\ptl_{xxxxtt}q}\big)\\
&+\infnorm{\ptl_{xx}q}\infnorm{\ptl_{ttt}\psi}\big]\dt^2+\infnorm{\ptl_t\psi}\infnorm{\ptl_{xxxx}q}\dx^2\Big).
\end{aligned}
\end{equation}

We now combine the linear and nonlinear estimates. Since $\psi$ solves \eqref{QGP}, differentiating that equation with respect to time yields
\begin{equation}\label{eq:QGP-time-increment}
-i\ptl_{tt}\psi+\ptl_{xxt}\psi+\gs\big(\ptl_tq\,\psi+q\,\ptl_t\psi\big)+\ka\big(\ptl_{xxt}q\,\psi+\ptl_{xx}q\,\ptl_t\psi\big)=0.
\end{equation}
Subtracting \eqref{eq:QGP-time-increment} from \eqref{eq:D-time-derivative-increment}, and using \eqref{ineq:time_app_D}, \eqref{ineq:D-linear-space}, \eqref{ineq:D-cubic-first}, \eqref{ineq:D-cubic-second}, \eqref{ineq:D-quasilinear-first} and \eqref{ineq:D-quasilinear-second}, we obtain
\begin{equation*}
\begin{aligned}
\norm{\ptl_t\mathfrak{R}_{\dx,\dt}(t)}_\infty\leq C\Big(&\Big[\infnorm{\ptl_{tttt}\psi}+\infnorm{\ptl_{xxttt}\psi}+\infnorm{\ptl_{xxxxttt}\psi}\\
&+\infnorm{\ptl_t\psi}\infnorm{\ptl_{ttt}q}+\infnorm{\ptl_tq}\infnorm{\ptl_{tt}\psi}\\
&+\infnorm{\ptl_t\psi}\infnorm{\ptl_{tt}q}+\infnorm{\abssqr{\psi}}\infnorm{\ptl_{ttt}\psi}\\
&+\infnorm{\psi}\big(\infnorm{\ptl_{xxttt}q}+\infnorm{\ptl_{xxxxttt}q}\big)\\
&+\infnorm{\ptl_{xxt}q}\infnorm{\ptl_{tt}\psi}+\infnorm{\ptl_{xx}q}\infnorm{\ptl_{ttt}\psi}\\
&+\infnorm{\ptl_t\psi}\big(\infnorm{\ptl_{xxtt}q}+\infnorm{\ptl_{xxxxtt}q}\big)\Big]\dt^2\\
&+\Big[\infnorm{\ptl_{xxxxt}\psi}+\infnorm{\psi}\infnorm{\ptl_{xxxxt}q}+\infnorm{\ptl_t\psi}\infnorm{\ptl_{xxxx}q}\Big]\dx^2\Big),
\end{aligned}
\end{equation*}
where $C$ also contains the coefficients $\gs$ and $\ka$. As in \eqref{eq:final_R}, we use the embedding $H^1(\mathcal{J})\hookrightarrow L^\infty(\mathcal{J})$ and the algebra property of $H^m(\mathcal{J})$, $m\geq1$, to reduce these norms. 
Consequently, we infer that
\begin{equation}\label{ineq:D-global-Sobolev}
\begin{aligned}
\norm{\ptl_t\mathfrak{R}_{\dx,\dt}(t)}_\infty\leq C\Big(&\Big[\CHnorm{\psi}{1}{4}+\CHnorm{\psi}{3}{3}+\CHnorm{\psi}{5}{3}\\
&+\CHnorm{\psi}{1}{}\CHnorm{\psi}{1}{3}^2+\CHnorm{\psi}{1}{1}^2\CHnorm{\psi}{1}{2}\\
&+\CHnorm{\psi}{1}{1}\CHnorm{\psi}{1}{2}^2+\CHnorm{\psi}{1}{}^2\CHnorm{\psi}{1}{3}\\
&+\CHnorm{\psi}{1}{}\big(\CHnorm{\psi}{3}{3}^2+\CHnorm{\psi}{5}{3}^2\big)\\
&+\CHnorm{\psi}{3}{1}^2\CHnorm{\psi}{1}{2}+\CHnorm{\psi}{3}{}^2\CHnorm{\psi}{1}{3}\\
&+\CHnorm{\psi}{1}{1}\big(\CHnorm{\psi}{3}{2}^2+\CHnorm{\psi}{5}{2}^2\big)\Big]\dt^2\\
&+\Big[\CHnorm{\psi}{5}{1}+\CHnorm{\psi}{1}{}\CHnorm{\psi}{5}{1}^2+\CHnorm{\psi}{1}{1}\CHnorm{\psi}{5}{}^2\Big]\dx^2\Big)\\
\leq{}&C\CHnorm{\psi}{5}{4}\big(1+\CHnorm{\psi}{5}{4}^2\big)(\dt^2+\dx^2).
\end{aligned}
\end{equation}
Recalling the norm in \eqref{eq:norms}, we have $\normpp{U}{2}\leq\sqrt{2L}\norm{U}_\infty$. Since $L$ is fixed, \eqref{ineq:D-global-Sobolev} implies
\begin{equation}\label{ineq:D-time-derivative-L2}
\sup_{t\in I_\dt}\normpp{\ptl_t\mathfrak{R}_{\dx,\dt}(t)}{2}\leq C\CHnorm{\psi}{5}{4}\big(1+\CHnorm{\psi}{5}{4}^2\big)(\dt^2+\dx^2)\leq C(\dt^2+\dx^2),
\end{equation}
where the dependence on $\CHnorm{\psi}{5}{4}$ is absorbed into $C$. Finally, from \eqref{ineq:residual-time-increment} and \eqref{ineq:D-time-derivative-L2}, we conclude that
\begin{equation}\label{ineq:D-residual-increment-final}
\normpp{\mathcal{R}^{n+\frac12}-\mathcal{R}^{n-\frac12}}{2}\leq C\dt(\dt^2+\dx^2).
\end{equation}
Taking the maximum over $1\leq n\leq\Nt-1$ in \eqref{ineq:D-residual-increment-final} proves \eqref{eq:residual-time-increment}.
\end{proof}
\section{Convergence}\label{sec:convergence}
In this section, we provide the necessary tools to prove the second-order convergence in Theorem~\ref{thm:error_convergence}, in addition to the previous results. To do so, our proof will use the following discrete Gronwall inequality. A proof of a standard discrete version of the inequality can be found in \cite{clark_1987}.
\begin{lemma}[Discrete Gronwall inequality]\label{lem:discrete_gronwall}
Let $\dt>0$, $N\in\mathbb{N}$, and $t_n=n\dt$, with $T=N\dt$. Let $(a_n)_{n=0}^N$ and $(b_n)_{n=0}^{N-1}$ be nonnegative sequences such that, for some constants $A,B\geq 0$, with $B$ independent of $\dt$ and $N$ and $a_0 \leq A$,
\begin{equation*}
    a_m \leq A + \dt\sum_{n=0}^{m-1} \big( B a_n+b_n \big), \quad \text{ for all } \quad m \in \dbra{1}{\Nt}.
\end{equation*}
Then, the following estimate holds:
\begin{equation*}
    \max_{0\leq m\leq N}a_m \leq e^{BT} \Big( A+\dt\sum_{n=0}^{N-1}b_n\Big).
\end{equation*}
\end{lemma}
We recall the definition of the error $e^n=\Psi^n-\psi^n$, and we define the density error as $d^n=\Dns^n-\dns^n$. We also introduce the midpoint notation $e^{n+\frac{1}{2}}=(e^{n+1}+e^n)/2$ and $d^{n+\frac{1}{2}}=(d^{n+1}+d^n)/2$. Subtracting \eqref{eq:main_remainder} to the scheme \eqref{scheme}, we obtain the following equation for $e^n$:
\begin{equation}\label{eq:error-equation}
    \begin{aligned}
 i\frac{e^{n+1}-e^n}{\dt} = & \Dxx e^{n+\frac{1}{2}} +\gs\Big(\Dns^{n+\frac{1}{2}}\odot \Psi^{n+\frac{1}{2}}-\dns^{n+\frac{1}{2}}\odot\psi^{n+\frac{1}{2}}\Big)\\
 &+\ka\Big((\Dxx \Dns^{n+\frac{1}{2}})\odot\Psi^{n+\frac{1}{2}}-(\Dxx \dns^{n+\frac{1}{2}})\odot\psi^{n+\frac{1}{2}}\Big)+\boR^{n+\frac{1}{2}}.
\end{aligned}
\end{equation}

Both, $e^n$ and $d^n$ satisfy some properties:
\begin{lemma}\label{lemma:H1_dn_en}
    For $n\in \N$, the density error $d^n$ satisfies
    \begin{equation}\label{ineq:H1_dn_en}
        \norm{d^n}_{2, \dx^+}\leq C_1\norm{e^n}_{2, \dx^+}\quad \text{ and } \quad \norm{\Dxp d^n}_{2, \dx^-} \leq C_2\norm{e^n}_{H^1_\dx},
    \end{equation}
for $C_1, C_2$ constants independent of $\dx$ and $\dt$.
\end{lemma}
\begin{proof}
    Let $n\in \N$. Notice that $d_j^n = \ip{\Psi_j^n+\psi_j^n}{e_j^n}$, for $j\in \dbra{0}{\J}$. Therefore, $|d_j^n| \leq |\Psi_j^n+\psi_j^n||e_j^n|$, and
    \begin{equation*}
        \norm{d^n}_{2,h,+} \leq \norm{\Psi^n+\psi^n}_\infty \norm{e^n}_{2,h,+}.
    \end{equation*} 
    We conclude the first inequality in \eqref{ineq:H1_dn_en} by using the stability property in Proposition~\ref{prop:stability}.
    For the second one, we use Corollary~\ref{corollary:product_rule} with $U=\Psi^n+\psi^n$ and $V=e^n$ to obtain
    \begin{equation}\label{eq:dx_dn_1}
        (\Dxp d^n)_j = \ip{\Psi^n_{j+\frac{1}{2}} + \psi^n_{j+\frac{1}{2}}}{\Dxp e^n} + \ip{e^n_{j+\frac{1}{2}}}{\Dxp(\Psi^n +\psi^n)_j}.
    \end{equation}
    Hence, completing the weights of the norm $\normpm{\cdot}{2}$ defined in \eqref{eq:norms}, we have
    \begin{align*}
        \norm{\Dxp d^n}_{2,h,-} \leq &(\norm{\Psi^n}_{2,h,+} + \norm{\psi^n}_{2,h,+} )\norm{\Dxp e^n}_{2,h,+}
        \\&+ (\norm{\Dxp \Psi^n}_{2,h,-} + \norm{\Dxp \psi^n}_{2,h,-} ) \norm{e^n}_{2,h,+},
    \end{align*}
    and the proof follows again by using the stability property in Proposition~\ref{prop:stability}.
    \end{proof}

\paragraph{Discrete action.} 
For $U\in \C^{\J+1}$, we recall the definition of the \textit{discrete action} 
operator $\boS_{\dx,\mu}$ introduced in \eqref{eq:action_op:intro},
which will be denoted simply by  $\boS_{\dx}$ in the sequel.
Notice that it can be recast as
\begin{equation}\label{eq:action_op}
    \boS_\dx(U) = E_\dx(U) + \frac{\mu}{2} N_\dx(U) = \frac{1}{2}\normpm{\Dxp U}{2}^2 - \frac{\gs}{4}\normpp{U}{4}^4 + \frac{\ka}{4}\normpm{\Dxp|U|^2}{2}^2 + \frac{\mu}{2} \normpp{U}{2}^2,
\end{equation}
where $\mu$ is a real positive constant, and the definition of the remainder $\Fh$ is given in \eqref{eq:F_remainder}.
The terms $D\boS_\dx(U)[V]$ and $D^2\boS_\dx(U)[V, V]$ in \eqref{eq:action_op} correspond to the 
derivatives of $\boS_\dx$ with respect to $U$ in the direction $V$, given by
\begin{equation*}
    D\boS_\dx(U)[V] = -\ipRp{\Dxx U + \gs |U|^2 \odot U + \ka \Dxx |U|^2\odot U - \mu U}{V},
\end{equation*}
and
\begin{equation}\label{eq:D2F}
\begin{aligned}
    D^2\boS_\dx(U)[V, V] = & \normpm{\Dxp V}{2}^2 - \gs \Big(2\normpp{\ip{U}{V}}{2}^2 + \ipRp{|U|^2}{|V|^2}\Big) \\+& \ka \Big(2\normpm{\Dxp \ip{U}{V}}{2}^2 + \ipRm{\Dxp |U|^2}{\Dxp |V|^2}\Big) + \mu \normpp{V}{2}^2.
\end{aligned}
\end{equation}
We recall the definition of the remainder $\boH_{\dx, \mu}$ defined in \eqref{eq:F_remainder}. For simplicity, we will denote it by $\Fh$ in the sequel. The next proposition establishes an upper bound estimate for $\Fh$:
\begin{proposition}\label{prop:F_remainder_estimate}
    Let $U, V\in \C^{\J+1}$ satisfying $\norm{U}_{H^1_\dx} \leq M_1$ and $\norm{U+V}_{H^1_\dx} \leq M_2$ with $M_1, M_2>0$ constants independent of $\dx$. Then
    \begin{equation*}
        |\Fh(U, V)| \leq C \norm{V}^2_{H^1_\dx},
    \end{equation*}
    where $C$ is a constant independent of $\dx$.
\end{proposition}

\begin{proof}
    It suffices to show that $|D^2\boS_\dx(U+\theta V)[V, V]| \leq C \norm{V}^2_{H^1_\dx}$ uniformly for $0\leq \theta \leq 1$, with $C$ independent of $\dx$. Let $\theta \in [0, 1]$ and define $W = U+\theta V = (1 - \theta) U + \theta(U + V)$. Then, $W$ satisfies the following bound:
    \begin{equation}\label{ineq:W_estimate}
        \norm{W}_{H^1_\dx} \leq (1 - \theta)\norm{U}_{H^1_\dx} + \theta \norm{U+V}_{H^1_\dx} \leq M,
    \end{equation}
    where $M$ contains the constants $M_1$ and $M_2$. Now, from \eqref{eq:D2F}, it is easy to check that
    \begin{equation}\label{eq:D2F_estimate_1}
    \begin{aligned}
        |D^2\boS_\dx(W)[V, V]| \leq & C_1 \norm{V}_{H^1_\dx}^2 + 2\normpp{\ip{W}{V}}{2}^2 + |\ipRp{|W|^2}{|V|^2}|\\
        + & |\ka| \Big(2\normpm{\Dxp \ip{W}{V}}{2}^2 + |\ipRm{\Dxp |W|^2}{\Dxp |V|^2}|\Big)
    \end{aligned}
    \end{equation}
    where $C_1$ contains the term $\mu$. We now estimate each term on the right-hand side of \eqref{eq:D2F_estimate_1}. First, we have $\ip{W_j}{V_j}\leq |W_j||V_j| \leq \norm{W}_\infty |V_j|$, for $j\in\dbra{0}{\J}$, and using Lemma~\ref{lemma:discrete_sobolev} and \eqref{ineq:W_estimate}, we obtain
    \begin{equation}\label{ineq:RF_term_1}
        2\normpp{\ip{W}{V}}{2}^2 \leq 2 \norm{W}_{H^1_\dx}^2 \normpp{V}{2}^2 \leq C_2 \norm{V}_{H^1_\dx}^2,
    \end{equation}
    where the constant $C_2$ contains the constant $M$. Similarly, using the definition of $\ipRp{\cdot}{\cdot}$ in \eqref{eq:inner}, we have
    \begin{equation}\label{ineq:RF_term_2}
        |\ipRp{|W|^2}{|V|^2}| \leq C_3 \norm{V}_{H^1_\dx}^2,
    \end{equation}
    where $C_3$ is a constant containing $M$. Passing to the fourth term on the right-hand side of \eqref{eq:D2F_estimate_1}, from the discrete product rule in Corollary~\ref{corollary:product_rule} and the discrete Sobolev inequality in Lemma~\ref{lemma:discrete_sobolev}, we have
    \begin{equation}\label{ineq:RF_term_3}
        \normpm{\Dxp \ip{W}{V}}{2}^2 \leq \norm{W}_\infty \normpm{\Dxp V}{2} + \norm{V}_\infty \normpm{\Dxp W}{2} \leq C_4 \norm{V}_{H^1_\dx}^2,
    \end{equation}
    where the constant $C_4$ contains the terms $M$ and the constants from the discrete Sobolev inequality. Finally, for the last term on the right-hand side of \eqref{eq:D2F_estimate_1}, we use the Cauchy--Schwarz inequality to get
    \begin{equation}\label{charmander}
        |\ipRm{\Dxp |W|^2}{\Dxp |V|^2}| \leq \normpm{\Dxp |W|^2}{2}\normpm{\Dxp |V|^2}{2}.
    \end{equation}
    From Corollary~\ref{corollary:product_rule}, we have $(\Dxp |W|^2)_j = 2\ip{W_{j+\frac{1}{2}}}{(\Dxp W)_j}$, and similarly for $(\Dxp |V|^2)_j$. Hence, using the discrete Sobolev inequality in Lemma~\ref{lemma:discrete_sobolev} again, we have
    \begin{equation}\label{ineq:RF_term_4}
        |\ipRm{\Dxp |W|^2}{\Dxp |V|^2}| \leq 4 \norm{W}_\infty \normpm{\Dxp W}{2} \norm{V}_\infty \normpm{\Dxp V}{2} \leq C_5 \norm{V}_{H^1_\dx}^2,
    \end{equation}
    where $C_5$ contains the constant $M$ and the ones from the discrete Sobolev inequality. From the definition of $\boH_\dx$ in \eqref{eq:F_remainder} and putting together \eqref{ineq:RF_term_1}--\eqref{ineq:RF_term_4}\silentref{ineq:RF_term_2}\silentref{ineq:RF_term_3}\silentref{charmander}
    in \eqref{eq:D2F_estimate_1}, we conclude that
    \begin{equation*}
        |\Fh(U, V)| \leq \int_0^1 |(1-\theta)| |D^2\boS_\dx(U + \theta V)[V, V]| \leq C \norm{V}_{H^1_\dx}^2,
    \end{equation*}
    where $C$ contains $C_1$, $C_2$, $C_3$, $C_4$, $C_5$ and $|\ka|$, which completes the proof.
\end{proof}
From \eqref{eq:taylor_action}, we can also write $\Fh(U, V)$ as
\begin{equation*}
    \Fh(U, V) = \boS_\dx(U) - \boS_\dx(V) - D\boS_\dx(V)[U-V].
\end{equation*}
In particular, for $U=\psi^n$ and $V=e^n$, we denote $\Fh^n = \Fh(\Psi^n, \psi^n)$. From Proposition~\ref{prop:stability} and the fact that $\psi \in \boC^4(\mathcal{I}, H^5(\mathcal{J}))$, $\psi^n$ and $e^n$ satisfy $\norm{\psi^n}_{H^1_\dx}<C_1$ and $\norm{\psi^n + e^n}_{H^1_\dx} = \norm{\Psi^n}_{H^1_\dx}<C_2$. 
Therefore, we have the following estimate as a result of Proposition~\ref{prop:F_remainder_estimate}:
\begin{equation}\label{ineq:Hn_remainder_estimate}
    |\Fh^n| \leq C \norm{e^n}_{H^1_\dx}^2~\text{ for all }~ n\in\dbra{0}{\Nt}.
\end{equation}
The next lemma provides us with another expression of the operator $\Fh^n$ in terms of $e^n$ and $d^n$, which will be useful for the following statements.
\begin{lemma}\label{lemma:Fh_Rh}
    The operator $\Fh^n$ can be rewritten as $\Fh^n = \Rh^n + \frac{\mu}{2}\normpp{e^n}{2}^2$, where
    \begin{equation}\label{eq:Rh^n}
         \Rh^n = \frac{1}{2}\norm{\Dxp e^n}_{2,h,-}^2-\frac{\gs}{4}\norm{d^n}_{2,h,+}^2 +\frac{\ka}{4}\norm{\Dxp d^n}_{2,h,-}^2 -\frac{1}{2}\ipRp{\gs \dns^n+\ka\Dxx \dns^n}{|e^n|^2}.
    \end{equation}
\end{lemma}
The proof is given in Section~\ref{appx:useful_equalities}. Using this lemma, we can now close the coercivity of $\Fh^n$ against $\norm{e^n}_{H^1_\dx}$:
\begin{proposition}\label{prop:coercivity}
    For $n\in \N$ and assuming the ellipticiy condition in \eqref{eq:disc_ellipticity} for $\ka<0$, there exists a $\mu^*>0$ such that for all $\mu \geq \mu^*$,
        \begin{equation}\label{eq:e_en_equivalence}
            C\norm{e^n}_{H_h^1}^2 \leq\Fh^n 
        \end{equation}
    where $C>0$ is a constant independent of $\dx$, $\dt$ and $n$.
\end{proposition}

\begin{proof}
Before going directly to the proof, we have the following useful estimates from Lemma~\ref{lemma:H1_dn_en} for the second and third terms on the right-hand side of \eqref{eq:Rh^n}:
\begin{equation}\label{ineq:d_bound}
    \Big|-\frac{\gs}{4} \norm{d^n}_{2, \dx^+}\Big| \leq C_1 \norm{e^n}_{2, \dx^+}\quad\text{ and }\quad \Big|\frac{\ka}{4}\norm{\Dxp d^n}_{2,h,-}^2 \Big|\leq C_2 \norm{e^n}^2_{H^1_\dx},
\end{equation}
where $C_1$ and $C_2$ are constants independent of $\dx$ and $\dt$, and $C_2$ contains the term $|\ka|$. 
Since $q\in C(\mathcal{I},H^5(\mathcal{J}))$ and $q_x=0$ at the endpoints, the integral formulas for the interior and boundary differences give $\norm{\Dxx q^n}_\infty\leq\infnorm{\ptl_{xx}q(\cdot,t_n)}\leq C$. Using these uniform maximum-norm bounds, the fourth term on the right-hand side of \eqref{eq:Rh^n} satisfies
\begin{equation}\label{ineq:ip_rem}
    \Big|-\frac{1}{2} \ipRp{\gs \dns^n+\ka\Dxx \dns^n}{|e^n|^2}\Big| \leq \Bigg(\frac{|\gs|}{2}\norm{\psi^n}_\infty^2 + \frac{|\ka|}{2}\norm{\Dxx |\psi^n|^2}_\infty\Bigg)\norm{e^n}^2_{2, \dx^+} \leq C_3\norm{e^n}^2_{2, \dx^+},
\end{equation}
recalling that $q^n = \abs{\psi^n}^2$, and where $C_3$ depends on $|\ka|$ but does not depend on $\dx$ nor $\dt$. 
Passing to the proof of the proposition, we first consider the case $\ka \geq 0$. Using the first inequality in \eqref{ineq:d_bound} and the one in \eqref{ineq:ip_rem}, we get
\begin{equation*}
    \Fh^n \geq\frac{1}{2}\norm{\Dxp e^n}_{2,h,-}^2 +\Big(\frac{\mu}{2}-C_4\Big)\norm{e^n}_{2,h,+}^2, 
\end{equation*}
where $C_4$ contains the term $|\ka|$ and the constants $C_1$ and $C_3$. Notice that, we neglected the term $\frac{\ka}{4}\norm{\Dxp d^n}_{2,h,-}^2$ since $\ka\geq 0$. Taking $\mu^* = 2+2C_4$, we conclude the lower bound in \eqref{eq:e_en_equivalence} for the case $\ka \geq 0$. Now, we consider $\ka<0$. From \eqref{eq:dx_dn_1} and after some algebraic calculations, we have
\begin{equation}\label{eq:dx_dn_2}
    (\Dxp d^n)_j = 2\ip{\Psi^n_{j+\frac{1}{2}}}{(\Dxp e^n)_j} + 2\ip{e^n_{j+\frac{1}{2}}}{(\Dxp \psi^n)_j}, \quad \text{ for } j\in\dbra{0}{\J-1}.
\end{equation}
Let $\varepsilon>0$. Using the inequality $|a+b|^2\leq(1+\varepsilon)|a|^2+(1+\varepsilon^{-1})|b|^2$ in \eqref{eq:dx_dn_2}, we obtain
\begin{align*}
    |(\Dxp d^n)_j| \leq & 4(1 + \varepsilon)\big|\Psi_{j+\frac{1}{2}}^n\big|^2|(\Dxp e^n)_j|^2 + 4(1 + \varepsilon^{-1})|(\Dxp \psi^n)_j|^2 |e_{j+\frac{1}{2}}^n|^2
    \\ \leq & 4(1 + \varepsilon)\norm{\Psi^n}^2_\infty |(\Dxp e^n)_j|^2 + 4(1 + \varepsilon^{-1})\norm{\Dxp \psi^n}_\infty^2 |e_{j+\frac{1}{2}}^n|^2.
\end{align*}
Hence,
\begin{equation}\label{ineq:norm_Dx_dn}
    \norm{\Dxp d^n}_{2,h,-}^2 \leq 4(1+\varepsilon)\norm{\Psi^n}_\infty^2 \norm{\Dxp e^n}_{2,h,-}^2 +4(1+\varepsilon^{-1})C_5 \norm{e^n}_{2,h,+}^2.
\end{equation}
Here, $C_5$ can be chosen independently of the mesh, since $\norm{\Dxp\psi^n}_\infty\leq\infnorm{\ptl_x\psi(t_n)}\leq C$ by the fundamental theorem of calculus and the assumed regularity. Since $\ka<0$, multiplying \eqref{ineq:norm_Dx_dn} by $\ka/4$ reverses the inequality. Therefore, using the first inequality in \eqref{ineq:d_bound} and \eqref{ineq:ip_rem}, we have
\begin{equation}\label{ineq:Rh_k_neg}
\begin{aligned}
    \Rh^n \geq & \frac{1}{2}\norm{\Dxp e^n}^2_{2, \dx^-} + \frac{\ka}{4} \norm{\Dxp d^n}^2_{2, \dx^-}  - C_6\norm{e^n}^2_{2, \dx^+}
    \\\geq & \Big(\frac{1}{2} + \ka(1+\epsilon)\norm{\Psi^n}_\infty^2\Big) \norm{\Dxp e^n}^2_{2, \dx^-} + (\ka(1+\varepsilon^{-1})C_5 - C_6)\norm{e^n}^2_{2, \dx^+},
\end{aligned}
\end{equation}
where $C_6$ is a constant obtained by adding the ones from \eqref{ineq:d_bound} and \eqref{ineq:ip_rem}. Notice that, from Lemma~\ref{lemma:discrete_sobolev} and Proposition~\ref{prop:stability}, it holds that $\norm{\Psi^n}_\infty \leq C_7$, where $C_7$ does not depend on $\dx$ and $\dt$. Using the ellipticity condition \eqref{eq:disc_ellipticity}, we obtain
\begin{equation*}
    \frac{1}{2} + \ka(1+\varepsilon)\norm{\Psi^n}_\infty^2 = \frac{1}{2} \Big(1 + 2\ka  \norm{\Psi^n}_\infty^2 \Big) + \varepsilon \ka \norm{\Psi^n}_\infty^2 \geq \frac{\delta}{2} - \varepsilon |\ka| C_7^2.
\end{equation*}
Choosing $\varepsilon>0$ such that $\varepsilon |\ka|C_7^2 \leq \delta/4$, we get from inequality \eqref{ineq:Rh_k_neg}, 
\begin{equation*}
    \Fh^n \geq \frac{\delta}{4} \norm{\Dxp e^n}^2_{2, \dx^-} + \Big(\frac{\mu}{2} + \ka(1+\varepsilon^{-1})C_5 - C_6)\Big)\norm{e^n}^2_{2, \dx^+}.
\end{equation*}
We conclude the lower bound in \eqref{eq:e_en_equivalence} by taking $\mu_* = \delta/2 - 2\ka(1+\varepsilon^{-1})C_5+ 2C_6>0$.
\end{proof}

The proof of Proposition~\ref{prop:equivalence} follows directly from \eqref{ineq:Hn_remainder_estimate} and Proposition~\ref{prop:coercivity}.

Now, we introduce several properties of the operators $\Fh$ and $\Rh$ before going to the convergence result. The following lemmas state an exact increment identity and estimates of the operator $\Rh^n$ from the error equation in \eqref{eq:error-equation}.

\begin{lemma}[Increment identity]\label{lemma:exact_increment}
For $n\in \N$, we have the following increment identity:
\begin{equation}\label{eq:exact-increment}
 \Rh^{n+1}-\Rh^n = \gs\mathcal{C}_n+ \ka \boQ_n+ \ipRp{\boR^{n+\frac{1}{2}}}{e^{n+1}-e^n},
\end{equation}
where $\mathcal{C}_n$ and $\boQ_n$ contain the cubic and quasilinear terms, respectively, and are given by
\begin{align}
 \mathcal{C}^n = & -\frac{1}{4}\ipRp{q^{n+1} - q^n}{|e^{n+1}|^2+|e^n|^2} -\ipRp{d^{n+\frac{1}{2}}}{\ip{e^{n+\frac{1}{2}}}{\psi^{n+1}-\psi^n}},\label{eq:Cn} \\
 \boQ^n = & -\frac{1}{4}\ipRp{\Dxx(q^{n+1}-q^n)}{|e^{n+1}|^2+|e^n|^2}-\ipRp{\Dxx d^{n+\frac{1}{2}}}{\ip{e^{n+\frac{1}{2}}}{\psi^{n+1}-\psi^n}}.\label{eq:Qn}
\end{align}
\end{lemma}
The proof of this lemma can be found in Section~\ref{appx:useful_equalities}.
\begin{lemma}
    Both $\mathcal{C}^n$ and $\boQ^n$ satisfy
    \begin{equation}\label{eq:algebraic-remainders}
    |\mathcal{C}^n|+|\boQ^n| \leq C\dt\left(\Fh^{n+1}+\Fh^n\right),
    \end{equation}
    where $C$ is a constant independent of $\dx$, $\dt$ and $n$.
\end{lemma}
\begin{proof}
    We start by estimating the first term in both $\mathcal{C}^n$ and $\boQ^n$ given by \eqref{eq:Cn} and \eqref{eq:Qn}, respectively, and then the second terms. From the definition of $q^n$, we have the identity
\begin{equation*}
    \dns^{n+1}-\dns^n =\ip{\psi^{n+1}+\psi^n}{\psi^{n+1} - \psi^n}.
\end{equation*}
By the fundamental theorem of calculus, 
\begin{equation}\label{ineq:psi_bound}
\norm{\psi^{n+1}-\psi^n}_\infty\leq\dt\sup_{t\in\mathcal{I}}\infnorm{\ptl_t\psi(t)}.
\end{equation}
Moreover, since $q_x=0$ at the endpoints, the integral formulas for $\Dxx$ give
$$
\norm{\Dxx(q^{n+1}-q^n)}_\infty\leq\dt\sup_{t\in\mathcal{I}}\infnorm{\ptl_{xxt}q(\cdot,t)}
$$
at every space point. Therefore, we have
\begin{equation*}
    \norm{\dns^{n+1}-\dns^n}_{\infty} \leq \big(\norm{\psi^{n+1}}_{\infty}+\norm{\psi^n}_{\infty}\big) \norm{\psi^{n+1}-\psi^n}_{\infty} \leq C_1\dt ~ \text{ and }~\norm{\Dxx(\dns^{n+1}-\dns^n)}_{\infty}\leq C_2\dt,
\end{equation*}
where $C_1$ and $C_2$ are constants independent of $\dt$. Consequently, using the definition of $\ipRp{\cdot}{\cdot}$ and the coercivity estimate in Proposition~\ref{prop:coercivity}, we obtain the bounds
\begin{equation}\label{ineq:first_estimates_Cn_Qn}
\begin{aligned}
    \Big|\ipRp{\dns^{n+1}-\dns^n}{|e^{n+1}|^2+|e^n|^2}\Big| &\leq C_1\dt\Big(\normpp{e^{n+1}}{2}^2+\normpp{e^n}{2}^2\Big) \leq C_3\dt\left(\Fh^{n+1}+\Fh^n\right),
    \\ \Big|\ipRp{\Dxx(\dns^{n+1}-\dns^n)}{|e^{n+1}|^2+|e^n|^2}\Big| &\leq C_2\dt\Big(\normpp{e^{n+1}}{2}^2+\normpp{e^n}{2}^2\Big)\leq  C_4\dt\left(\Fh^{n+1}+\Fh^n\right),
\end{aligned}
\end{equation}
where $C_3$ and $C_4$ are constants containing $C_1$ and $C_2$, and the ones from \eqref{eq:e_en_equivalence}. These estimates control the first term in both $\mathcal{C}^n$ and $\boQ^n$. We now proceed to the second term in $\mathcal{C}^n$ in \eqref{eq:Cn}. Using Cauchy--Schwarz for the product $\ipRp{\cdot}{\cdot}$ gives
\begin{equation}\label{ineq:second_Cn_estimate}
    \left|\ipRp{d^{n+\frac{1}{2}}}{\ip{e^{n+\frac{1}{2}}}{\psi^{n+1} - \psi^n}}\right|\leq\normpp{d^{n+\frac{1}{2}}}{2}\normpp{\ip{e^{n+\frac{1}{2}}}{\psi^{n+1} - \psi^n}}{2}.
\end{equation}
The term $\normpp{d^{n+\frac{1}{2}}}{2}$ is bounded using the first estimate in \eqref{ineq:H1_dn_en} at $n$ and $n+1$, obtaining
\begin{equation}\label{ineq:Cn_first_term_estimate}
    \normpp{d^{n+\frac{1}{2}}}{2}
    \leq \frac{1}{2}\left(\normpp{d^{n+1}}{2}+\normpp{d^n}{2}\right)
    \leq C\left(\normpp{e^{n+1}}{2}+\normpp{e^n}{2}\right).
\end{equation}
On the other hand, for the term $\normpp{\ip{e^{n+\frac{1}{2}}}{\psi^{n+1}-\psi^n}}{2}$, we use the inequality $|\ip{z}{w}|\leq |z||w|$ and the bound \eqref{ineq:psi_bound} to get
\begin{equation}\label{ineq:Cn_second_term_estimate}
    \normpp{\ip{e^{n+\frac{1}{2}}}{\psi^{n+1}-\psi^n}}{2} \leq \norm{\psi^{n+1}-\psi^n}_{\infty} \normpp{e^{n+\frac{1}{2}}}{2} \leq C\dt\left(\normpp{e^{n+1}}{2}+\normpp{e^n}{2}\right).
\end{equation}
Combining \eqref{ineq:Cn_first_term_estimate} and \eqref{ineq:Cn_second_term_estimate} with \eqref{ineq:second_Cn_estimate}, and using Young's inequality gives
\begin{equation}\label{ineq:second_Cn_estimate_final}
    \left|\ipRp{d^{n+\frac{1}{2}}}
    {\ip{e^{n+\frac{1}{2}}}{\delta\psi^n}}\right|
    \leq C\dt\left(\Fh^{n+1}+\Fh^n\right).
\end{equation}
Therefore, adding the first inequality in \eqref{ineq:first_estimates_Cn_Qn} with \eqref{ineq:second_Cn_estimate_final}, we conclude the bound for $\boC_n$
\begin{equation}\label{ineq:C_n_estimate}
    |\mathcal{C}^n|\leq C\dt\left(\Fh^{n+1}+\Fh^n\right).
\end{equation}
It remains to estimate the second term in $\boQ^n$. By the summation-by-parts formula
of Lemma~\ref{lemma:sum_by_parts} and Cauchy--Schwarz, we have
\begin{equation*}
    \Big|\ipRp{\Dxx d^{n+\frac{1}{2}}}{\ip{e^{n+\frac{1}{2}}}{\psi^{n+1}-\psi^n}}\Big|\leq\normpm{\Dxp d^{n+\frac{1}{2}}}{2}\normpm{\Dxp\ip{e^{n+\frac{1}{2}}}{\psi^{n+1}-\psi^n}}{2}.
\end{equation*}
The term $\normpm{\Dxp d^{n+\frac{1}{2}}}{2}$ is estimated in a similar manner as in \eqref{ineq:Cn_first_term_estimate}, but using the second estimate in \eqref{ineq:H1_dn_en}. For the remaining factor, the integral formula for the forward difference gives 
$$
\norm{\Dxp(\psi^{n+1}-\psi^n)}_\infty\leq\dt\sup_{t\in\mathcal{I}}\infnorm{\ptl_{xt}\psi(t)}\leq C\dt.
$$
Applying Corollary~\ref{corollary:product_rule} and the definition of the discrete norms in \eqref{eq:norms}, we have
\begin{align*}
    \normpm{\Dxp\ip{e^{n+\frac{1}{2}}}{\psi^{n+1}-\psi^n}}{2} \leq& \norm{\Dxp(\psi^{n+1}-\psi^n)}_{\infty}\normpp{e^{n+\frac{1}{2}}}{2}+\norm{\psi^{n+1}-\psi^n}_{\infty}\normpm{\Dxp e^{n+\frac{1}{2}}}{2}\\
    \leq& C\dt\norm{e^{n+\frac{1}{2}}}_{H^1_\dx}.
\end{align*}
As for $\boC^n$, it follows that
\begin{equation}\label{ineq:Q_n_second_term_estimate}
    \Big|\ipRp{\Dxx d^{n+\frac{1}{2}}}{\ip{e^{n+\frac{1}{2}}}{\delta\psi^n}}\Big| \leq C\dt\left(\Fh^{n+1}+\Fh^n\right),
\end{equation}
and adding \eqref{ineq:Q_n_second_term_estimate} with the second inequality in \eqref{ineq:first_estimates_Cn_Qn}, we obtain the bound for $\boQ^n$
\begin{equation}\label{ineq:Q_n_estimate}
    |\boQ^n|\leq C\dt\left(\Fh^{n+1}+\Fh^n\right).
\end{equation}
We conclude the proof of \eqref{eq:algebraic-remainders} by adding both estimates \eqref{ineq:C_n_estimate} and \eqref{ineq:Q_n_estimate}.
\end{proof}

\begin{lemma}\label{lem:residual-pairing}
For every $m\in\dbra{1}{N}$ and every $\varepsilon>0$,
\begin{equation}\label{eq:residual-pairing}
    \Big| \sum_{n=0}^{m-1}\ipRp{\boR^{n+\frac12}}{e^{n+1}-e^n}\Big| \leq \varepsilon\Big( \Fh^m + \dt\sum_{n=1}^{m-1}\Fh^n\Big)+ C(\dt^2+\dx^2)^2.
\end{equation}
where $C$ is independent of $\dx$, $\dt$ and $m$, and contains $\varepsilon$ and $T$.
\end{lemma}

\begin{proof}
Let us fix $m\in\dbra{1}{N}$. Performing a change of indices in the summations, the left-hand side of \eqref{eq:residual-pairing} satisfies
\begin{equation}\label{eq:residual_1}
    \sum_{n=0}^{m-1}\ipRp{\boR^{n+\frac{1}{2}}}{e^{n+1}-e^n} = \ipRp{\boR^{m-\frac{1}{2}}}{e^m} - \ipRp{\boR^{\frac{1}{2}}}{e^0} - \sum_{n=1}^{m-1}\ipRp{\boR^{n+\frac{1}{2}} - \boR^{n-\frac{1}{2}}}{e^n}.
\end{equation}
Since we start with the exact initial condition, $e^0=0$, the second term on the right-hand side vanishes. We now estimate the first and third terms separately. By Cauchy--Schwarz, the inequality \eqref{eq:residual-L2-bounds} in Proposition~\ref{prop:second_order} and Proposition~\ref{prop:coercivity}, we have
\begin{equation*}
    \Big|\ipRp{\boR^{m-\frac{1}{2}}}{e^m}\Big| \leq C_1(\dt^2 + \dx^2)\normpp{e^m}{2} \leq  C_2(\dt^2 + \dx^2)(\Fh^m)^{1/2},
\end{equation*}
where $C_1$ and $C_2$ are constants independent of $\dx$, $\dt$ and $m$.
Hence, Young's inequality, $ab \leq \varepsilon a^2 + b^2/(4\varepsilon)$, implies that, for every $\varepsilon>0$,
\begin{equation}\label{eq:residual_1_estimate}
    \Big|\ipRp{\boR^{m-\frac{1}{2}}}{e^m}\Big| \leq \varepsilon \Fh^m + C_3(\dt^2 + \dx^2)^2,
\end{equation}
where $C_3$ contains the constant $C_2$ and $\varepsilon$. For the third term on the right-hand side of \eqref{eq:residual_1}, we use again Cauchy--Schwarz, the inequality \eqref{eq:residual-time-increment} in Proposition~\ref{prop:second_order}, Proposition~\ref{prop:coercivity} and Young's inequality to get
\begin{equation*}
    \Big|\ipRp{\boR^{n+\frac{1}{2}} - \boR^{n-\frac{1}{2}}}{e^n}\Big| \leq C_4\dt(\dt^2 + \dx^2)(\Fh^n)^{1/2} \leq \varepsilon\dt\Fh^n + C_5\dt(\dt^2 + \dx^2)^2, \text{ for } n=\dbra{1}{m-1},
\end{equation*}
where $C_4$ and $C_5$ are constants independent of $\dx$, $\dt$ and $m$, and $C_5$ contains $\varepsilon$. Summing from $n=1$ to $n=m-1$ gives
\begin{equation}\label{eq:residual_2_estimate}
    \sum_{n=1}^{m-1}\Big|\ipRp{\boR^{n+\frac{1}{2}} - \boR^{n-\frac{1}{2}}}{e^n}\Big| \leq \varepsilon\dt\sum_{n=1}^{m-1}\Fh^n + C_5 (m-1)\dt(\dt^2 + \dx^2)^2 \leq \varepsilon\dt\sum_{n=1}^{m-1}\Fh^n + C_5T(\dt^2 + \dx^2)^2,
\end{equation}
where we used that $(m-1)\dt \leq T$ in the last step. Finally, we conclude the proof of the inequality \eqref{eq:residual-pairing} by combining both estimates in \eqref{eq:residual_1_estimate} and \eqref{eq:residual_2_estimate} with \eqref{eq:residual_1}.
\end{proof}

\begin{lemma}[$L^2$-increment]\label{lem:L2-increment}
For every $n\in\dbra{0}{N-1}$, the operator $\Fh^n$ satisfies the following increment inequality:
\begin{equation}\label{eq:L2-increment}
 \norm{e^{n+1}}_{2,h,+}^2-\norm{e^n}_{2,h,+}^2\leq C\dt\Big(\Fh^{n+1}+\Fh^n+\normpp{\boR^{n+\frac{1}{2}}}{2}^2\Big),
\end{equation}
where $C$ is a constant independent of $\dx$ and $\dt$.
\end{lemma}

\begin{proof}
    Applying the product $\ipCp{\cdot}{e^{n+\frac{1}{2}}}$ to \eqref{eq:error-equation} and taking the imaginary part, we obtain
    \begin{equation}\label{eq:L2_error_equation}
    \begin{aligned}
         \Im \frac{i}{\dt}\ipCp{e^{n+1}-e^n}{e^{n+\frac{1}{2}}} = & \Im\ipCp{\Dxx e^{n+\frac{1}{2}}}{e^{n+\frac{1}{2}}} \\ & + \gs \Im\ipCp{  \Dns^{n+\frac{1}{2}}\odot e^{n+\frac{1}{2}} + d^{n+\frac{1}{2}}\odot\psi^{n+\frac{1}{2}}}{e^{n+\frac{1}{2}}}
         \\ & + \ka \Im\ipCp{(\Dxx \Dns^{n+\frac{1}{2}})\odot e^{n+\frac{1}{2}}+(\Dxx d^{n+\frac{1}{2}})\odot\psi^{n+\frac{1}{2}}}{e^{n+\frac{1}{2}}}
         \\ & + \Im\ipCp{\boR^{n+\frac{1}{2}}}{e^{n+\frac{1}{2}}},
    \end{aligned}
    \end{equation}
    where we used that $\Dns^{n+\frac{1}{2}}\odot \Psi^{n+\frac{1}{2}} - \dns^{n+\frac{1}{2}}\odot\psi^{n+\frac{1}{2}} = \Dns^{n+\frac{1}{2}}\odot e^{n+\frac{1}{2}} + d^{n+\frac{1}{2}}\odot\psi^{n+\frac{1}{2}}$, and similarly for the last term on the right-hand side of \eqref{eq:error-equation}. We examine each term in \eqref{eq:L2_error_equation} separately. The left-hand side can be rewritten as
    \begin{equation*}
        \Im \frac{i}{\dt}\ipCp{e^{n+1}-e^n}{e^{n+\frac{1}{2}}} = \frac{1}{2\dt}\left(\norm{e^{n+1}}_{2,h,+}^2-\norm{e^n}_{2,h,+}^2\right).
    \end{equation*}
    Using summation by parts in Lemma~\ref{lemma:sum_by_parts}, the first term on the right-hand side satisfies
    \begin{equation*}
        \Im \ipCp{\Dxx e^{n+\frac{1}{2}}}{e^{n+\frac{1}{2}}} = -\Im\ipCm{\Dxp e^{n+\frac{1}{2}}}{\Dxp e^{n+\frac{1}{2}}} = \Im \normpm{\Dxp e^{n+\frac{1}{2}}}{2}^2 = 0.
    \end{equation*}
    For the second and third terms, notice that, by Lemma~\ref{lemma:distributive}, we have
    \begin{equation*}
        \Im \ipCp{\Dns^{n+\frac{1}{2}}\odot e^{n+\frac{1}{2}}}{e^{n+\frac{1}{2}}} = \Im \ipCp{\Dns^{n+\frac{1}{2}}}{|e^{n+\frac{1}{2}}|^2} = 0,
    \end{equation*}
    and
    \begin{equation*}
        \Im \ipCp{(\Dxx \Dns^{n+\frac{1}{2}})\odot e^{n+\frac{1}{2}}}{e^{n+\frac{1}{2}}} = \Im \ipCp{\Dxx \Dns^{n+\frac{1}{2}}}{|e^{n+\frac{1}{2}}|^2} = 0,
    \end{equation*}
    since $\Dns^{n+\frac{1}{2}}$ and $\Dxx \Dns^{n+\frac{1}{2}}$ are real-valued. Therefore, equation \eqref{eq:L2_error_equation} reduces to
    \begin{equation}\label{ineq:L2_error_equation_reduced}
    \begin{aligned}
        \frac{1}{2\dt}\left(\norm{e^{n+1}}_{2,h,+}^2-\norm{e^n}_{2,h,+}^2\right) = & \gs \Im\ipCp{d^{n+\frac{1}{2}}\odot\psi^{n+\frac{1}{2}}}{e^{n+\frac{1}{2}}} \\ + & \ka \Im\ipCp{(\Dxx d^{n+\frac{1}{2}})\odot\psi^{n+\frac{1}{2}}}{e^{n+\frac{1}{2}}} + \Im\ipCp{\boR^{n+\frac{1}{2}}}{e^{n+\frac{1}{2}}}.
    \end{aligned}
    \end{equation}
We now bound the terms on the right-hand side of \eqref{ineq:L2_error_equation_reduced}. First, the term $\Im\ipCp{d^{n+\frac{1}{2}}\odot\psi^{n+\frac{1}{2}}}{e^{n+\frac{1}{2}}}$ can be bounded using Cauchy--Schwarz, Lemma~\ref{lemma:discrete_holder} and Lemma~\ref{lemma:H1_dn_en} 
to get
\begin{equation}\label{ineq:linear_estimate_en}
\begin{aligned}
    \Big|\ipCp{d^{n+\frac{1}{2}}\odot\psi^{n+\frac{1}{2}}}{e^{n+\frac{1}{2}}}\Big| \leq & \norm{\psi^{n+\frac{1}{2}}}_{\infty} \normpp{d^{n+\frac{1}{2}}}{2} \normpp{e^{n+\frac{1}{2}}}{2} \leq C \Big(\normpp{e^{n+1}}{2} + \normpp{e^{n}}{2}\Big)^2 \\
    \leq & C \Big(\norm{e^{n+1}}_{H_h^1}^2 + \norm{e^{n}}_{H_h^1}^2\Big) \leq C (\Fh^{n+1} + \Fh^n),
\end{aligned}
\end{equation}
where in the last two inequalities we used that $(a + b)^2 \leq 2(a^2 + b^2)$ and the coercivity of $\Fh^n$ in Proposition~\ref{prop:coercivity}. The constant $C$ contains all the constants from the used lemmas and propositions, and is independent of $\dx$ and $\dt$. Similarly, for the quasilinear term $\Im\ipCp{(\Dxx d^{n+\frac{1}{2}})\odot\psi^{n+\frac{1}{2}}}{e^{n+\frac{1}{2}}}$, we first use the distributive property of $\ipCp{\cdot}{\cdot}$ and then the summation by parts in Lemma~\ref{lemma:sum_by_parts} to obtain
\begin{equation*} 
    \ipCp{(\Dxx d^{n+\frac{1}{2}})\odot\psi^{n+\frac{1}{2}}}{e^{n+\frac{1}{2}}} = -\ipCm{\Dxp d^{n+\frac{1}{2}}}{\Dxp (\overline{\psi^{n+\frac{1}{2}}} \odot e^{n+\frac{1}{2}})}.
\end{equation*}
Hence, by Cauchy--Schwarz, Lemma~\ref{lemma:discrete_H1_algebra} and Lemma~\ref{lemma:H1_dn_en}
, we have
\begin{equation}\label{ineq:quasilinear_estimate_en}
\begin{aligned}
    \Big|\ipCp{(\Dxx d^{n+\frac{1}{2}})\odot\psi^{n+\frac{1}{2}}}{e^{n+\frac{1}{2}}}\Big| \leq & C \normpm{\Dxp d^{n+\frac{1}{2}}}{2} \norm{\psi^{n+\frac{1}{2}}}_{H^1_\dx} \norm{e^{n+\frac{1}{2}}}_{H_h^1} \\
    \leq & C \Big(\norm{e^{n+1}}_{H_h^1}^2 + \norm{e^{n}}_{H_h^1}^2\Big) \leq C (\Fh^{n+1} +    \Fh^n),
\end{aligned}
\end{equation}
where in the last two inequalities we used again that $(a + b)^2 \leq 2(a^2 + b^2)$ and the coercivity of $\Fh^n$ in Proposition~\ref{prop:coercivity}. As before, the constant $C$ contains all the constants from the used lemmas and propositions, and is independent of $\dx$ and $\dt$. Finally, for the residual term $\ipCp{\boR^{n+\frac{1}{2}}}{e^{n+\frac{1}{2}}}$, we use Cauchy--Schwarz, Young's inequality and Proposition~\ref{prop:coercivity} to get
\begin{equation}\label{ineq:residual_estimate_en}
\Big| \ipCp{\boR^{n+\frac{1}{2}}}{e^{n+\frac{1}{2}}} \Big| \leq \frac{1}{2} \normpp{\boR^{n+\frac{1}{2}}}{2}^2 + \frac{1}{2} \normpp{e^{n+\frac{1}{2}}}{2}^2 \leq C \Big(\normpp{\boR^{n+\frac{1}{2}}}{2}^2 + \Fh^{n+1} + \Fh^n \Big).    
\end{equation}
Combining the estimates \eqref{ineq:linear_estimate_en}--\eqref{ineq:residual_estimate_en}\silentref{ineq:quasilinear_estimate_en} with \eqref{ineq:L2_error_equation_reduced}, and multiplying by $2\dt$, gives the desired result \eqref{eq:L2-increment}, concluding the proof.
\end{proof}
We are finally in position to prove our main result concerning the second-order convergence in the discrete $H^1$-norm established in Theorem~\ref{thm:error_convergence}:
\begin{proof}[Proof of Theorem~\ref{thm:error_convergence}]
Let us choose $\mu>0$ as in Proposition~\ref{prop:coercivity}. For this value, the relative energy $\Fh^n$ controls $\norm{e^n}_{H^1_\dx}^2$, and therefore, it is enough to obtain a uniform estimate for $\Fh^n$. Let us fix $m\in\dbra{1}{N}$. Summing the increment identity of Lemma~\ref{lemma:exact_increment} from $n=0$ to $n=m-1$, we obtain
\begin{equation*}
    \Rh^m-\Rh^0 =\gs\sum_{n=0}^{m-1}\mathcal C^n +\ka\sum_{n=0}^{m-1}\mathcal \boQ^n +\sum_{n=0}^{m-1}\ipRp{\boR^{n+\frac{1}{2}}}{e^{n+1}-e^n},
\end{equation*}
and the estimate \eqref{eq:algebraic-remainders} in the same lemma yields
\begin{equation}\label{ineq:sum_Rh_increment}
    \Rh^m-\Rh^0 \leq C\dt\sum_{n=0}^{m-1}\Big(\Fh^{n+1}+\Fh^n\Big) +\Bigg|\sum_{n=0}^{m-1}\ipRp{\boR^{n+\frac{1}{2}}}{e^{n+1}-e^n} \Bigg|.
\end{equation}
On the other hand, summing the $L^2$-error increment of Lemma~\ref{lem:L2-increment} from $n=0$ to $n=m-1$ and multiplying it by $\mu/2$ gives
\begin{equation}\label{ineq:L2_mu_increment}
    \frac{\mu}{2}\left(\normpp{e^m}{2}^2-\normpp{e^0}{2}^2\right) \leq C\dt\sum_{n=0}^{m-1}\Big(\Fh^{n+1}+\Fh^n\Big) + C\dt\sum_{n=0}^{m-1}\normpp{\boR^{n+\frac{1}{2}}}{2}^2.
\end{equation}
Adding both inequalities \eqref{ineq:sum_Rh_increment} and \eqref{ineq:L2_mu_increment}, and using Lemma~\ref{lemma:Fh_Rh}, we obtain
\begin{equation}\label{ineq:Fh_increment}
    \Fh^m-\Fh^0 \leq C\dt\sum_{n=0}^{m-1}\Big(\Fh^{n+1}+\Fh^n\Big) + \Bigg|\sum_{n=0}^{m-1}\ipRp{\boR^{n+\frac{1}{2}}}{e^{n+1}-e^n} \Bigg| + C\dt\sum_{n=0}^{m-1}\normpp{\boR^{n+\frac{1}{2}}}{2}^2.
\end{equation}
Since we start from the exact initialization, we have $e^0=0$ and hence $\Fh^0=0$. We now apply our results above to estimate the terms on the right-hand side of \eqref{ineq:Fh_increment}. First, by Lemma~\ref{lem:residual-pairing}, for every $\varepsilon>0$, we have
\begin{equation}\label{ineq:residual_pairing_estimate_conv}
    \Big| \sum_{n=0}^{m-1}\ipRp{\boR^{n+\frac{1}{2}}}{e^{n+1}-e^n}\Big| \leq \varepsilon\Big( \Fh^m + \dt\sum_{n=1}^{m-1}\Fh^n\Big) + C(\dt^2+\dx^2)^2.
\end{equation}
On the other hand, summing the estimate \eqref{eq:residual-L2-bounds} in Proposition~\ref{prop:second_order} from $n=0$ to $n=m-1$ gives
\begin{equation}\label{ineq:residual_L2_estimate_squared}
    \dt\sum_{n=0}^{m-1}\normpp{\boR^{n+\frac{1}{2}}}{2}^2 \leq Cm\dt (\dt^2+\dx^2)^2 \leq C(\dt^2+\dx^2)^2,
\end{equation}
where we used that $m\dt\leq T$ and $C$ contains $T$ in the last inequality. Putting together \eqref{ineq:residual_pairing_estimate_conv} and \eqref{ineq:residual_L2_estimate_squared} in \eqref{ineq:Fh_increment}, and using elementary algebra of sums, we obtain
\begin{equation*}
    \Fh^m \leq (C\dt+\varepsilon)\Fh^m + (C+\varepsilon)\dt\sum_{n=1}^{m-1}\Fh^n + C(\dt^2+\dx^2)^2.
\end{equation*}
Choosing $\varepsilon=1/4$ and $\dt_0=\min\{1, 1/(4C)\}$, we have, for $\dt\leq \dt_0$, $C\dt + \varepsilon < 1/2$. Therefore, the first term on the right-hand side can be absorbed into the left-hand side, obtaining
\begin{equation*}
    \Fh^m \leq C\dt\sum_{n=1}^{m-1}\Fh^n + C(\dt^2+\dx^2)^2.
\end{equation*}
This estimate holds for every $m\in\dbra{1}{N}$. The discrete Gronwall inequality in Lemma~\ref{lem:discrete_gronwall} gives
\begin{equation*}
    \max_{0\leq m\leq N}\Fh^m \leq C(\dt^2+\dx^2)^2.
\end{equation*}
Finally, the coercivity property of $\Fh^n$ in Proposition~\ref{prop:coercivity} implies
\begin{equation*}
    \norm{e^m}_{H^1_\dx} \leq C(\Fh^m)^{\frac{1}{2}} \leq C(\dt^2+\dx^2),
\end{equation*}
uniformly for $0\leq m\leq N$, obtaining \eqref{eq:global-stability} and finishing the proof.
\end{proof}


\section{Numerical experiments}\label{sec:numerics}
One of the main features of \eqref{QGP} in the real line and with $\gs=1$ and $\ka>0$ is the existence of localized traveling-wave solutions called solitons or \textit{bright solitons}, which are localized profiles decaying exponentially to zero at infinity. These solutions generalize the family of localized solitons known in the semilinear case $\ka=0$; see \cite{fibich,bethuel2008existence}. Precisely, they correspond to smooth \emph{finite-energy} traveling wave solutions of the form 
\begin{equation}
    \label{def:bright}
    \psi(x,t)=u(x-ct-x_0)e^{i\vartheta(t)}
\end{equation}
for some speed $c\in\R$, and initial position $x_0\in \R$, where $u$ is a complex-valued function with a localized derivative,  referred to as the \emph{profile} of the soliton, and $\vartheta$ is a real-valued function describing the temporal frequency.

The explicit formula for the profile and temporal frequency of bright solitons to \eqref{QGP}, given in  \cite{krolikowski_solitons_2000}, is:
\begin{equation}
\label{bright}
    u_{\omega,\ka}(x) =F_{\omega,\ka}^{-1}(|x|)e^{-icx/2} \quad \text{ and }\quad  \vartheta_{c, \omega}(t)={-\Big(\frac{c^2}{4}+\omega   \Big)}t,
\end{equation}
for all $x,t,c\in\R$, $\omega>0$, where $F_{\omega, \ka}: (0, \sqrt{2\omega}] \to [0, \infty)$ is the strictly decreasing function given by
\begin{equation}\label{eq:implibright}
    F_{\omega,\ka}(y)=\frac{1}{\sqrt{\omega}}\atanh\Big(\frac{1}{\sqrt{2\omega}}\sqrt{\frac{2\omega-y^2}{1+2\ka y^2}}\Big)+2\sqrt{\ka}\atan\Big( \sqrt{2\ka} \sqrt{\frac{2\omega-y^2}{1+2\ka y^2}}\Big).
\end{equation}
Here, $u_{\omega, \ka}$ is the solution to equation \eqref{eq:soliton_eq} obtained from the Fr\'echet derivative of the action $\boS$. Using the results in \cite{Colin2} or in \cite{ilievquasilinstab}, we deduce that the bright solitons in \eqref{bright}--\eqref{eq:implibright} are orbitally stable for all $\omega>0$ and $c\in \R$.

We use bright solitons as benchmark solutions to validate our numerical scheme \eqref{scheme}. These solutions are suitable for this purpose, since they have an explicit form, up to the inversion of a scalar function, and their exponential decay to zero allows us to consider a finite computational domain with homogeneous Neumann boundary conditions without introducing significant errors. 

\begin{figure}[ht]
	\centering
	\begin{tabular}{cc}
		\resizebox{0.45\textwidth}{!}{
			\includegraphics{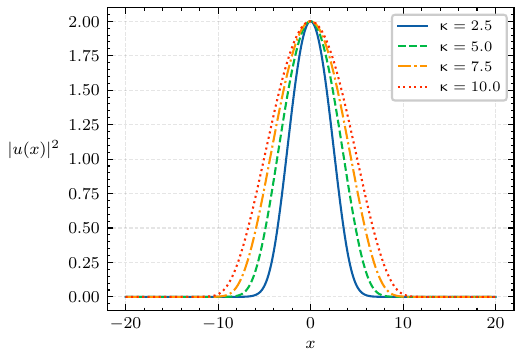}}
	&
	\resizebox{0.43\textwidth}{!}{\includegraphics{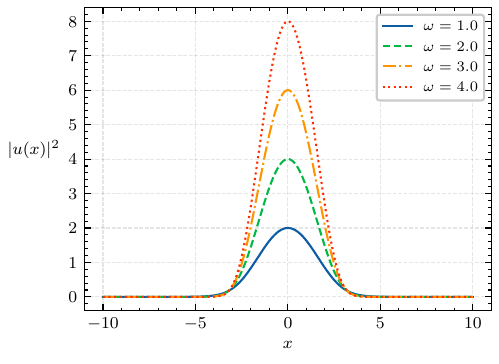}}
\end{tabular}
\caption{Plots of several squared moduli of bright soliton profiles obtained using \eqref{bright}. The left panel shows profiles with $\omega=1$ and $\ka\in \{2.5, 5, 7.5, 10\}$. The right panel shows profiles with $\ka=1$ and $\omega\in\{1, 2, 3, 4\}$.}
\label{fig:Bright_plots}
\end{figure}

\begin{figure}[ht]
	\centering
	\begin{tabular}{cc}
		\resizebox{0.45\textwidth}{!}{
			\includegraphics{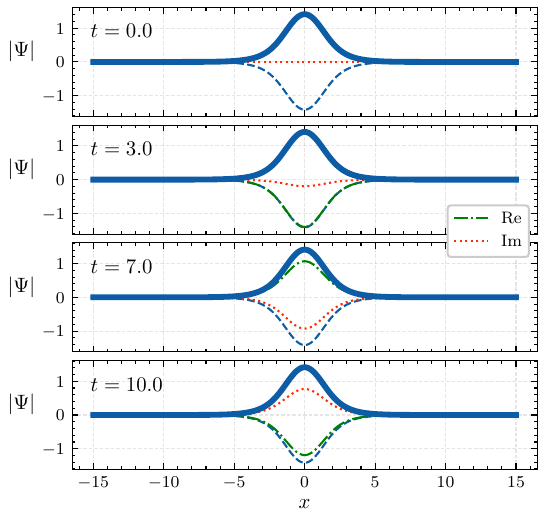}}
	&
	\resizebox{0.45\textwidth}{!}{\includegraphics{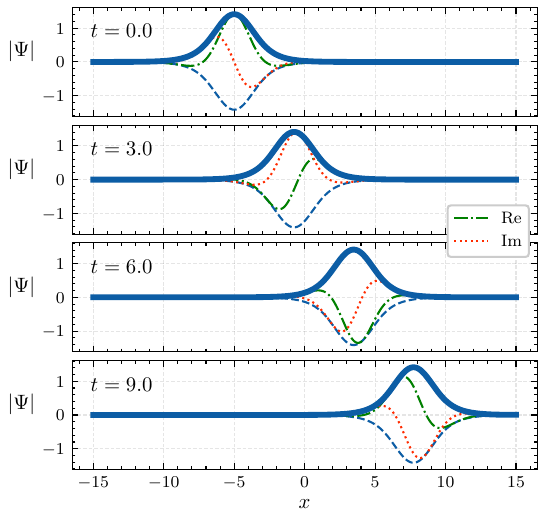}}
\end{tabular}
\caption{Snapshots of the absolute value of bright soliton solutions with $\ka = 1/4$ and $\omega = 1$, along with their real and imaginary parts, obtained using \eqref{def:bright} and \eqref{bright}. The left panel shows a standing wave with $c=0$, and the right panel shows a traveling wave with $c=\sqrt{2}$.}
\label{fig:snap_bright}
\end{figure}

We approximate $F^{-1}_{\omega,\ka}$ given in \eqref{eq:implibright} at each spatial grid point by solving a least-squares problem in the squared modulus, with tolerance $10^{-12}$. The resulting values are interpolated using a cubic spline. For traveling solutions, interpolation is performed on an auxiliary grid with twice the number of spatial points and twice the domain length as in the original simulation, to avoid evaluating outside the interpolation range. Figures~\ref{fig:Bright_plots} and \ref{fig:snap_bright} show various examples of solitons obtained by this method.

\paragraph{Complex Newton iteration.} Since the scheme \eqref{scheme} is implicit, each time step requires solving a nonlinear system of the form $F(z,\conj{z})=0$, with $F:\,\C^{2(\J+1)}\to\C^{2(\J+1)}$, where $z\in\mathbb{C}^{\J+1}$ contains the variables at the new time step. We employ a Newton method adapted to complex-valued maps via Wirtinger derivatives, as in \cite{dutto_extending_2019}. We define the blocks $A, B$ and the Jacobian $J_z$ as 
\begin{equation}\label{eq:blockJac}
A:=\frac{\partial F}{\partial z}(z_k,\conj{z}_k), \qquad
B:=\frac{\partial F}{\partial \conj{z}}(z_k,\conj{z}_k), \qquad J_z(z_k, \conj{z}_k):=\begin{pmatrix}
A & B\\[2pt]
\conj{B} & \conj{A}
\end{pmatrix}.
\end{equation}
The Newton direction $(\Delta z_k,\Delta\conj{z_k})$ is computed by solving
\begin{equation*}
J_z(z_k, \conj{z}_k)\begin{pmatrix}\Delta z_k\\[7pt] \Delta\conj{z_k}\end{pmatrix}
= - \begin{pmatrix} F(z_k,\conj{z}_k)\\[7pt] \conj{F(z_k,\conj{z}_k)}\end{pmatrix}.
\end{equation*}
We notice that $z_k$ and $\conj{z}_k$ are treated as independent variables. This can cause the Jacobian to lose its symmetry. Therefore, we enforce the update by setting
\begin{equation*}
z_{k+1}=z_k+\Delta z_k,\qquad \conj{z}_{k+1}=\conj{z_{k+1}}.
\end{equation*}
This also applies to blocks $\conj{A},\conj{B}$ in \eqref{eq:blockJac}. The iteration stops at the first $k$ such that
\begin{equation*}
\normp{\Delta z_k}{2}\leq \delta_{\mathrm{abs}},
\end{equation*}
or if the algorithm reaches the maximum number of iterations $\Nt_{\mathrm{max}}$. For all simulations, we set $\delta_{\mathrm{abs}} = 10^{-10}$ and $\Nt_{\mathrm{max}} = 10^3$.

\paragraph {Boundary control.} To ensure that no significant perturbation reaches the boundary, we monitor the boundary values at each time step $n$. Specifically, we check that
\begin{equation}\label{eq:cond_bound}
|\Psi_j^n|<\delta, \qquad j\in\{0,\J\}, \text{ for all } n\in \dbra{0}{\Nt}.
\end{equation}
If this condition is not satisfied, we stop the simulation and return the last computed step.

Notice that, on truncated domains, the boundary conditions may generate artificial reflections when outgoing radiation reaches the boundary, due to the preservation of mass and energy. Transparent or absorbing boundary conditions have been developed for the cubic NLS equation with zero background \cite{antoine_besse_descombes_2006,zheng_2006,antoine_arnold_besse_ehrhardt_schaedle_2008}. However, to the best of our knowledge, no analogous construction is currently available for the quasilinear Gross--Pitaevskii equation.

\begin{figure}[ht]
	\centering
	\begin{tabular}{cc}
		\resizebox{0.458\textwidth}{!}{
			\includegraphics{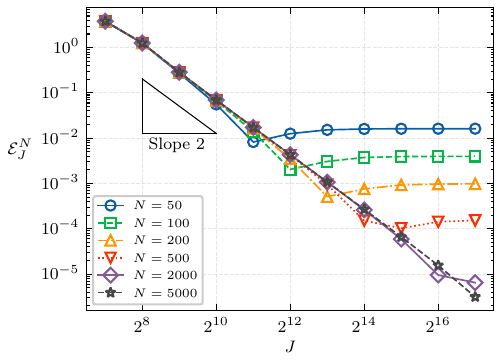}}
	&
	\resizebox{0.45\textwidth}{!}{\includegraphics{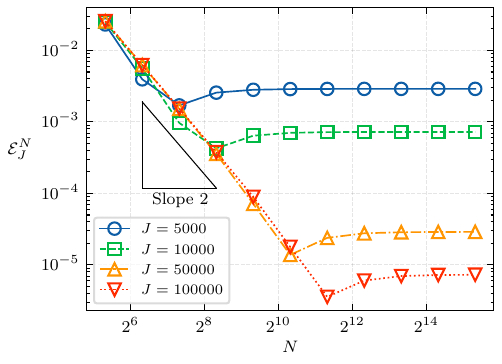}}
\end{tabular}
\caption{Errors of the numerical scheme \eqref{scheme}. The left panel shows the space error, and the right panel the time error.}
\label{fig:Bright_errors}
\end{figure}

We set as initial condition a bright soliton with parameters $\ka = 1/4$, $c=\sqrt{2}/2$, $\omega = 1$ and $x_0=0$, and denote the travelling wave solution as $\psi_{\omega, \ka}(x, t)$ given in \eqref{def:bright}. We run various instances of the simulation, refining the spatial and temporal meshes, until the final time $T=5$ on a grid with length $L=100$. In fact, we consider $\J=2^p$ for $p\in \dbra{7}{17}$, and $\Nt=T \cdot 2^q$ for $q\in \dbra{5}{15}$. We compute the discrete $H^1$-norm of the difference between the simulation and the theoretical solution at time $T$:
\begin{equation}\label{error}
    \boE^{\Nt}_\J = \norm{\Psi^{\Nt} - \psi_{\omega, \ka}(\cdot_j,T)}_{H^1_\dx},
\end{equation}
where $\cdot_j$ denotes the evaluation on the discrete space grid. The values of $c$ and $L$ were chosen to ensure that the soliton remains away from the boundary, satisfying the condition \eqref{eq:cond_bound}, and to ensure that the boundary conditions do not introduce significant errors. The results are displayed in Figure~\ref{fig:Bright_errors}, where we observe that the error \eqref{error} decreases following a two-slope line in both cases. Performing a linear regression for the biggest values of $\J$ and $\Nt$ from $p=7$ to $p=17$, and from $q=5$ to $q=11$, respectively, we obtain the slopes $m\approx-2.0417$ in space and $m\approx -2.1559$ in time. These results numerically validate the second-order convergence of the scheme in both variables, which is in accordance
with the convergence result proved in Theorem~\ref{thm:error_convergence}.

To show the order of convergence for bright solitons with more precision, we define the following experimental order of convergence (EOC) for $\Nt$ and $\J$ by 
\begin{equation*}
    \mathrm{EOC}_{\J} = \log_2\left(\mathcal{E}_\J^{\Nt}/\mathcal{E}_{2\J}^{\Nt}\right) \quad \text{ and } \quad \mathrm{EOC}_{\Nt} = \log_2\left(\mathcal{E}_\J^{\Nt}/\mathcal{E}_{\J}^{2\Nt}\right).
\end{equation*}
The numerical values are shown in  Tables \ref{tab:EOC_bright_N} and \ref{tab:EOC_bright_NT}, which corroborate that the errors are of order $2$ in time and space.
\begin{table}[ht]
    \centering
    \begin{tabular}{cccccccc}
\toprule
 $J$ & $N=50$ & $N=100$ & $N=200$ & $N=500$ & $N=2000$ & $N=5000$ \\
\midrule
$2^{7}$ & 1.806 & 1.791 & 1.788 & 1.787 & 1.787 & 1.787 \\
$2^{8}$ & 2.215 & 2.164 & 2.152 & 2.149 & 2.148 & 2.148 \\
$2^{9}$ & 2.321 & 2.099 & 2.050 & 2.037 & 2.034 & 2.034 \\
$2^{10}$ & 2.754 & 2.291 & 2.073 & 2.019 & 2.010 & 2.009 \\
$2^{11}$ & -0.858 & 2.773 & 2.281 & 2.043 & 2.005 & 2.003 \\
$2^{12}$ & -0.255 & -0.845 & 2.784 & 2.171 & 2.011 & 2.002 \\
$2^{13}$ & -0.059 & -0.261 & -0.837 & 2.775 & 2.041 & 2.007 \\
$2^{14}$ & -0.015 & -0.060 & -0.264 & 0.184 & 2.170 & 2.026 \\
$2^{15}$ & -0.004 & -0.015 & -0.061 & -0.437 & 2.771 & 2.106 \\
$2^{16}$ & -0.001 & -0.004 & -0.015 & -0.097 & 0.193 & 2.468 \\
\bottomrule
\end{tabular}

    \caption{$\mathrm{EOC}_{\J}$ for increasing $\Nt$ and $\J$.}
    \label{tab:EOC_bright_N}
\end{table}
\begin{table}[ht]
    \centering
    \begin{tabular}{lrrrr}
\toprule
$N$ & $J=5000$ & $J=10000$ & $J=50000$ & $J=100000$ \\
\midrule
$T \cdot 2^{4}$ & 2.512 & 2.139 & 2.037 & 2.034 \\
$T \cdot 2^{5}$ & 1.578 & 2.501 & 2.034 & 2.021 \\
$T \cdot 2^{6}$ & -0.714 & 1.547 & 2.077 & 2.025 \\
$T \cdot 2^{7}$ & -0.155 & -0.709 & 2.300 & 2.072 \\
$T \cdot 2^{8}$ & -0.037 & -0.154 & 2.705 & 2.298 \\
$T \cdot 2^{9}$ & -0.009 & -0.037 & -0.876 & 2.694 \\
$T \cdot 2^{10}$ & -0.002 & -0.009 & -0.248 & -0.870 \\
$T \cdot 2^{11}$ & -0.001 & -0.002 & -0.058 & -0.248 \\
$T \cdot 2^{12}$ & -0.000 & -0.001 & -0.014 & -0.058 \\
\bottomrule
\end{tabular}

    \caption{$\mathrm{EOC}_{\Nt}$ for increasing $\J$ and $\Nt$}
    \label{tab:EOC_bright_NT}
\end{table}

\paragraph{Computational setup.} All experiments across the work were implemented in Python, using NumPy and SciPy for array operations and sparse linear algebra. The largest simulations were performed on a 64-bit virtual machine with 32 GB of RAM. The simplest visualizations, such as snapshots or soliton plots, were performed on a local 64-bit M4-chip MacBook Pro with 24 GB of unified memory.

\paragraph{Data availability.}
All the code used for the simulations can be found in \url{https://github.com/ElTapia/Quasilinear-GP.git}.

\appendix
\section{Well-posedness of the Cauchy problem}\label{appx:well-posedness}
In the case of periodic boundary conditions, Feola and Iandoli in \cite{feola2022local} provide a local well-posedness result for quasilinear  Hamiltonian Sch\"odringer equations in $H^p(\T_L)$, where $\T_L=\R/(2L\Z)$ is the torus identified with $[-L,L]$,
for $p>13$ Schr\"odinger equations, under some ellipticity condition for the density of the energy. In the next result, we show how to apply their framework to \eqref{QGP}, 
to obtain a well-posedness result, including a blow-up alternative for maximal solutions.
For this purpose, we use the following ellipticity margin function that measures the loss of ellipticity:
	\begin{equation*}
		\boD_L(\psi)
:=\inf_{x\in[-L,L]}
\bigl(1+2{\upkappa}|\psi(x)|^2\bigr),
\quad \text{for }  \psi\in C([-L,L];\C).
	\end{equation*}

\begin{theorem}
\label{thm:periodic-maximal}
	Let $\upkappa,\mathfrak{s}\in\R$, $L>0$, $p>13$, and
	$\psi_0\in H^p(\T_L)$ satisfy
	$\boD_L(\psi_0)>0$.
	Then there is a unique maximal solution to 
     \begin{equation}\label{eq:qgp}
   i\partial_t\psi=\partial_{xx}\psi
   +\mathfrak{s}|\psi|^2\psi+\upkappa\psi\partial_{xx}|\psi|^2,
   \qquad x\in\T_L,
 \end{equation}
	on $[0,T^*)$, with $T^*\in(0,\infty]$, such
		that, for every $0<T<T^*$,
	\begin{equation*}
		\psi\in C([0,T];H^p(\T_L))
		\cap C^1([0,T];H^{p-2}(\T_L)),
		\quad \text{and}\quad
		{\boD_L(\psi(t))}>0
		\ \text{for all } {t\in[0,T^*)}.
	\end{equation*}
    Also, the solution depends continuously on its initial datum in
	$H^p$ on every compact subinterval of
	its existence interval. Moreover, if $T^*<\infty$, the following blow-up alternative holds:
		\begin{equation}\label{eq:periodic-blowup-alternative}
		\lim_{t\uparrow T^*}
		\Big(\norm{\psi(t)}_{H^p(\T_L)}
		+{(\boD_L(\psi(t)))^{-1}}\Big)=+\infty.
	\end{equation}
	In particular, for $\upkappa\geq0$, a finite maximal time implies
	$
		\lim_{t\uparrow T^*}\norm{\psi(t)}_{H^p(\T_L)}=+\infty.
	$
	\end{theorem}

\begin{proof}
We first establish a uniform local existence property: for every
$M>0$ and $\eta>0$, there exists $\tau=\tau(M,\eta)>0$
such that any initial condition $\psi_0\in H^p(\T_L)$ satisfying
\begin{equation}
    \label{hyp:IC}
    \|\psi_0\|_{H^p}\le M\quad  \text{ and } \quad \boD_L(\psi_0)\geq\eta
\end{equation} 
admits a unique solution $\psi$ on $[0,\tau]$ with
$\boD_L(\psi(t ))>0$, for  $t\in [0,\tau].$

In the case $\upkappa>0$, this is a simple consequence of Theorem~1.2 of Feola--Iandoli\cite{feola2022local}, 
as we verify below.
For  $\upkappa<0$,  the argument is more involved, and we need to 
 modify the Hamiltonian density associated with \eqref{eq:qgp}, 
by introducing a well-chosen cut-off function, in the spirit of \cite{Pluschke1996, ChenFeldman2010}, 
to prevent a possible loss of ellipticity.

More precisely,  we  consider the Hamiltonian density
\[
 F_\chi(z,w)=\frac{\upkappa}{2}\chi(|z|^2)
 (\bar zw+z\bar w)^2-\frac{\mathfrak{s}}{2}|z|^4, \quad \text { for }z,w\in\C,
\]
 where $\chi\in C^\infty([0,\infty);[0,1])$ is a function to be chosen.
 The corresponding
Feola--Iandoli ellipticity quantities are
\begin{equation}
\label{hyp1.1}
 1+(F_\chi)_{w\bar w}=1+\upkappa\chi(|z|^2)|z|^2,
 \qquad
 (1+(F_\chi)_{w\bar w})^2-|(F_\chi)_{\bar w\bar w}|^2
 =1+2\upkappa\chi(|z|^2)|z|^2.   
\end{equation}
For $\upkappa\geq 0$, we simply take $\chi\equiv1$.
 Hence, $F(z,w)=O((|z|+|w|)^4$ at the origin, 
 	so  $F_\chi$ has a zero of order at least
	three, as required in \cite{feola2022local}.
 Also, it satisfies Hypothesis~1.1 in \cite{feola2022local},
 since both functions in \eqref{hyp1.1} are bounded below by $1$.
 Thus, Theorem~1.2 in \cite{feola2022local} provides the local existence of a solution.
 In the case $a:=-\upkappa>0$, let $0<\eta\le1$ such that $\boD_L(\psi_0)\geq \eta$
 so that 
 \begin{equation}\label{cota0}
 \|\psi_0\|_\infty^2\le(1-\eta)/(2a).
 \end{equation}
 We set the 
 cut-off function $\chi\in C^\infty([0,\infty);[0,1])$ satisfying
\[
 \chi(r)=1\quad\text{for }r\le\frac{1-\eta/2}{2a},
 \qquad
 \chi(r)=0\quad\text{for }r\geq\frac{1-\eta/4}{2a}.
\]
Notice that this cut-off depends only on
	$\eta$ and $a$, and is therefore the same for all initial
	data satisfying the above assumptions.
In this manner, 
$1+\upkappa\chi(|z|^2)|z|^2\geq(1+\eta/4)/2$
and $1+2\upkappa\chi(|z|^2)|z|^2\geq\eta/4$
for all $z\in \C$,
 and both functions in \eqref{hyp1.1} are bounded  below by $\eta/8$.
 Thus, we can apply again Theorem~1.2 of \cite{feola2022local}.
	Since the Hamiltonian density and its ellipticity constants
	are fixed, the estimates in the proof of that theorem
	allow us to choose a common time $\tau_0>0$ and a constant
	$R>0$, depending only on $M,\eta$, such that for every
	initial condition satisfying \eqref{hyp:IC}, 
  there exists a unique
solution
  $\psi\in C([0,\tau_0];H^p(\T_L))
		\cap C^1([0,\tau_0];H^{p-2}(\T_L))$. 
  Moreover, we have the bound 
\begin{equation}
    \label{cotaR}
     \sup_{0\le t\le\tau_0}\|\psi(t)\|_{H^p}\le R,
\end{equation}
where $\tau_0$ and $R$ depend only on $M$ and $\eta$.
Writing the equation associated with the modified Hamiltonian
as 
\begin{equation}
    \label{modfied}
i    \partial_t\psi=i\mathcal F_\chi(\psi):=
\partial_{xx}\psi+\mathfrak{s}|\psi|^2\psi
+\upkappa\chi(|\psi|^2)\psi\partial_{xx}|\psi|^2\\
+\frac{\upkappa}{2}\chi'(|\psi|^2)
\psi\bigl(\partial_x|\psi|^2\bigr)^2,
\end{equation}
we observe that $\mathcal F_\chi$ maps bounded subsets of
$H^p(\T_L)$ into bounded subsets of
$H^{p-2}(\T_L)$. Hence, since $p-2>1/2$,
the fundamental theorem of calculus and Sobolev embedding give
$$
 \|\psi(t)-\psi_0\|_\infty
 \leq C_p\int_0^t
 \|\mathcal F_\chi(\psi(r))\|_{H^{p-2}}\,dr
 \leq C(M,\eta)t,
 \qquad 0\leq t\leq\tau_0.
$$
Using \eqref{cota0} and \eqref{cotaR},  we deduce that
$$
 \begin{aligned}
 \bigl|\|\psi(t)\|_\infty^2-\|\psi_0\|_\infty^2\bigr|
 \leq
 \bigl(\|\psi(t)\|_\infty+\|\psi_0\|_\infty\bigr)
 \|\psi(t)-\psi_0\|_\infty
 \leq C(M,\eta)t.
 \end{aligned}
$$
Hence, choosing $\tau\in(0,\tau_0]$ 
small enough, depending only on $M$, $\eta$, we obtain
$$ \|\psi(t)\|_\infty^2
 <\frac{1-\eta/2}{2a},
 \qquad 0\leq t\leq\tau.
$$
This implies that 
$\chi(|\psi(t,x)|^2)=1$ and
$\chi'(|\psi(t,x)|^2)=0$ for all
$(t,x)\in[0,\tau]\times\T_L$.
Thus, the modified equation \eqref{modfied} coincides with
\eqref{eq:qgp} on $[0,\tau]$.
Moreover,
$$
 \boD_L(\psi(t))
 =1-2a\|\psi(t)\|_\infty^2
 \geq {\eta}/{2}>0,
 \qquad 0\leq t\leq\tau.
$$
This concludes the proof of existence of $\tau=\tau(M,\eta)$.

The uniqueness and continuous dependence follow from \cite{feola2022local} by using the same cutoff. 
Indeed, two solutions of \eqref{eq:qgp} with the same initial condition and positive ellipticity margin satisfy the same modified equation \eqref{modfied} on a sufficiently small common time interval. 
The uniqueness for equation \eqref{modfied} then gives local uniqueness for \eqref{eq:qgp}. This argument can be repeated by restarting the equation at any time at which the ellipticity margin remains positive.

We now construct the maximal solution.
By uniqueness, the local solutions can be
concatenated, so we can obtain a unique maximal solution
$\psi$ on an interval $[0,T^*)$,
where $T^*\in(0,\infty]$.
Continuous dependence on compact subintervals
of $[0,T^*)$ follows by successive applications of the local well-posedness result.

Now, suppose that $T^*<\infty$ and that
\eqref{eq:periodic-blowup-alternative} does not hold.
Then there exist $M, \eta>0$ and a sequence
$t_n\uparrow T^*$ such that
$
 \|\psi(t_n)\|_{H^p(\T_L)}\leq M,
 $ and 
 $\boD_L(\psi(t_n))\geq\eta.
$
By the uniform local existence property established
above, we can restart \eqref{eq:qgp} at each time
$t_n$, obtaining a solution on
$[t_n,t_n+\tau(M,\eta)]$.
Since $\tau(M,\eta)>0$ is independent of $n$,
we have
$
 t_n+\tau(M,\eta)>T^*
$
for all sufficiently large $n$.
This contradicts the maximality of $T^*$ and
proves \eqref{eq:periodic-blowup-alternative}.
\end{proof}

We finally deduce a maximal well-posedness result for the Neumann problem
\eqref{QGP} on the interval $\mathcal J=[-L,L]$, by using an even reflection.
Indeed, consider the torus $\T_{2L}=\R/(4L\Z)$. For $\psi :  \mathcal J \to \C$, let $\E \psi$ denote its
even, $4L$-periodic extension, defined on one period by
\begin{equation*}
 (\E \psi )(x)=
 \begin{cases}
   \psi (x),&-L\le x\le L,\\
   \psi(2L-x),&L\le x\le3L.
 \end{cases}
 \end{equation*}
Notice that, for $p\in\mathbb N$ and
$\psi\in H^p(\mathcal J)$,
the even periodic extension $\E\psi$ belongs to
$H^p(\T_{2L})$ if and only if the following
compatibility conditions hold:
\begin{equation}\label{comp_cond_init}
    \partial_x^{2\ell+1}\psi(-L)
    =\partial_x^{2\ell+1}\psi(L)=0,
    \quad \text{ for all }\ell\geq0
    \text{ such that }
     2\ell+1\leq p-1.
\end{equation}

The condition in \eqref{comp_cond_init} is exactly the classical necessary condition 
for the existence of a solution $ \psi$ in $C([0,T];H^p(\mathcal J))\cap C^1([0,T];H^{p-2}(\mathcal J))$
 of \eqref{QGP}, satisfying
$\boD_L(\psi(t))>0$ for all $t\in[0,T]$.
Indeed, we can argue by induction on the odd integer
$m=2\ell+1\leq p-1$.
The case $m=1$ follows from the Neumann boundary conditions.
Suppose that all lower odd-order derivatives vanish at
$x=\pm L$. Differentiating \eqref{QGP} $m-2$ times with
respect to $x$ and evaluating at the endpoints, we obtain
\begin{equation}
    \label{borde}
 \partial_x^m\psi
 +\upkappa\psi\,\partial_x^m|\psi|^2=0.
\end{equation}
Here we have used the induction hypothesis and the fact
that $\partial_t\partial_x^{m-2}\psi(t,\pm L)=0$.
Moreover, the induction hypothesis gives
\begin{equation}
    \label{borde2}
 \partial_x^m|\psi|^2
 =\bar\psi\,\partial_x^m\psi
 +\psi\,\overline{\partial_x^m\psi}
 \qquad\text{at }x=\pm L,
\end{equation}
so that,  multiplying \eqref{borde2} by $\upkappa\psi$, and adding \eqref{borde},
\begin{equation}
\label{borde3}
 (1+\upkappa|\psi|^2)\partial_x^m\psi
 +\upkappa\psi^2\overline{\partial_x^m\psi}=0.   
\end{equation}
Multiplying  \eqref{borde3} by $1+\upkappa \abs{\psi}^2$, 
and subtracting the complex conjugate of \eqref{borde3} multiplied by $\upkappa {\psi}^2$,
we eliminate $\overline{\partial_x^m\psi}$,  and obtain
\[
 (1+2\upkappa|\psi|^2)\partial_x^m\psi=0
 \qquad\text{at }x=\pm L.
\]
Since $\boD_L(\psi(t))>0$, it follows that
$\partial_x^m\psi(t,\pm L)=0$, completing the induction.

We are now in a position to state the well-posedness of \eqref{QGP}.
\begin{corollary}
\label{lucario}
Let $\boJ=[-L,L]$, $p\in\N$ with $p>13$, and
$\mathfrak{s},\upkappa\in\R$.
Assume that $\psi_0\in H^p(\mathcal J)$
satisfies the compatibility conditions
\eqref{comp_cond_init} and
$
 \boD_{L}(\psi_0)>0,
$
Then the Neumann problem \eqref{QGP}
admits a unique maximal solution $\psi$
on $[0,T^*)$, with $T^*\in(0,\infty]$,
such that, for every $0<T<T^*$,
\[
 \psi\in C([0,T];H^p(\mathcal J))
 \cap C^1([0,T];H^{p-2}(\mathcal J)),
\]
and
$ \boD_{L}(\psi(t))>0,$
 for all $t\in[0,T^*)$
The solution depends continuously on its
initial datum in $H^p(\mathcal J)$
on every compact subinterval of its
existence interval.
Moreover, if $T^*<\infty$, then
\begin{equation}\label{eq:neumann-blowup-alternative}
 \lim_{t\uparrow T^*}
 \left(
 \|\psi( t)\|_{H^p(\mathcal J)}
 +\bigl(\boD_{L}(\psi( t))\bigr)^{-1}
 \right)=+\infty.
\end{equation}
\end{corollary}

\begin{proof}
In view of the compatibility conditions \eqref{comp_cond_init},
the even extension $\E\psi_0$ belongs to $H^p(\T_{2L})$, as remarked above. In addition,
$
 \boD_{2L}(\E\psi_0)=\boD_L(\psi_0)>0.
$
Thus, Theorem~\ref{thm:periodic-maximal} gives the existence of a unique
maximal periodic solution $\widetilde\psi$ on $[0,T^*)$, with $T^*\in(0,\infty]$.

Since \eqref{QGP} and the $\psi_0$ are invariant
under the reflection $x\mapsto 2L-x$, the uniqueness implies that 
$$
 \widetilde\psi(t,2L-x)=\widetilde\psi(t,x),
 \qquad 0\leq t<T^*.
$$
Thus, by periodicity, we deduce that  $\partial_x\widetilde\psi(t,\pm L)=0$.
Therefore, the function $\psi:=\widetilde\psi|_{\mathcal J}$
is a solution of the Neumann problem \eqref{QGP},
with $\boD_L(\psi(t))>0$.

To prove the uniqueness of the equation, let $\phi$ be another solution
of \eqref{QGP} satisfying the properties in the statement. 
As remarked above, $\phi$ satisfies the compatibility conditions
\eqref{comp_cond_init}. Hence, its even extension $\E\phi$ is a solution
of the periodic problem in $H^p(\T_{2L})$.
By invoking the uniqueness of the periodic problem, we infer that  $\E\phi=\widetilde\psi$,
and therefore $\phi=\psi$.
The same argument shows that $\psi$ is maximal.
Continuous dependence follows from the periodic result
and the boundedness of the even extension operator
on the compatible Sobolev subspace.

Finally, since
\[
 \boD_{2L}(\widetilde\psi( t))=\boD_L(\psi( t))
\]
and the norms $\|\widetilde\psi( t)\|_{H^p(\T_{2L})}$
and $\|\psi( t)\|_{H^p(\mathcal J)}$
are equivalent, the blow-up alternative
\eqref{eq:neumann-blowup-alternative}
follows directly from
\eqref{eq:periodic-blowup-alternative}.
\end{proof}

\section{Useful equalities proofs}\label{appx:useful_equalities}
This appendix is devoted to proving the equalities established in Lemma~\ref{lemma:Fh_Rh} and Lemma~\ref{lemma:exact_increment}.

\begin{proof}[Proof of Lemma~\ref{lemma:Fh_Rh}]
    Let us write $\Fh^n$ as follows:
    \begin{equation}\label{eq:H_h_terms}
        \Fh^n = \frac{1}{2} \mathcal{H}_1^n + \frac{\gs}{4} \mathcal{H}_{2}^n + \frac{\ka}{4} \mathcal{H}_{3}^n + \frac{\mu}{2}\mathcal{H}_{4}^n 
    \end{equation}
    where
    \begin{equation*}
        \mathcal{H}_{1}^n = \normpm{\Dxp \Psi^n}{2}^2 - \normpm{\Dxp \psi^n}{2}^2 + 2\ipRp{\Dxx \psi^n}{e^n}, 
    \end{equation*}
    \begin{equation*}
        \mathcal{H}_{2}^n = -\normpp{\Psi^n}{4}^4 + \normpp{\psi^n}{4}^4 + 4\ipRp{\dns^n \odot \psi^n}{e^n}, 
    \end{equation*}
    \begin{equation*}
        \mathcal{H}_{3}^n =\normpm{\Dxp \Dns^n}{2}^2 -\normpm{\Dxp \dns^n}{2}^2 + 4\ipRp{\Dxx \dns^n \odot \psi^n}{e^n}, 
    \end{equation*}
    \begin{equation*}
        \mathcal{H}_{4}^n = \normpp{\Psi^n}{2}^2 - \normpp{\psi^n}{2}^2 + 2\ipRp{\psi^n}{e^n}, 
    \end{equation*}
    Using Lemma~\ref{lemma:sum_by_parts}, the first term $\mathcal{H}_{1}^n$ becomes
    \begin{equation}\label{eq:H_kin}
    \begin{aligned}
        \mathcal{H}_{1}^n = & \normpm{\Dxp \Psi^n}{2}^2 - 2\ipRm{\Dxp \psi^n}{\Dxp e^n} - \normpm{\Dxp \psi^n}{2}^2
        \\ = & \normpm{\Dxp \Psi^n}{2}^2 - 2\ipRm{\Dxp \psi^n}{\Dxp \Psi^n} + \normpm{\Dxp \psi^n}{2}^2 = \normpm{\Dxp e^n}{2}^2
    \end{aligned}
    \end{equation}
    For the other terms, we use the fact that, from the definition of $d^n$, we have $d^n = |\psi^n + e^n|^2 - |\psi^n|^2 = |e^n|^2 + 2\ip{\psi^n}{e^n}$. Therefore, from Lemma~\ref{lemma:distributive}, the second term $\mathcal{H}_{\text{cub}}^n$ becomes
    \begin{equation}\label{eq:H_cub}
    \begin{aligned}
        \mathcal{H}_{2}^n =& -\normpp{\Dns^n}{2}^2 + \normpp{\dns^n}{2}^2 + 4\ipRp{\dns^n}{\ip{\psi^n}{e^n}_\odot} \\
        = & -\ipRp{\Dns^n+\dns^n}{d^n} + 2\ipRp{\dns^n}{d^n} - 2\ipRp{\dns^n}{|e^n|^2} \\
        = & -\normpp{d^n}{2}^2 - 2\ipRp{\dns^n}{d^n} + 2\ipRp{\dns^n}{d^n} - 2\ipRp{\dns^n}{|e^n|^2} \\
        = & -\normpp{d^n}{2}^2 - 2\ipRp{\dns^n}{|e^n|^2}.
    \end{aligned}
    \end{equation}
    Similarly, the third term $\mathcal{H}_{\text{quasi}}^n$ becomes
    \begin{equation}\label{eq:H_quasi}
        \mathcal{H}_{3}^n = \normpm{\Dxp d^n}{2}^2 + 2 \ipRm{\Dxp \dns^n}{\Dxp |e^n|^2}.
    \end{equation}
    Finally, following the same procedure as in $\mathcal{H}_{\text{kin}}^n$, the last term $\mathcal{H}_{\text{mass}}^n$ satisfies
    \begin{equation}\label{eq:H_mass}
        \mathcal{H}_{4}^n = \normpp{\Psi^n}{2}^2 + 2\ipRp{\psi^n}{e^n} - \normpp{\psi^n}{2}^2 = \normpp{\Psi^n}{2}^2 - 2\ipRp{\Psi^n}{\psi^n} + \normpp{\psi^n}{2}^2 =  \normpp{e^n}{2}^2.
    \end{equation}
    Replacing \eqref{eq:H_kin}--\eqref{eq:H_mass}\silentref{eq:H_cub}\silentref{eq:H_quasi} in \eqref{eq:H_h_terms} completes the proof
\end{proof}

\begin{proof}[Proof of Lemma~\ref{lemma:exact_increment}]
Applying $\ipRp{\cdot}{e^{n+1}-e^n}$ in equation \eqref{eq:error-equation} gives
\begin{equation}\label{eq:error_product}
    \frac{1}{\dt} \ipRp{i(e^{n+1}-e^n)}{e^{n+1}-e^n} =  \ipRp{\Dxx e^{n+\frac{1}{2}}}{e^{n+1}-e^n} + \gs \boG_{1}^n + \ka \boG_{2}^n + \ipRp{R^{n+\frac{1}{2}}}{e^{n+1}-e^n},
\end{equation}
where
\begin{equation*}
    \boG_{1}^n = \ipRp{\Dns^{n+\frac{1}{2}}\odot \Psi^{n+\frac{1}{2}}-\dns^{n+\frac{1}{2}}\odot\psi^{n+\frac{1}{2}}}{e^{n+1}-e^n}
\end{equation*}
and
\begin{equation*}
    \boG_{2}^n = \ipRp{(\Dxx \Dns^{n+\frac{1}{2}})\odot\Psi^{n+\frac{1}{2}}-(\Dxx \dns^{n+\frac{1}{2}})\odot\psi^{n+\frac{1}{2}}}{e^{n+1}-e^n}.
\end{equation*}
First, we notice that the left-hand side of \eqref{eq:error_product} vanishes since $ \ipRp{i(e^{n+1}-e^n)}{e^{n+1}-e^n} =  \ipRp{i}{|e^{n+1}-e^n|^2}$ and $|e^{n+1}-e^n|^2\in \R$. We now proceed to the right-hand side of \eqref{eq:error_product}. From Corollary~\ref{corollary:product_rule}, we can rewrite the first term as
\begin{align*}
    \ipRp{\Dxx e^{n+\frac{1}{2}}}{e^{n+1}-e^n} = & -\frac{1}{2}\ipRm{\Dxp (e^{n+1}+e^n)}{\Dxp (e^{n+1}-e^n)} \\
    = &  -\frac{1}{2}\Big( \normpm{\Dxp e^{n+1}}{2}^2 -\normpm{\Dxp e^{n}}{2}^2 \Big).
\end{align*}
For the second term $\boG_{2}^n$, the definition of $d^n$ yields
\begin{equation*}
    \Dns^{n+\frac{1}{2}} \odot \Psi^{n+\frac{1}{2}} - \dns^{n+\frac{1}{2}} \odot \psi^{n+\frac{1}{2}} = \Dns^{n+\frac{1}{2}} \odot e^{n+\frac{1}{2}} + d^{n+\frac{1}{2}} \odot \psi^{n+\frac{1}{2}}.
\end{equation*}
Therefore, Lemma~\ref{lemma:distributive} yields
\begin{equation}\label{eq:H_cub_exp}
\begin{aligned}
    \boG_{1}^n = & \ipRp{\Dns^{n+\frac{1}{2}} \odot e^{n+\frac{1}{2}}}{e^{n+1}-e^n} + \ipRp{ d^{n+\frac{1}{2}} \odot \psi^{n+\frac{1}{2}} }{e^{n+1}-e^n} \\
    & = \frac{1}{2}\ipRp{\Dns^{n+\frac{1}{2}}}{|e^{n+1}|^2 - |e^n|^2} + \ipRp{d^{n+\frac{1}{2}}}{\ip{\psi^{n+\frac{1}{2}}}{e^{n+1}- e^n}_\odot}.
\end{aligned}
\end{equation}
We recall that $d^n = |e^n|^2 + 2\ip{\psi^n}{e^n}$ from the proof of Lemma~\ref{lemma:Fh_Rh}. Subtracting $d^{n+1}$ to $d^n$, and applying Lemma~\ref{lemma:discrete_product} with $U=\psi^{n+1}$, $V=e^{n+1}$, $W=\psi^n$ and $Z=e^n$, we obtain
\begin{equation}\label{eq:psi_n_med_d_n}
    2\ip{\psi^{n+\frac{1}{2}}}{e^{n+1} - e^n} = d^{n+1} - d^n - (|e^{n+1}|^2 - |e^n|^2) - 2\ip{e^{n+\frac{1}{2}}}{\psi^{n+1} - \psi^n}.
\end{equation}
Substituting \eqref{eq:psi_n_med_d_n} into \eqref{eq:H_cub_exp} and cancelling similar terms, we get
\begin{equation}\label{eq:H_cub_2}
    \boG_{1}^n = \frac{1}{4}\Big(\normpp{d^{n+1}}{2}^2 - \normpp{d^n}{2}^2\Big) + \frac{1}{2}\ipRp{q^{n+\frac{1}{2}}}{|e^{n+1}|^2 - |e^
    n|^2} - \ipRp{d^{n+\frac{1}{2}}}{\ip{e^{n+\frac{1}{2}}}{\psi^{n+1}-\psi^n}}.
\end{equation}
From Lemma~\ref{lemma:discrete_product} with $U=q^{n+1}$, $V=|e^{n+1}|^2$, $W=q^n$ and $Z=|e^n|^2$, we can rewrite the second term on the right-hand side of \eqref{eq:H_cub_2} as
\begin{equation*}
    \frac{1}{2}\ipRp{q^{n+\frac{1}{2}}}{|e^{n+1}|^2 - |e^n|^2} = \frac{1}{2}\ipRp{q^{n+1}}{|e^{n+1}|^2} - \frac{1}{2}\ipRp{q^n}{|e^n|^2} - \frac{1}{4} \ipRp{q^{n+1} - q^n}{|e^{n+1}|^2 + |e^n|^2}.
\end{equation*}
Therefore,
\begin{equation}\label{eq:H_cub_final}
    \boG_{1}^n = \frac{1}{4}\Big(\normpp{d^{n+1}}{2}^2 - \normpp{d^n}{2}^2\Big) + \frac{1}{2}\Big(\ipRp{q^{n+1}}{|e^{n+1}|^2} -\ipRp{q^n}{|e^n|^2} \Big) + \boC^n.
\end{equation}
For the term $\boG_{2}^n$, we proceed analogously as the term above, applying an extra discrete summation by parts step in Lemma~\ref{lemma:sum_by_parts} for $\Dxx d^{n+\frac{1}{2}}$ to obtain
\begin{equation}\label{eq:H_quasi_final}
    \boG_{2}^n = -\frac{1}{4}\Big(\normpm{\Dxp d^{n+1}}{2}^2 - \normpm{\Dxp d^{n}}{2}^2\Big) + \frac{1}{2}\Big(\ipRp{\Dxx q^{n+1}}{|e^{n+1}|^2} -\ipRp{\Dxx q^n}{|e^n|^2} \Big) + \boQ^n.
\end{equation}
Replacing \eqref{eq:H_cub_final} and \eqref{eq:H_quasi_final} in \eqref{eq:error_product}, gives \eqref{eq:exact-increment}.
\end{proof}

\begin{merci}
The author acknowledges the support of the CDP C2EMPI, as well as the French State under the France-2030 program, the University of Lille, the Initiative of Excellence of the University of Lille, the European Metropolis of Lille for their funding and support of the R-CDP-24-004-C2EMPI project.
The author gratefully acknowledges the many helpful suggestions of A.~de Laire and G.~Dujardin during the preparation of the paper.
\end{merci}

\bibliographystyle{abbrv}  
\bibliography{ref} 
\end{document}